\documentclass[11pt,leqno]{amsart}
\usepackage{amssymb}
\usepackage{mathrsfs}
\usepackage[all]{xy}
\usepackage{color}
\usepackage{soul}
\usepackage{hyperref}

\newcommand\N{{\mathbb N}}
\newcommand\R{{\mathbb R}}

\def\BB{{\mathcal B}}
\def\CC{{\mathcal C}}
\def\DD{{\mathcal D}}
\def\EE{{\mathcal E}}
\def\FF{{\mathcal F}}
\def\GG{{\mathcal G}}

\def\KK{{\mathcal K}}
\def\LL{{\mathcal L}}
\def\MM{{\mathcal M}}

\def\OO{{\mathcal O}}

\def\SS{{\mathcal S}}
\def\TT{{\mathcal T}}
\def\UU{{\mathcal U}}
\def\VV{{\mathcal V}}
\def\WW{{\mathcal W}}
\def\VV{{\mathcal V}}

\def\ZZ{{\mathcal Z}}

\def\eps{{\varepsilon}}

\newcommand{\wto}{\rightharpoonup}

\def\SN{\mathfrak{S}_N}
\def\SSS{\mathfrak{S}}

\newtheorem{theo}{Theorem}

\newtheorem{lem}[theo]{Lemma}
\newtheorem{cor}[theo]{Corollary}
\newtheorem{rem}[theo]{Remark}
\newtheorem{rems}[theo]{Remarks}
\newtheorem{defin}[theo]{Definition}
\newtheorem{ex}[theo]{Example}
\newtheorem{notation}[theo]{Notations}

\newcommand{\beqn}{\begin{equation}}
\newcommand{\eeqn}{\end{equation}}
\newcommand{\bear}{\begin{eqnarray}}
\newcommand{\eear}{\end{eqnarray}}
\newcommand{\bean}{\begin{eqnarray*}}
\newcommand{\eean}{\end{eqnarray*}}

\newcommand{\Black}{\color{black}}

\def\signsm{\bigskip \begin{center} {\sc St\'ephane Mischler\par\vspace{3mm}
Universit\'e Paris-Dauphine\par
CEREMADE, UMR CNRS 7534\par
Place du Mar\'echal de Lattre de Tassigny
75775 Paris Cedex 16\par
FRANCE\par\vspace{3mm}
e-mail:} \tt{mischler@ceremade.dauphine.fr} \end{center}}

\def\signcm{\bigskip \begin{center} {\sc 
Cl\'ement Mouhot\par\vspace{3mm}
University of Cambridge\par
DAMTP, Centre for Mathematical Sciences\par
Wilberforce Road\par
Cambridge CB3 0WA\par
ENGLAND\par
\textit{On leave from:}\par
CNRS \& \'Ecole Normale Sup\'erieure\par
DMA, UMR CNRS 8553\par
45, rue d'Ulm\par
F 75230 Paris cedex 05\par
FRANCE\par\vspace{3mm}
e-mail:} \tt{Clement.Mouhot@ens.fr} \end{center}}

\begin{document}

\title[On chaos propagation of Boltzmann collision processes] 
{Quantitative uniform in time chaos propagation for Boltzmann
  collision processes}


\author{S. Mischler}
\author{C. Mouhot}

\maketitle

\begin{center} {\bf Preliminary version of \today}
\end{center}

\begin{abstract}
  This paper is devoted to the study of mean-field limit for systems
  of indistinguables particles undergoing collision processes. As
  formulated by Kac~\cite{Kac1956} this limit is based on the {\em
    chaos propagation}, and we (1) prove and quantify this property
  for Boltzmann collision processes with unbounded collision rates
  (hard spheres or long-range interactions), (2) prove and quantify
  this property \emph{uniformly in time}. This yields the first chaos
  propagation result for the spatially homogeneous Boltzmann equation
  for true (without cut-off) Maxwell molecules whose ``Master
  equation'' shares similarities with the one of a L\'evy process and
  the first {\em quantitative} chaos propagation result for the
  spatially homogeneous Boltzmann equation for hard spheres
  (improvement of the 
  convergence result of Sznitman~\cite{S1}). Moreover our chaos
  propagation results are the first uniform in time ones for Boltzmann
  collision processes (to our knowledge), which partly answers the
  important question raised by Kac of relating the long-time behavior
  of a particle system with the one of its mean-field limit, and we
  provide as a surprising application a new proof of the well-known
  result of gaussian limit of rescaled marginals of uniform measure on
  the $N$-dimensional sphere as $N$ goes to infinity (more
  applications will be provided in a forthcoming work). Our results
  are based on a new method which reduces the question of chaos
  propagation to the one of proving a purely functional estimate on
  some generator operators ({\em consistency estimate}) together with
  fine stability estimates on the flow of the limiting non-linear
  equation ({\em stability
    estimates}). 
\end{abstract}

\vspace{0.3cm}

processes. \smallskip


\medskip

\textbf{Acknowledgments}: The authors would like to thank the
mathematics departement of Chalmers University for the invitation in
november 2008, where the abstract method was devised and the related
joint work~\cite{MMW} with Bernt Wennberg was initiated. The second
author would like to thank Cambridge University who provided repeated
hospitality in 2009 thanks to the Award No. KUK-I1-007-43, funded by
the King Abdullah University of Science and Technology (KAUST). They
also thank F. Bolley, J. A. Ca\~nizo, N. Fournier, A. Guillin,
P.-L. Lions, J. Rousseau and C. Villani for fruitful
discussions. Finally they wish to mention the inspirative courses of
P.-L. Lions at Coll\`ege de France on ``Mean-Field Games'' in
2007-2008 and 2008-2009.

\tableofcontents



\section{Introduction and main results}
\label{sec:intro}
\setcounter{equation}{0}
\setcounter{theo}{0}


\subsection{The Boltzmann equation}
\label{sec:introEB}

The Boltzmann equation (Cf. \cite{Ce88} and \cite{CIP}) describes the
behavior of a dilute gas when the only interactions taken into account
are binary collisions. It writes
\begin{equation}
  \label{eq:Boltzmann-complete}
  \frac{\partial f}{\partial t} + v \cdot \nabla_x f = Q(f,f)
\end{equation}
where $Q=Q(f,f)$ is the bilinear {\em Boltzmann collision operator}
acting only on the velocity variable.

In the case when the distribution function is assumed to be
independent on the position $x$, we obtain the so-called {\it
  spatially homogeneous Boltzmann equation}, which reads
 \begin{equation}\label{el}
   \frac{\partial f}{\partial t}(t,v) 
   = Q(f,f)(t,v), \qquad  v \in \R^d, \quad t \geq 0,
 \end{equation}
where $d \ge 2$ is the dimension. 

Let us now focus on the collision operator $Q$. It is defined by the
bilinear symmetrized form
 \begin{equation}\label{eq:collop}
 Q(g,f)(v) = \frac12\,\int _{\R^d \times \mathbb{S}^{d-1}} B(|v-v_*|, \cos \theta)
         \left(g'_* f' + g' f_* '- g_* f - g f_* \right) \, dv_* \, d\sigma,
 \end{equation}
where we have used the shorthands $f=f(v)$, $f'=f(v')$, $g_*=g(v_*)$ and
$g'_*=g(v'_*)$. Moreover, $v'$ and $v'_*$ are parametrized by
 \begin{equation}\label{eq:rel:vit}
   v' = \frac{v+v_*}2 + \frac{|v-v_*|}2 \, \sigma, \qquad 
   v'_* = \frac{v+v_*}2 - \frac{|v-v_*|}2 \, \sigma, \qquad 
   \sigma \in \mathbb{S}^{d-1}. 
 \end{equation}
 Finally, $\theta\in [0,\pi]$ is the deviation angle between $v'-v'_*$
 and $v-v_*$ defined by $\cos \theta = \sigma \cdot \hat u$, $u = v-v_*$, $\hat u = u/|u|$, 
  and $B$ is the Boltzmann collision
 kernel determined by physics (related to the cross-section
 $\Sigma(v-v_*,\sigma)$ by the formula $B=|v-v_*| \, \Sigma$).  


Boltzmann's collision operator has the fundamental properties of
conserving mass, momentum and energy
  \begin{equation}
  \int_{\R^d}Q(f,f) \, \phi(v)\,dv = 0, \qquad
  \phi(v)=1,v,|v|^2, \label{CON}
  \end{equation}
  and satisfying the so-called Boltzmann's $H$ theorem

  We shall consider collision kernels $B=\Gamma(|v-v_{*}|) \, b(\cos
  \theta)$ (with $\Gamma,b$ nonnegative functions).
 Typical physical interesting kernels are in dimension $3$:
 \begin{itemize}
 \item {\bf (HS)}: the hard spheres collision kernel $B(|v-v_*|, \cos
   \theta)= \mbox{{\rm cst}} \, |v-v_*|$;
\item collision kernels deriving from interaction potentials $V(r) =
  \mbox{{\rm cst}} \, r^{-s}$, $s >2$: $\Gamma(z)=|z|^\gamma$ with 
  $\gamma = (s-4)/s$, and $b$ is $L^1$ apart from
  $\theta \sim 0$, where
  \[ b(\cos \theta) \sim_{\theta \sim 0} C_b \, \theta^{-(d-1)-\nu} \mbox{
    with } \nu = 2/s \mbox{ (see~\cite{Ce88}) including in particular:} \]
  {\bf (tMM)}: the true Maxwell molecules collision kernel when
  $\gamma=0$
  and $\nu = 1/2$;
\item {\bf (GMM)}: Grad's cutoff Maxwell molecules when $B=1$.
 \end{itemize}
 

\subsection{Deriving the Boltzmann equation from many-particle
   systems}
\label{sec:questionderivation}

The question of deriving the Boltzmann equation from particles systems
(interacting {\it via} Newton's laws) is a famous problem, related to
the so-called $6$-th Hilbert problem mentionned by Hilbert at the
International Congress of Mathematics at Paris in 1900.

At least at the formal level, the correct limiting procedure has been
identified by Grad~\cite{Grad1949} in the late fourties (see
also~\cite{Cerc1972} for mathematical formulation of the open
question): it is now called the {\em Boltzmann-Grad} or {\em low
  density} limit. However the original question of Hilbert remains
largely open, in spite of a striking breakthrough due to
Lanford~\cite{Lanford}, who proved the limit for short times. The
tremendous difficulty underlying this limit is the {\em
  irreversibility} of the Boltzmann equation, whereas the particle
system interacting {\it via} Newton's laws is a reversible Hamiltonian
system.

In 1954-1955, Kac~\cite{Kac1956} proposed a simplified problem in
order to make mathematical progress on the question: start from the
Markov process corresponding to collisions only, and try to prove the
limit towards the {\em spatially homogeneous} Boltzmann
equation. Going back to the idea of Boltzmann of ``stosszahlansatz''
(molecular chaos), he formulated the by now standard notion of {\em
  chaos propagation}.

Let us first define the key notion of {\em chaoticity} for a sequence
$(f^N)_{N \ge 1}$ of probabilities on $E^N$, where $E$ is some given
Polish space (and we will take $E = \R^d$ in the applications):
roughly speaking it means that
\[ 
f^N \sim f^{\otimes N}\quad\hbox{when}\quad N \to \infty
\]
for some given one-particle probability $f$ on $E$. It was clear since
Boltzmann that in the case when the joint probability density $f^N$ of
the $N$-particles system is tensorized into $N$ copies $f^{\otimes N}$
of a $1$-particle probability density, the latter would satisfy the
limiting Boltzmann equation. Then Kac made the key remark that
although in general coupling between a finite number of particles
prevents any possibility of propagation of the ``tensorization''
property, the weaker property of chaoticity can be propagated
(hopefully!) in the correct scaling limit. The application example of
\cite{Kac1956} was a simplified one-dimensional collision model
inspired from the spatially homogeneous Boltzmann equation.

The framework set by Kac is our starting point. Let us emphasize that
the limit performed in this setting is different from the
Boltzmann-Grad limit. It is in fact a {\em mean-field limit}. This
limiting procedure is most well-known for deriving Vlasov-like
equations. In a companion paper~\cite{MMW} we shall develop
systematically our new functional approach for Vlasov equations,
McKean-Vlasov equations, and granular gases Boltzmann collision
models.

\subsection{Goals, existing results and method}
\label{sec:review}

Our goal in this paper is to prove (and set up a general robust method
for proving) chaos propagation with {\em quantitative rate} in terms of
the number of particles $N$ and of the final time of observation $T$.
Let us explain briefly what it means. The original formulation of Kac
\cite{Kac1956} of chaoticity is: a sequence $f^N \in P_{sym}(E^N)$ of
symmetric probabilities on $E^N$ is $f$-chaotic, for a given
probability $f \in P(E)$, if for any $\ell \in \N^*$ and any $\varphi
\in C_b(E)^{\otimes \ell}$ there holds
\[
\lim_{N \to \infty} \left\langle f^N, \varphi \otimes {\bf 1}^{N-\ell}
\right\rangle = \left\langle f^{\otimes \ell}, \varphi \right\rangle
\]
which amounts to the weak convergence of any marginals 
(see also \cite{CCLLV} for another stronger notion of ``entropic''
  chaoticity).  Here we will deal with {\em quantified chaoticity},
in the sense that we measure precisely the rate of convergence in the
above limit. Namely, we say that $f^N$ is $f$-chaotic with rate
$\eps(N)$, where $\eps(N) \to 0$ when $N \to\infty$ (typically
$\eps(N) = N^{-r}$, $r >0$), if for some normed space of {\em smooth
  functions} $\FF \subset C_b(E)$ (to be precised) and for any $\ell
\in \N^*$ there exists $K_\ell \in (0,\infty)$ such that for any
$\varphi \in \FF^{\otimes\ell}$, $\| \varphi \|_\FF \le 1$, there holds
\begin{equation}
  \label{eq:chaos}
  \big| \left\langle \Pi_\ell \left[   f^N  \right] - f ^{\otimes \ell}, 
  \varphi  \right\rangle \big| 
  \le K_\ell \, \eps(N),
\end{equation}
where $\Pi_\ell \left[ f^N \right] $ stands for the $\ell$-th marginal
of $f^N$.

Now, considering a sequence of densities of a $N$-particles system
$f^N \in C([0,\infty);P_{sym}(E^N))$ and a $1$-particle density of the
expected mean field limit $f \in C([0,\infty);P(E))$, we say that
there is {\em propagation of chaos} on some time interval $[0,T]$ if
the $f_0$-chaoticity of the initial family $f^N_0$ implies the
$f_t$-chaoticity of the family $f_t ^N$ for any time $t \in [0,T]$. 

Moreover one can roughly classify the different questions around chaos
propagation into the following layers (in parenthesis, the
corresponding probabilistic interpretation for the empirical measure,
see below):
\begin{enumerate}
\item Proof of propagation in time of the convergence in the
  chaoticity definition (propagation of the law of large numbers along
  time).
\item Same result with a rate $\eps(N)$ as above (estimates of the
  rates in the law of large numbers, estimates on the size of
  fluctuations around the deterministic limit along time).
\item Proof of propagation in time of the convergence to a ``universal
  behavior'' around the deterministic limit (central limit theorem).
  See for instance~\cite{McKean-TCL,Rezakhanlou-98} for related
  results.
\item Proof of propagation in time of bounds of exponential type on
  the ``rare'' events far from chaoticity (large deviation estimates).
\item Estimations (2)--(4) can be made uniformly on intervals $[0,T)$
  with $T$ finite or $T=+\infty$.
\end{enumerate}

For Boltzmann collision processes, Kac \cite{Kac1956}-\cite{Kac1957}
has proved the point (1) in the case of his baby one-dimensional
model. The key point in his analysis is a clever combinatorial use of
a semi-explicit form of the solution (Wild sums). It was generalized
by McKean \cite{McKean1967} to the Boltzmann collision operator but
only for ``Maxwell molecules with cutoff'', {\it i.e.}, roughly when
the collision kernel $B$ above is constant. In this case the
combinatorial argument of Kac can be extended. Kac raised in
\cite{Kac1956} the question of proving chaos propagation in the case
of hard spheres and more generally unbounded collision kernels,
although his method seemed impossible to extend (no semi-explicit
combinatorial formula of the solution exists in this case).

In the seventies, Gr\"unbaum~\cite{Grunbaum} then proposed in a very
compact and abstract paper another method for dealing with hard spheres,
based on the Trotter-Kato formula for semigroups and a clever functional
framework (partially remindful of the tools used for mean-field limit for
McKean-Vlasov equations). Unfortunately this paper was incomplete for two
reasons: (1) It was based on two ``unproved assumptions on the Boltzmann
flow'' (page 328): (a) existence and uniqueness for measure solutions and
(b) a smoothness assumption. Assumption (a) was indeed recently proved in
\cite{Fo-Mo} using Wasserstein metrics techniques and in \cite{EM} adapting 
the classical DiBlasio trick \cite{DiB74}, but concerning
assumption (b), although it was inspired by cutoff maxwell molecules (for
which it is true), it fails for hard spheres (cf. the counterexample built
by Lu and Wennberg in \cite{LuWe}) and is somehow ``too rough'' in this
case.  (2) A key part in the
proof in this paper is the expansion of the ``$H_f$'' function, which is an
a clever idea of Gr\"unbaum (and the starting point for our idea of
developing an abstract differential calculus in order to control
fluctuations) --- however it is again too rough and is adapted for cutoff
Maxwell molecules but not hard spheres.

A completely different approach was undertaken by Sznitman in the eighties
\cite{S6} (see also Tanaka \cite{T2}). Starting from the observation that
Gr\"unbaum's proof was incomplete, he gave a full proof of chaos
propagation for hard spheres. His work was based on: (1) a new uniqueness
result for measures for the hard spheres Boltzmann equation (based on a
probabilistic reasoning on an enlarged space of ``trajectories''); (2) an
idea already present in Gr\"unbaum's approach: reduce by a combinatorial
argument on symmetric probabilities the question of chaos propagation to a
law of large numbers on measures; (3) a new compactness result at the level
of the empirical measures; (4) the identification of the limit by an
``abstract test function'' construction showing that the (infinite
particle) system has trajectories included in the chaotic ones. Hence the
method of Sznitman proves convergence but does not provide any rate for
chaoticity.  Let us also emphasize that Graham and M\'el\'eard in \cite{GM}
have obtained a rate of convergence (of order $1/ \sqrt{N}$) on any bounded
finite interval of the $N$-particles system to the deterministic Boltzmann
dynamic in the case of Maxwell molecules under Grad's cut-off hypothesis,
and that Fournier and M\'el\'eard in \cite{FM7,FM10} have obtained the
convergence of the Monte-Carlo approximation (with numerical cutoff) of the
Boltzmann equation for true Maxwell molecules with a rate of convergence
(depending on the numerical cutoff and on the number $N$ of particles).

Our starting point was Gr\"unbaum's paper \cite{Grunbaum}. Our original
goal was to construct a general and robust method able to deal with mixture
of jump and diffusion processes, as it occurs for granular gases (see for
this point the companion paper \cite{MMW}).  It turns out that it lead us
to develop a new theory, inspiring from more recent tools such as the
course of Lions on ``Mean-field games'' at Coll\`ege de France, and the
master courses of M\'el\'eard \cite{Meleard1996} and Villani
\cite{VillaniMF} on mean-field limits. One of the byproduct of our paper is
that we make fully rigorous the original intuition of Gr\"unbaum in order
to prove chaos propagation for the Boltzmann velocities jump process
associated to hard spheres contact interactions.


As Gr\"unbaum \cite{Grunbaum} we shall use a duality argument.  We
introduce $S^N_t$ the semigroup associated to the flow of the
$N$-particle system and $T^N_t$ its {\em ``dual'' semigroup}.  We also
introduce $S^{N \! L}_t$ the (nonlinear) semigroup associated to the
meanfield dynamic (the exponent ``NL'' recalling that the limit
semigroup is nonlinear in the most physics interesting cases) as well
as $T^\infty_t$ the associated (linear) {\em ``pushforward''
  semigroup} (see below for the definition). Then we will prove the
above kind of convergence on the linear semigroups $T^N_t$ and
$T^\infty_t$.
\Black

The first step consists in defining a common functional framework in which
the $N$-particles dynamic and the limit dynamic make sense so that we can
compare them. Hence we work at the level of the ``full'' limit space
$P(P(E))$ (see below).
Then we shall identify the regularity required in order to prove the
``consistency estimate'' between the generators $G^N$ and $G^\infty$ of the
dual semigroups $T^N _t$ and $T^\infty _t$, and then prove a corresponding
``stability estimate'' at the level of the limiting semigroup $S^{N \!
  L}_t$. The latter crucial step shall lead us to introduce an abstract
differential calculus for functions acting on measures endowed with various
metrics.

In terms of existing open questions, this paper solves two related
problems. First it proves quantitative rates for chaos propagation for
hard spheres and for (non cutoff) Maxwell molecules. These two results
can be seen as two advances into proving chaos propagation for
collision processes with unbounded kernels (in the two physically
relevant ``orthogonal'' directions: unboundedness is either due to
growth at large velocities, or to grazing collisions). Second we
provide the first \emph{uniform in time} chaos propagation results
(moreover quantitative), which answers partly the question raised by
Kac of relating the long-time behavior of the $N$-particle system with
the one of its mean-field limit. 

Finally the general method that we provide for solving these problems
is, we hope, interesting by itself for several reasons: (1) it is
fully quantitative, (2) it is highly flexible in terms of the
functional spaces used in the proof, (3) it requires a minimal amount
of informations on the $N$-particles systems but more stability
information on the limiting PDE  (we intentionally presented the
assumption as for the proof of the convergence of a numerical scheme,
which was our ``methodological model''), (4) the ``differential
stability'' conditions that are required on the limiting PDE seem (to
our knowledge) new, at least at the level of Boltzmann or more
generally transport equations.

\subsection{Main results}
\label{sec:intromainresults}

Without waiting for the full abstract framework, let us give a slightly
fuzzy version of our main results, gathered in a single theorem. The full
definitions of all the objects considered shall be given in the forthcoming
sections, together with fully rigorous statements (see Theorems
\ref{theo:abstract}, \ref{theo:tMM} and \ref{theo:HS}).

\begin{theo}
  \label{theo:intro}
  Let $S^{N \! L} _t$ denotes the semigroup of the spatially homogeneous
  Boltzmann equation (acting on $P(\R^d)$) and $S^N _t$, $N \ge 1$
  denotes the semigroup of the $N$-particle system satisfying a
  collision Markov process. Then for any (one particle) initial datum
  $f_0 \in P(\R^d)$ with compact support and the corresponding
  tensorized $N$ particle initial data $f_0 ^{\otimes N}$, we have

 \begin{itemize}
 \item[(i)] In the case where $S^{N \! L} _t$ and $S^N_t$ correspond to
   the {\bf hard spheres collision kernels (HS)}, we have for any
   $\ell \in \N^*$ and any $N \ge 2 \ell$:
    \[ 
    \sup_{t \in (0,\infty)} \,\, \sup_{\varphi \in
      W^{1,\infty}(\R^d)^{\otimes\ell}, \, \|Ê\varphi
      \|_{W^{1,\infty}(\R^d)^{\otimes\ell}} \le 1}
    \left\langle \Pi_\ell \left[ S^N_t\left(f_0^{\otimes N}\right)
      \right] - S^{N \! L}_t(f_0)^{\otimes \ell}, \varphi
    \right\rangle \le \ell \, \eps(N) 
    \] 
    
  \item[(ii)] In the case where $S^\infty_t$ and $S^N_t$ correspond to the
    {\bf true (or cutoff) Maxwell molecules (tMM)-(GMM)}, we have for
    any $\ell \in \N^*$, any $N \ge 2 \ell$ and any $\delta \in
    (0,2/d)$:
   \[ 
   \sup_{t \in (0,\infty)} \,\, \sup_{\varphi \in \FF^{\otimes\ell},
     \, \|\varphi \|_{\FF^{\otimes\ell}} \le 1} \left\langle
     \Pi_\ell \left[ S^N_t\left(f_0^{\otimes N}\right) \right] - S^{N
       \! L} _t(f_0)^{\otimes \ell}, \varphi \right\rangle \le \ell^2
   \, {C_\delta \over N^{{2 \over d} - \delta}} 
    \] 
    for some constant $C_\eta \in (0,\infty)$ which may blow up when
    $\eta \to 0$ and where
   $$
   \FF := \left\{ \varphi : \R^d \to \R; \,\, \| \varphi \|_\FF :=
     \int_{\R^d} (1 + |\xi|^4) \, |\hat \varphi (\xi)|\, d\xi <
     \infty \right\},
   $$
   $\hat \varphi$ standing for the Fourier transform of $\varphi$. 
   
 \item[(iii)] In the case where $S^\infty_t$ and $S^N_t$ correspond to the
   {\bf cutoff Maxwell molecules (GMM)}, we have for any $\ell \in
   \N^*$, any $N \ge 2 \ell$, any $s \in (d/2,d/2+1)$ and any $T \in
   (0,\infty)$:
   \[ 
   \sup_{t \in (0,T)} \,\, \sup_{\varphi \in
     H^{s}(\R^d)^{\otimes\ell}, \, \| \varphi
     \|_{H^{s}(\R^d)^{\otimes\ell}} \le 1} \left\langle \Pi_\ell
     \left[ S^N_t\left(f_0^{\otimes N}\right) \right] - S^{N \! L}
     _t(f_0)^{\otimes \ell}, \varphi \right\rangle \le \ell^2 \,
   {C_{s,T} \over N^{{1 \over 2}}} 
    \] 
    for some constant $C_{s,T} \in (0,\infty)$ which may blow up when
    $s \to d/2$ or $s \to d/2+1$ or when $T\to\infty$.
 \end{itemize}
\end{theo}


\begin{rems}
\begin{itemize}
 
\item To be more precise the constants depend on the the initial datum
  $f_0$ through polynomial moments bounds for \textbf{(GMM)} and
  \textbf{(tMM)}, and exponential moment bound for
  \textbf{(HS)}. However some one needs also a bound on the support of
  the energy of the $N$-particle empirical measure. When the
  $N$-particle initial data are simply tensor products of the
  $1$-particle initial datum, the simplest sufficient condition is the
  compact support of $f_0$. It could be relaxed at the price of an
  additional error term. This restriction can also be relaxed by
  considering $N$-particle initial datum conditioned to the sphere
  $S^{Nd-1}(\sqrt{N})$ of constant energy as in \cite{Kac1956}: these
  extensions are considered in
  Section~\ref{sec:extensions}. 

\item Note that the two first estimates (i) and (ii) are {\em global
    in time}. This is an important qualitative improvement over
  previous chaos propagation results for Boltzmann collision
  processes.

\item In the third estimate (iii), the rate of convergence
  $\OO(1/N^{1/2})$ is optimal in the case of Maxwellian kernels, as
  predicted by the law of large numbers.
  \end{itemize}
\end{rems}

\subsection{Some open questions and extensions}
\label{sec:plan}

\begin{itemize}
\item What about larger classes of initial data, say with only and as few
  as possible polynomial moments?
\item   What about optimizing the rate?   Is
the rate $N^{-1/2}$ always optimal and how to obtain in general when no
Hilbert structure seems available for the estimates (such as for the hard
spheres case {\bf (HS)})? 
\item What about more general Boltzmann
models? A work is in progress for applying the method to inelastic hard
spheres and inelastic Maxwell molecules Boltzmann equation with thermal or
stochastic baths. Another issue is the true (without cut-off) Boltzmann
equation for hard or soft potential: only assumption {\bf (A4)} (see below)
remains to be proved in that case. The question of proving uniform in time
chaos propagation also remains open in this case.
\item What about a central limit theorem and/or a large deviation result
  with the help of our setting? A natural guess would be that a
  higher-order ``differential stability'' (see below) is required on the
  limit system, and this point should be clarified. 
\Black
\item Finally let us mention that works are in progress to apply our
  method to Vlasov and McKean-Vlasov equations.
\end{itemize}

\subsection{Plan of the paper}
\label{sec:plan}

In Section~\ref{sec:abstract-setting} we set the abstract functional
framework together with the general assumption and in
Section~\ref{sec:abstract-theo} we state and prove the abstract
Theorem~\ref{theo:abstract}. In Section~\ref{sec:BddBoltzmann} we
apply this method to the (true) Maxwell molecules: we show how to
choose metrics so that the general assumptions can be proved
(Theorem~\ref{theo:tMM}). In Section~\ref{sec:hardspheres} we
apply the method to hard spheres molecules and quantitative chaos
propagation in Theorem~\ref{theo:HS}. Finally in
Section~\ref{sec:extensions} we address several important extensions
and applications: we relaxe the assumption of compactly supported
initial datum by considering $N$-particle initial data conditioned on
the energy sphere as suggested by Kac, and as a consequence of our
uniform in time chaos propagation result we study the chaoticity of
the steady state. Finally we conclude with some remarks and
computations on the BBGKY hierarchy and its link with our work.


\section{The abstract setting}
\label{sec:abstract-setting}
\setcounter{equation}{0}
\setcounter{theo}{0}


In this section we shall state and prove the key abstract result. This will
motivate the introduction of a general functional framework.

\subsection{The general functional framework of the duality approach}
\label{framework}

Let us set the framework. Here is a diagram which sums up the duality
approach (norms and duality brackets shall be precised in
Subsections~\ref{sec:funct-set}):
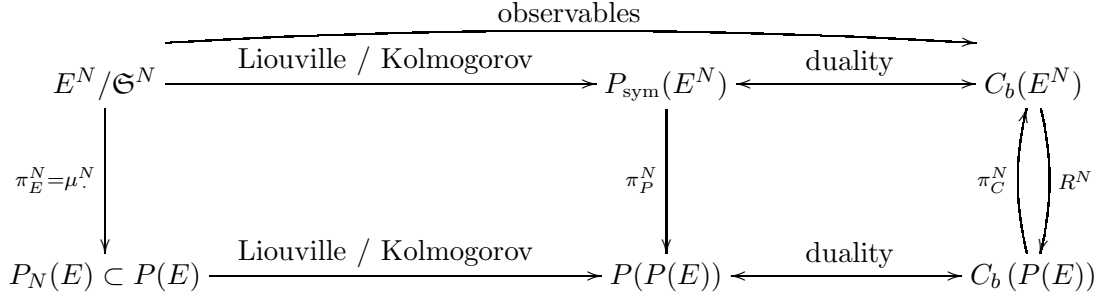
\begin{displaymath}
  \xymatrix{
    E^N / \SSS^N \ar[dd]_{\pi^N_E = \mu^N _{\, \cdot}} 
    \ar[rrrrr]^{\mbox{Liouville / Kolmogorov}} 
    \ar@/^/@<+3ex>[rrrrrrrr]^{\mbox{observables}}
    & & & & & P_{\mbox{{\tiny sym}}}(E^N) \ar[dd]_{\pi^N_P}
    \ar[rrr]^{\mbox{duality}} 
    & & &  C_b(E^N) \ar[lll] \ar@/^/[dd]^{R^N} \\
    & & & & & & & \\
    P_N(E) \subset P(E)  
    \ar[rrrrr]^{\mbox{Liouville / Kolmogorov}} & & & & & 
    P(P(E)) \ar[rrr]^{\mbox{duality}} & & &  
    C_b\left(P(E)\right) \ar[lll] \ar@/^/[uu]^{\pi^N_C} 
  }
\end{displaymath}

In this diagram:
\begin{itemize}
\item $E$ denotes a Polish space. 
\item $\SSS^N$ denotes the $N$-permutation group.
\item $P_{\mbox{{\tiny sym}}}(E^N)$ denotes the set of symmetric
  probabilities on $E^N$: For a given permutation $\sigma \in \SSS^N$, a
  vector $V = (v_1, ..., v_N) \in E^N$, a function $\phi \in C_b(E^N)$ and
  a probability $\rho^N \in P(E^N)$ we successively define $V_\sigma =
  (v_{\sigma(1)}, ..., v_{\sigma(N)}) \in E^N $, $\phi_\sigma \in C_b(E^N)$
  by setting $\phi_\sigma(V) = \phi(V_\sigma)$ and $\rho^N_\sigma \in P(E^N)$
  by setting $\langle \rho^N_\sigma, \phi \rangle = \langle \rho^N, \phi_\sigma
  \rangle$.  We say that a probability $\rho^N$ on $E^N$ is symmetric or
  invariant under permutations if $\rho^N_\sigma = \rho^N$ for any permutation
  $\sigma \in \SSS^N$.
\item For any $V \in E^N$ the probability measure $\mu^N _V$ denotes the
  {\em empirical measure}:
  \[
  \mu^N _V = \frac1N \, \sum_{i=1} ^N \delta_{v_i},
  \quad V=(v_1,\dots, v_N)
  \]
  where $\delta_{v_i}$ denotes the Dirac mass on $E$ at point $v_i$.
\item $P_N(E)$ denotes the subset $\{\mu^N_V, \, V \in E^N \}$ of
  $P(E)$.
\item $P(P(E))$ denotes the set of probabilities on the polish space
  $P(E)$.
\item $C_b\left(P(E)\right)$ denotes the space of continuous and bounded
  functions on $P(E)$, the latter space being endowed with the weak or
  strong topologies (see Subsection~\ref{sec:funct-set}).
\item The arrow pointing from $E^N / \SSS^N$ to $P_N(E)$ denotes the
  map $\pi^N_E$ defined by 
  \[
  \forall \, V \in E^N / \SSS^N, \quad \pi^N_E(V)  := \mu^N _V.
  \]
\item The arrow pointing from $C_b(P(E))$ to $C_b(E^N)$ denotes the
  following map $\pi^N_C$
  \[ 
  \forall \, \Phi \in C_b\left(P(E)\right), \ \forall \, V \in E^N, \
  \left(\pi^N_C\Phi\right)(V) := \Phi\left( \mu^N _V \right).
  \]
\item The counter arrow pointing from $C_b(E^N)$ to $C_b(P(E))$
  denotes the transformation $R^N$ defined by:
  \[ 
  \forall \, \phi \in C_b(E^N), \ \forall \, \rho \in P(E), \quad
  R^N[\phi](\rho) := \left\langle \rho^{\otimes N}, \phi
  \right\rangle.
  \]
  In the sequel we shall sometimes use the shorthand notation
  $R^\ell_\phi$ instead of $R^\ell[\phi]$ for any $\ell \in \N^*$ and
  $\phi \in C_b(E^\ell)$.
\item The arrow pointing from $P_{\mbox{{\tiny sym}}}(E^N)$ to
  $P(P(E))$ denotes the following transformation: consider a symmetric
  probability $\rho^N \in P_{\mbox{{\tiny sym}}}(E^N)$ on $E^N$, we
  define the probability on probability $\pi^N_P \rho^N \in P(P(E))$ by
  setting
  \[ 
  \forall \, \varphi \in \Phi \in C_b(P(E)), \quad 
  \left\langle \pi^N_P \rho^N,\Phi \right\rangle = 
  \left\langle \rho^N,\pi^N_C \Phi \right\rangle =
  \left\langle \rho^N, \Phi\left( \mu^N _V \right) \right\rangle,
  \] 
  where the first bracket means $\langle \cdot , \cdot
  \rangle_{P(P(E)),C_b(P(E))}$ and the second bracket means $\langle \cdot ,
  \cdot \rangle_{P(E^N),C_b(E^N)}$.
\item The arrows pointing from the first column to the second one denote
  the procedure of the writing either of the {\em Liouville} transport
  equation associated with the set of ODEs of a particle system (as for
  mean-field limits for Vlasov equations), or the writing of the {\em
    Kolmogorov} equation associated with a stochastic Markov process of a
  particles system (e.g. jump or diffusion processes).
\item Finally the dual spaces of the spaces of probabilities on the
  phase space can be interpreted as the spaces of observables on the
  original systems. We shall discuss this point later.
\end{itemize}

\begin{rem} Consider a random variable $V$ on $E^N$ with law $\rho^N \in
  P(E^N)$. Then it is often denoted by $\mu^N _V$ the
  random variable on $P(E)$ with law $\pi^N_P \rho^N \in P(P(E))$. Our
  notation is slightly less compact and intuitive, but at the same
  time more accurate.
\end{rem}

\begin{rem}\label{rem:poids}
  Our functional framework shall be applied to weighted probability
  spaces rather than directly in $P(E)$. More precisely, for a given
  weight function $m : E \to \R_+$ we shall use (subsets) of the
  weighted probability space
  \begin{equation}\label{def:Mmrho}
    \left\{ \rho \in P(E); \ 
      M_m(\rho) := \langle \rho, m \rangle < \infty \right\}
  \end{equation}
  as our core functional space. Typical examples are $m(v) := \tilde m
  (\hbox{dist}_E (v,v_0))$, for some fixed $v_0 \in E$ and $\tilde
  m(z) = z^k$ or $\tilde m(z) = e^{a \, z^k}$, $a,k>0$. We shall
  sometimes abuse notation by writing $M_k$ for $M_m$ when $\tilde
  m(z) = z^k$ in the above example.
\end{rem}

\subsection{The evolution semigroups}
\label{sec:semigroups}

Let us introduce the mathematical objects living in these spaces, for any
$N \ge 1$.  \smallskip

\noindent {\em Step~1.} Consider a process $(\VV^N_t)$ on $E^N$ which
describes the trajectories of the particles (Lagrangian viewpoint). The
evolution can correspond to stochastic ODEs (Markov process), or
deterministic ODEs (deterministic Hamiltonian flow). We make the
fundamental assumption that this flow commutes with permutations: for any
$\sigma \in \SSS^N$, the solution at time $t$ starting from
$\left(\VV^N_0\right)_\sigma$ is $\left(\VV^N_t\right)_\sigma$. This
reflects mathematically the fact that particles are indistinguishable.
\smallskip

\noindent {\em Step~2.} One naturally derives from this flow on $E^N$
a corresponding semigroup $S^N _t$ acting on $P_{\mbox{{\tiny
      sym}}}(E^N)$ for the presence density of particles in the phase
space $E^N$. This corresponds to a linear evolution equation
\begin{equation}
  \label{eq:evolN}
  \partial_t f^N = A^N f^N, \qquad f^N \in P_{\mbox{{\tiny sym}}}(E^N), 
\end{equation}
which can be interpreted as the forward Kolmogorov equation on the law
in case of a Markov process at the particle level, or the Liouville
equation on the probability density in case of an Hamiltonian process
at the particle level. As a consequence of the previous assumption that the flow $(\VV^N_t)$
commutes with permutation, we have that $S^N _t$ acts on
$P_{\mbox{{\tiny sym}}}(E^N)$. In other words, if the law $f^N_0$ of
$\VV^N_0$ belongs to $P_{\mbox{{\tiny sym}}}(E^N)$, then for any times
the law $f^N_t$ of $\VV^N_t$ also belongs to $P_{\mbox{{\tiny sym}}}(E^N)$. 

\smallskip

\noindent {\em Step~3.} We then define the {\em dual} semigroup $T^N
_t$ of $S^N _t$ acting on functions $\phi \in C_b(E^N)$ by
\[ 
\forall \, f^N \in P(E^N), \ \phi \in C_b(E^N), \quad \left\langle
  f^N, T^N _t (\phi) \right\rangle := \left\langle S^N_t(f^N), \phi
\right\rangle
\]
and we denote its generator by $G^N$, which corresponds to the following
{\em linear} evolution equation:
\begin{equation}\label{eq:dualN}
  \partial_t \phi = G^N(\phi), \qquad \phi \in C_b(E^N).
\end{equation}
This is the semigroup of the {\em observables} on the evolution system
$(\VV_t ^N)$ on $E^N$.  \smallskip

Let us state precisely our assumptions on the $N$-particles dynamics.  Here
and below, for a given $1$-particle weight function $m: E \to \R_+$,
we define the $N$-particle weight function  
\begin{equation}
\label{eq:defMm}
\forall \, V = (v_1, ... , v_N) \in E^N, \qquad M_m^N(V) := {1 \over
  N} \sum_{i=1}^N m(v_i) = \langle \mu^N_V, m \rangle = M_m (\mu^N_V).
\end{equation}
Again, we shall sometimes abuse notation by writing $M^N_k$ and $M_k$
instead of $M^N_m$ and $M_k$ with $\tilde m(z) = z^k$ in the example
of Remark~\ref{rem:poids}. In particular, when $E=\R^d$ endowed with
the euclidian structure, we have
\begin{equation}
\label{eq:defMmV}
\forall \, V= (v_1, ... , v_N) \in \R^{dN}, \qquad M_k^N(V) := {1 \over
  N} \sum_{i=1}^N |v_i|^k = \langle \mu^N_V, | \cdot |^k \rangle = M_k (\mu^N_V).
\end{equation}

\fbox{
\begin{minipage}{0.9\textwidth}
\begin{itemize}
\item[{\bf (A1)}] {\bf On the $N$-particle system.} Together with the fact
  that $G^N$ and $T^N_t$ are well defined and satisfy the symmetry
  condition introduced in Step 2 above, we assume that the following
  moments conditions hold:
  \begin{itemize}
  \item[(i)] {\it Energy bounds}: There exists a weight function $m_e$
    and a constant $\EE \in (0,\infty)$ such that
    \begin{equation}\label{eq:pNUnifEnergy} 
      \hbox{supp} \, f^N_t \subset \mathbb{E}_N :=
      \left\{V \in E^N; \,\, M^N_{m_e}(V) \le \EE\right\}.  
  \end{equation}
 
\item[(ii)] {\it Integral moment bound}: There exists a weight
  function $m_1$, a time $T \in (0,\infty]$ and a constant
  $C^N_{T,m_1} \in (0,\infty)$, possibly depending on $T$, and $m_1$,
  $m_e$ and $\EE$, but not on the number of particles $N$, such that
    \begin{equation}\label{eq:pNBornes0} 
      \sup_{0 \le t < T} \left\langle f^N
      _t, M^N _{m_1} \right \rangle \le C^N_{T,m_1}.  
  \end{equation}
  
\item[(iii)] {\it Support moment bound at initial time}: There exists
  a weight function $m_3$ and a constant $C^N_{0,m_3} \in
  (0,+\infty)$, possibly depending on the number of particles $N$ and
  on $\EE$, such that
    \begin{equation}\label{eq:pNm30} 
      \hbox{supp}
      \, f^N_0 \subset \left\{V \in E^N; M^N_{m_3} (V) \le C^N_{0,m_3} \right\}.
    \end{equation}
  \end{itemize}
 \end{itemize}
\end{minipage}
 }
\medskip

 \begin{rem}
   Note that these assumptions on the $N$-particle system are very
   weak, and do not require any precise knowledge of the $N$-particle
   dynamics.
 \end{rem}

\smallskip 
\noindent {\em Step~4.} Consider a (possibly nonlinear) semigroup $S^{N \!
  L}_t$ acting on $P(E)$ associated with the limit (possibly nonlinear)
kinetic equation: for any $\rho \in P(E)$, $S^{N \! L}_t(\rho) := f_t$ where $f_t
\in C(\R_+, P(E))$ is the solution to
\begin{equation}\label{eq:limit}
\partial_t f_t  = Q(f_t), \quad f_0 = \rho.
\end{equation}
\smallskip

\noindent {\em Step~5.}
Then we consider its {\em pushforward} semigroup $T^\infty _t$ acting
on $C_b(P(E))$ defined by:
\[ 
\forall \, \rho \in P(E), \ \Phi \in C_b(P(E)), \quad T^\infty _t
[\Phi](\rho) := \Phi\left( S^{N\! L}_t(\rho)\right).
\]
Note carefully that $T_t ^\infty$ is always {\em linear} as a function
of $\Phi$ (although of course be careful that $T_t ^\infty[\Phi](\rho)$
is not linear as a function of $\rho$). We denote its generator by
$G^\infty$, which corresponds to the following {\em linear} evolution
equation on $C_b(P(E))$:
\begin{equation}\label{eq:dual-limit}
  \partial_t \Phi = G^\infty(\Phi).
\end{equation}

\begin{rem}
  The semigroup $T_t ^\infty$ can be interpreted physically as the
  semigroup of the evolution of {\em observables} of the nonlinear
  equation \eqref{eq:limit}:\\
  Given a nonlinear ODE $V'=F(V)$ on $\R^d$, one can define (at least
  formally) the {\em linear} Liouville transport PDE
  \[ 
  \partial_t \rho + \nabla_v \cdot (F \, \rho) =0,
  \]
  where $\rho=\rho_t(v)$ is a time-dependent probability density. When
  the trajectories $(V_t(v))$ of the ODE are properly defined, the
  solution of the associated transport equation is given by $\rho_t(v)
  = V_{-t} ^* (\rho_0)$, that is the pullback of the initial measure
  $\rho_0$ by $V_{-t}$ (for smooth functions, this amounts to
  $\rho_t(v) = \rho_0(V_{-t}(v))$). Now, instead of the Liouville
  viewpoint, one can adopt the viewpoint of {\em observables}, that is
  functions depending on the position of the system in the phase space
  (e.g. energy, momentum, etc.) For some observable function $\Phi_0$
  defined on $\R^d$, the evolution of the value of this observable
  along the trajectory is given by $\phi_t(v)=\phi_0(V_t(v))$. In
  other words we have $\phi_t = V_t ^* \phi_0 = \phi_0 \circ
  V_t$. Then $\phi_t$ is solution to the following {\em dual} linear PDE
  \[
  \partial_t \phi - F \cdot \nabla_v \phi = 0,
  \]
  and it satisfies
  \[ 
  \langle \phi_t, \rho_0 \rangle = \langle V_t ^* \phi_0, \rho_0 \rangle =
  \langle \phi_0, V_{-t} ^* \rho_0 \rangle = \langle \phi_0, \rho_t \rangle.
  \]

  Now let us consider a nonlinear evolution system $V'=Q(V)$ in an
  asbtract space $H$. Then (keeping in mind the analogy with ODE/PDE
  above) we see that one can formally naturally define two linear
  evolution systems on the larger functional spaces $P(H)$ and
  $C_b(H)$: first the transport equation (trajectories level)
\[ 
\partial_t \pi + \nabla \cdot (Q(v) \, \pi) = 0, \qquad \pi \in P(H)
\]
(where the second term of this equation has to be properly defined)
and second the {\em dual} equation for the evolution of observables
\[
\partial_t \Phi - Q(v) \cdot \nabla \Phi = 0, \qquad \Phi \in C_b(H).
\]

Taking $H=P(E)$, this provides an intuition for our functional
construction, and also for the formula of the generator $G^\infty$ below
(compare the previous equations with formula \eqref{eq:formulaGinfty}) and
the need for developping a differential calculus in $P(E)$. Be careful that
when $H=P(E)$, the word ``trajectories'' refers to trajectories {\em in the
  space of probabilites} $P(E)$ ({\em i.e.}, solutions to the nonlinear
equation \eqref{eq:limit}), and not trajectories of a particle in $E$.

Another important point to notice is that for a {\em dissipative equation}
at the level of $H$, one cannot define reverse ``characteristics'' (the
flow in $H$ is not defined backwards), and therefore for the nonlinear
Boltzmann equation, only the equation for the observables can be defined in
terms of trajectories.
\end{rem}

Summing up we obtain the following picture for the semigroups:
\begin{displaymath}
  \xymatrix{
    P^N_t \ \mbox{on} \ E^N / \SSS^N \ar[dd]_{\mu^N _V} 
    \ar[rrr]^{\mbox{observables}} 
    & & & T^N _t \ \mbox{on} \ C_b(E^N) \ar@/^/[dd]^{R^N} \\
    & & &  \\
    P_N(E) \subset P(E)  
    \ar[rrr]^{\mbox{observables}} 
    & & & \boxed{T^\infty _t \ \mbox{on} \ C\left(P(E)\right)} 
    \ar@/^/[uu]^{\pi^N _C}
    \\
    & & & \\
    & & & S^{N\! L}_t\ \mbox{on} \ P(E) \ar[uu]_{\mbox{observables}} 
  }
\end{displaymath}

Hence we see that the key point of our construction is that, through the
evolution of {\em observables} one can ``interface'' the two evolution
systems (the nonlinear limit equation and the $N$-particles system) {\em
  via} the applications $\pi^N _C$ and $R^N$. From now on we shall denote
$\pi^N = \pi^N _C$.


\subsection{The metric issue}
\label{sec:funct-set}

$P(E)$ is our fundamental space, where we shall compare (through their
observables) the marginals of the $N$-particle density $f^N_t$ and the
marginals of the chaotic $\infty$-particle dynamic $f_t^{\otimes
  \infty}$. Let us make precise the topological and metric structure
used on $P(E)$. At the topological level there are two canonical
choices (which determine two different sets $C(P(E))$): (1) the strong
topology (associated to the total variation norm that we denote by
$\|\cdot \|_{M^1}$) and (2) the weak topology ({\it i.e.}, the trace
on $P(E)$ of the weak topology $M^1(E)$, the space of Radon measures
on $E$ with finite mass, induced by $C_b(E)$).



The set $C_b(P(E))$ depends on the choice of the topology on
$P(E)$. In the sequel, we will denote $C_b(P(E),w)$ the space of
continuous and bounded functions on $P(E)$ {\em endowed with the weak
  topology}, and $C_b(P(E),TV)$ the similar space on $P(E)$ {\em
  endowed with the total variation norm}. It is clear that
$C_b(P(E),w) \subset C_b(P(E),TV)$. 

The supremum norm $\| \Phi \|_{L^\infty(P(E))}$ {\em does not} depend on
the choice of topology on $P(E)$, and induces a Banach topology on the
space $C_b(P(E))$. The transformations $\pi^N$ and $R^N$ satisfy:
\begin{equation}\label{eq:compat:infty}
  \left\| \pi^N \Phi \right\|_{L^\infty(E^N)} \le \| \Phi
  \|_{L^\infty(P(E))} \ \mbox{ and } \ \| R^N[\phi] \|_{L^\infty(P(E))} \le
  \| \phi \|_{L^\infty(E^N)}.
\end{equation}

The transformation $\pi^N$ is well defined from $C_b(P(E),w)$ to
$C_b(E^N)$, but in general, it does not map $C_b(P(E),TV)$ into
$C_b(E^N)$ since $V \in E^N \mapsto \mu^N_V \in (P(E),TV)$ is not
continuous.

In the other way round, the transformation $R^N$ is well defined from
$C_b(E^N)$ to $C_b(P(E),w)$, and therefore also from $C_b(E^N)$ to
$C_b(P(E),TV)$: for any $\phi \in C_b(E^N)$ and for any sequence
$f_k \wto f$ weakly, we have $f_k^{\otimes N} \wto
f^{\otimes N}$ weakly, and then $R^N[\phi](f_k) \to
R^N[\phi](f)$.

\smallskip The different possible metric structures inducing the weak
topology are not seen at the level of $C_b(P(E),w)$. However any H\"older
or $C^k$-like space will strongly depend on this choice, as we shall see.


\begin{defin}\label{defGG} 
  For a given weight functions $m_\GG : E \to \R_+$
 and some constant $\EE \in (0,\infty)$,
  we define the subspaces of probabilities:
 \begin{equation}\label{def:PGG}
  P_{\GG} := \{ f \in P(E); \,\, \langle f, m_\GG \rangle < \infty, \,\, 
  \langle f, m_e \rangle \le \EE \rangle ,
  \end{equation}
   where $m_e$ is introduced in {\bf (A1)}. 
  For a  given constraint function ${\bf m}_\GG : E \to \R^D$ such that the
  components ${\bf m}_\GG$ are controlled by $m_\GG$, we
  also define the corresponding constrained subsets
  $$
  P_{\GG,\hbox{\small\bf r}} := \{ f \in P_\GG; \,\, \langle
  f, {\bf m}_\GG \rangle = {\bf r}\}, \quad {\bf r} \in \R^D,
  $$
  the corresponding bounded subsets for $a > 0$
  $$
  \mathcal{B}P_{\GG,a} := \{ f \in P_GG); \,\, \langle f, m_\GG
  \rangle < a \}, \quad P_{\GG,\hbox{\small\bf r},a} := \{ f \in
  \mathcal{B}P_{\GG,a}; \,\, \langle f, {\bf m}_\GG \rangle = {\bf r}\},
  $$
  and the corresponding vectorial space of ``increments''
  $$
  \mathcal{I} P_\GG := \left\{f_2 - f_1; \,\, \exists \, {\bf r }Ê\in \R^D \,\,\hbox{s.t.}\,\, f_1,f_2 \in P_{\GG,{\bf r}} \right\}
  $$
  and $\mathcal{I}P_{\GG,a}$ as expected. Now, we shall encounter two
  situations:
  \begin{itemize}
  \item Either $\mbox{dist}_\GG$ denotes a distance defined on the
    whole space $P_{\GG}(E)$, and thus on $P_{\GG,\hbox{\small\bf r}}$
    for any ${\bf r} \in \R^D$.
  \item Either there is a vectorial space $\GG \supset \mathcal{I} P_\GG$
    endowed with a norm $\| \cdot \|_\GG$ such that we can define a
    distance $\mbox{dist}_\GG$ on $P_{\GG,\hbox{\small\bf r}}$ for any
    ${\bf r} \in \R^D$ by setting
  $$
  \forall \, f, \, g \in P_{\GG,\hbox{\small\bf r}}, \qquad
  \mbox{dist}_\GG(f,g) := \| g- f \|_\GG.
  $$
  \end{itemize}
  
  Finally, we say that two metrics $d_0$ and $d_1$ on $P_\GG$ are {\it
    topologically uniformly equivalent on bounded sets} if there exists
  $\kappa \in (0,\infty)$ and for any $a \in (0,\infty)$ there exists $C_a
  \in (0,\infty)$ such that
  $$
  \forall \, f, \, g \in \BB P_{\GG,a} \quad
  d_0(f,g) \le C_a \, \left[d_1(f,g)\right]^\kappa, \quad 
  d_1(f,g) \le C_a \, \left[d_0(f,g)\right]^\kappa.
  $$
  If $d_0$ and $d_1$ are resulting from some normed spaces $\GG_0$ and
  $\GG_1$, we abusively say that $\GG_0$ and $\GG_1$ are {\it topologically
    uniformly equivalent (on bounded sets)}.
  \end{defin}

\begin{ex}\label{expleTV} 
  The choice $m_\GG := 1$, ${\bf m}_\GG := 0$, $\| \cdot \|_\GG := \|\cdot
  \|_{M^1}$ recovers $P_\GG(E) = P(E)$. More generally on can choose
  $m_{\GG_k}(v) := \hbox{dist}_E(v,v_0)^{k}$, ${\bf m}_{\GG_k} := 0$, $\|
  \cdot \|_{\GG_k} := \|\cdot \, \hbox{dist}_E(v,v_0)^{k} \|_{M^1}$. For
  $k_1 > k_2,k_3 >0$, the spaces $P_{\GG_{k_2}}$ and $P_{\GG_{k_3}}$
  are topologically uniformly equivalent on bounded sets of
  $P_{\GG_{k_1}}$.
\end{ex}

\begin{ex}\label{expleWeak} 
  There are many distances on $P(E)$ which induce the weak topology, see
  for instance \cite{BookRachev}. In section~\ref{subsec:ExpleMetrics}  below, we
  will present some of them which have a practical interest for us, and
  which are all topologically uniformly equivalent on ``bounded sets'' of
  $P(E)$, when the bounded sets are defined thanks to a convenient (strong
  enough) weight function.
\end{ex}


\subsection{Differential calculus for functions of probability measures }
\label{sec:diff-measures}

We start with a purely metric definition in the case of usual H\"older
regularity.
\begin{defin}\label{def:Holdercalculus}
  For some metric spaces $\tilde\GG_1$ and $\tilde\GG_2$, some weight function $\Lambda : \tilde\GG_1 \mapsto \R_+^*$
  and  some $\eta \in (0,1]$, we denote by $C^{0,\eta}_\Lambda(\tilde \GG_1,\tilde\GG_2)$ the weighted
  space of functions from $\tilde \GG_1$ to $\tilde\GG_2$ with
  $\eta$-H\"older regularity, that is the functions $\SS  : \tilde  \GG_1 \to \tilde\GG_2$ such that there exists a constant $C >0$ so
  that
  \begin{equation}
    \label{eq:devphi}
    \forall \, f_1, \, f_2 \in \tilde \GG_1 \qquad 
    \mbox{{\em dist}}_{\GG_2} (\SS(f_1) , \SS(f_2) ) 
    \le C \,  \Lambda(f_1,f_2) \, 
    \mbox{{\em dist}}_{\GG_1} (f_1, f_2)^{\eta},
  \end{equation}
  with $\Lambda(f_1,f_2) := \max \{ \Lambda(f_1), \Lambda(f_2) \}$ and $ \mbox{{\em dist}}_{\GG_k}$ denotes the metric of $\tilde \GG_k$ (the tilde sign in the notation of the distance has been removed in order to present unified notation with the next definition). 
  We define the semi-norm $[ \cdot ]_ {C_\Lambda^{0,\eta}(\tilde \GG_1,\tilde\GG_2)}$ in $C_\Lambda^{0,\eta}(\tilde \GG_1,\tilde\GG_2)$ as the infimum of the constants $C > 0$ such that \eqref{eq:devphi} holds. 
 \end{defin}

Second we define a first order differential calculus, for which we
require a norm structure on the functional spaces.

\begin{defin}\label{def:C1kcalculus}
  For some normed spaces $\GG_1$ and $\GG_2$, some metric sets $\tilde
  \GG_1$ and $\tilde \GG_2$ such that $\tilde \GG_i - \tilde \GG_i \subset
  \GG_i$, some weight function $\Lambda : \tilde\GG_1 \mapsto [1,\infty)$
  and some $\eta \in (0,1]$, we denote $C_\Lambda^{1,\eta}(\tilde
  \GG_1, \GG_1; \tilde \GG_2, \GG_2)$ (or simply
  $C_\Lambda^{1,\eta}(\tilde \GG_1; \tilde \GG_2))$, the space of
  continuously differentiable functions from $\tilde \GG_1$ to $\tilde
  \GG_2$, whose derivative satisfies some weighted $\eta$-H\"older
  regularity. More explicitely, these are the continuous functions $\SS :
  \tilde \GG_1 \to \tilde \GG_2$ such that there exists a continuous
  fonction $D \SS : \tilde \GG_1 \to \BB(\GG_1,\GG_2)$ (where
  $\BB(\GG_1,\GG_2)$ denotes the space of bounded linear applications from
  $\GG_1$ to $\GG_2$ endowed with the usual norm operator), and some
  constants $C_i >0$, $i=1, \, 2, \, 3$, so that for any $f_1, \, f_2
  \in \tilde \GG_1$:
  \bear
    \label{eq:devdist1}\qquad
    \left\| \SS(f_2) - \SS(f_1)\right\|_{\GG_2} \!\!&\le&\!\!  C_1 \,
    \Lambda(f_1,f_2)\, \|f_2 - f_1 \|_{\GG_1}
    \\
    \label{eq:devdist2}\qquad
    \left\| \left \langle D \SS[f_1] , f_2 - f_1 \right \rangle
    \right\|_{\GG_2} \!\!&\le&\!\!  C_2 \, \Lambda(f_1,f_2)\,
    \|f_2 - f_1 \|_{\GG_1}
    \\
    \label{eq:devdist3}\qquad
    \left\| \SS(f_2) - \SS(f_1) - \left \langle D \SS[f_1] ,
        f_2 - f_1 \right \rangle \right\|_{\GG_2} \!\!&\le&\!\! C_3
    \, \Lambda(f_1,f_2) \, \|f_2 - f_1 \|_{\GG_1}^{1+\eta}.
    \eear 
 \end{defin}
 
 \begin{notation} 
  For $\SS \in C^{1,\eta}_\Lambda$, we define $C_i^\SS$, $i=1, \, 2, \, 3$, as the infimum of the
    constants $C_i > 0$ such that \eqref{eq:devdist1}
    (resp. \eqref{eq:devdist2}, \eqref{eq:devdist3}) holds.
  We then denote
  \[ [\SS
  ]_{C^{0,1}_\Lambda} := C_1^\SS, \quad [\SS ]_{C^{1,0}_\Lambda} :=
  C_2^\SS, \quad [\SS ]_{C^{1,\eta}_\Lambda} := C_3^\SS, \quad \| \SS
  \|_{C^{1,\eta}_\Lambda} := C_1^\SS + C_2^\SS+ C_3^\SS,
  \]
  and we will remove the subscript $\Lambda$ when $\Lambda \equiv
  1$. 
 \end{notation}
 
 \begin{rem} In the sequel, we shall apply this differential calculus
   with some suitable subspaces $\tilde \GG_i \subset P(E)$. This
   choice of subspaces is crucial in order to make rigorous the
   intuition of Gr\"unbaum \cite{Grunbaum} (see the --- unjustified
   --- expansion of $H_f$ in \cite{Grunbaum}). It is worth emphasizing
   that our differential calculus is based on the idea of considering
   $P(E)$ (or subsets of $P(E)$) as ``plunged sub-manifolds'' of some
   larger normed spaces $\GG_i$. Our approach thus differs from the
   approach of P.-L.  Lions recently developed in his course at
   Coll\`ege de France \cite{PLL-cours} or the one developed by
   L. Ambrosio et al in order to deal with gradient flows in probability
   measures spaces, see for instance \cite{AmbrosioBook}. In the
   sequel we develop a differential calculus in probability measures
   spaces into a simple and robust framework, well suited to deal with
   the different objects we have to manipulate ($1$-particle
   semigroup, polynomial, generators, ...). And the main innovation
   from our work is the use of this differential calculus to state
   some subtle ``differential'' stability conditions on the limiting
   semigroup. Roughly speaking the latter estimates measure how this
   limiting semigroup handles fluctuations departing from
   chaoticity. They are the corner stone of our analysis. Surprisingly
   they seem new, at least for Boltzmann type equations.
\end{rem}

This differential calculus behaves well for composition in the sense that 
for any given $\UU \in C^{1,\eta}_{\Lambda_\UU} (\tilde \GG_1, \tilde \GG_2)$ 
and $\VV \in  C^{1,\eta}_{\Lambda_\VV} (\tilde \GG_2, \tilde \GG_3)$ there holds
$\SS:= \VV \circ \UU \in C^{1,\eta} _{\Lambda_\SS} (\tilde \GG_1, \tilde \GG_3)$ 
for some appropriate weight function $\Lambda_\SS$. 
We conclude the section by stating a precise result well adapted to our applications.
The proof is straightforward by writing and compounding the expansions of $\UU$ and 
$\VV$ provided by Definition~\ref{def:C1kcalculus} and we then skip it. 

\begin{lem}\label{lem:DL} 
For any given $\UU \in C^{1,\eta}
  _{\Lambda} (\tilde \GG_1, \tilde \GG_2)$ and $\VV \in
  C^{1,\eta} (\tilde \GG_2, \tilde \GG_3)$ there holds
  $\SS:= \VV \circ \UU \in C^{1,\eta} _{\Lambda^{1+\eta}} (\tilde \GG_1,
  \tilde \GG_3)$ and $D \SS [f] = D \VV [\UU(f)] \circ D \UU[f]$. 
More precisely, there holds
 \begin{eqnarray*}
    [ \SS ]_{C^{0,1}_{\Lambda}} \le [\VV
    ]_{C^{0,1}} \, [\UU]_{C^{0,1}_\Lambda}, \quad [\SS
    ]_{C^{1,0}_{\Lambda}} \le [\VV ]_{C^{1,0}} \,
    [\UU]_{C^{1,0}_\Lambda} 
  \end{eqnarray*}
  and 
  \begin{eqnarray*} 
   [\SS ]_{C^{1,\eta}_{\Lambda^{1+\eta}}} \le [\VV
    ]_{C^{1,0}} \, [\UU ]_{C^{1,\eta}_\Lambda} + [\VV ]_{C^{1,\eta}} \,
    [\UU ]_{C^{0,1}_\Lambda}^{1+\eta}.
  \end{eqnarray*}
When further $\VV \in C^{1,1} (\tilde \GG_2, \tilde \GG_3)$, we also 
have  $\SS:= \VV \circ \UU \in C^{1,\eta}_{\Lambda^2} (\tilde \GG_1, \tilde \GG_3)$ 
with
  \begin{eqnarray*} 
   [\SS ]_{C^{1,\eta}_{\Lambda^2}} \le [\VV
    ]_{C^{1,0}} \, [\UU ]_{C^{1,\eta}_\Lambda} + [\VV ]_{C^{1,1}} \,
    [\UU ]_{C^{0,(1+\eta)/2}_\Lambda}^{2}.
  \end{eqnarray*}
  \end{lem}

\subsection{The pushforward generator}
\label{sec:calculus-gen}

As a first example of application of this differential calculus, let
us compute the generator of the pushforward limiting semigroup.  
Assume the following:

\fbox{
\begin{minipage}{0.9\textwidth}
\begin{itemize}
\item[{\bf (A2)}] {\bf Existence of the generator of the pushforward
    semigroup.} For some normed space $\GG_1$ and some probability
  space $P_{\GG_1}(E)$ (defined as above) associated to a weight
  function $m_1$ and constraint function ${\bf m_1}$, and
  endowed with the metric induced from $\GG_1$, for some $\delta
  \in (0,1]$ and some  $\bar a\in (0,\infty)$ we have for any $a\in (\bar a,\infty)$:
  \begin{itemize}
  \item[(i)] For any $t \in [0,+\infty)$, $S^{N \! L}_t : \BB P_{\GG_1,a}
    \to \BB P_{\GG_1,a}$ is continuous (for the metric
    $\mbox{dist}_{\GG_1}$), uniformly in time.
  
  \item[(ii)] The application $Q$ is $\delta$-H\"older continuous from
    $\BB P_{\GG_1,a}$ into $\GG_1$.

  \item[(iii)] For any $f \in \BB P_{\GG_1,a}$, for some $\tau >0$ the
    application $[0,\tau) \to \BB P_{\GG_1,a}$, $t
    \mapsto \SS^{N \! L}_t(f)$ is $C^{1,\delta}([0,\tau);\BB P_{\GG_1,a})$, with $\SS(f)'(0) =
    Q(f)$.
\end{itemize}
\end{itemize}
\end{minipage}
}

\begin{lem}\label{lem:H0} Under assumption {\bf (A2)} the generator
  $G^\infty$ of the pushforward semigroup $T^\infty_t$ exists as an
  unbounded linear operator on $C(\BB P_{\GG_1,a}(E);\R)$ for any $a
  \in (\bar a,\infty)$ with domain including $C^{1,\delta}(\BB
  P_{\GG_1,a}(E); \R)$, and on this domain it is defined by the
  formula
  \begin{equation}
    \label{eq:formulaGinfty}
    \forall \, \Phi \in C^{1,\delta}(\BB P_{\GG_1,a}(E); \R), \ 
    \forall \, f \in \BB P_{\GG_1,a}(E), \quad 
    \left( G^\infty \Phi \right) (f) :=
    \left \langle D\Phi[f], Q(f)\right\rangle.
  \end{equation}
\end{lem}

\noindent {\sl Proof of Lemma~\ref{lem:H0}.} We split the proof
in several steps.

\smallskip \noindent {\sl Step 1.} First $(T^\infty_t)$ is a
continuous semigroup on $C(\BB P_{\GG_1,a}(E);\R)$: consider $\Phi \in
C(\BB P_{\GG_1,a}(E);\R)$ and some sequence $(f_n)$ of $\BB
P_{\GG_1,a}(E)$ such that $\mbox{dist}_{\GG_1}(f_n,f) \to 0$,
then thanks to {\bf (A2)}-(i) we deduce $(T^\infty_t \Phi)(f_n) =
\Phi(S^{N\!L}_t(f_n)) \to \Phi(S^{N\!L}_t(f)) = (T^\infty_t
\Phi)(f)$, and by composition $T^\infty_t \Phi \in C(\BB
P_{\GG_1,a}(E);\R)$. 

Next, we have
$$
\| T^\infty_t \| = \sup_{\|\Phi \| \le 1} \| T^\infty_t \Phi \| =
\sup_{\|\Phi \| \le 1} \sup_{f \in \BB P_{\GG_1,a} } |
\Phi(S^{N\!L}_t(f))| \le 1, \qquad \|\Phi \| = \sup_{g \in
  \BB P_{\GG_1,a} } |\Phi(g)|.
$$
Now, from {\bf (A2)}-(iii) there exists a modulus of continuity
$\omega$ such that $\| S^{N\!L}_t f - f \|_{\GG_1} \le \omega
(t) \to 0$ as $t \to 0$ for any $f \in \BB P_{\GG_1,a}$
(e.g. $\omega(t) = t^\alpha$, $\alpha \in (0,1)$), which implies
$$
\forall \, \Phi \in C(\BB P_{\GG_1,a}(E);\R), \qquad \| T^\infty_t \Phi -
\Phi \| = \sup_{f \in \BB P_{\GG_1,a} } |\Phi(S^{N\!L}_t(f)) -
\Phi(f)|\to 0.
$$
We therefore deduce from Hille-Yosida's Theorem that $(T^\infty_t)$
has a closed generator $G^\infty$ with dense domain included in 
$C(\BB P_{\GG_1,a}(E);\R)$.

\smallskip \noindent {\sl Step 2.} 
Let us define $\tilde G^\infty \Phi $ by
$$ 
\forall \, \Phi \in C^{1,0}(\BB P_{\GG_1,a}; \R), \ \forall \, f \in
\BB P_{\GG_1,a}, \quad ( \tilde G^\infty \Phi ) (f) := \left \langle
  D\Phi[f], Q(f)\right\rangle.
$$
The right-hand side quantity is well defined since $D\Phi(f) \in
\BB(\GG_1,\R) = \GG_1'$ and $Q(f) \in \GG_1$. Moreover,
since both $f \mapsto D\Phi[f]$ and $f \mapsto
Q(f)$ are continuous 
we have $\tilde G^\infty \Phi \in C(\BB P_{\GG_1,a}; \R)$.

\smallskip \noindent {\sl Step 3.} Consider $\Phi \in C^{1,\delta}(\BB
P_{\GG_1,a}; \R)$. By the composition rule of Lemma~\ref{lem:DL}, for any
fixed $f \in P_{\GG_1}(E)$, $t \mapsto T^\infty _t \Phi (f) =
\Phi \circ S^{N\! L}_t (f)$ is $C^{1,\delta}([0,\eta);\R)$ and
\begin{eqnarray*} 
  {d \over dt} (T^\infty_t \Phi) (f) \big|_{t=0}
  &:=&  {d \over dt} \big(\Phi \circ \SS(f)(t)\big) \big|_{t=0} \\
  &=& \left\langle D\Phi (\SS(f))\big|_{t=0}, 
    {d \over dt} \SS(f)\big|_{t=0} \right\rangle \\
  &=& \left\langle D\Phi[f], Q(f) \right\rangle = \left( G^\infty \Phi
  \right) (f),
\end{eqnarray*}
which precisely means that any such $\Phi$ belongs to
$\mbox{domain}(G^\infty)$ and that \eqref{eq:formulaGinfty} holds.
\qed


\subsection{Compatibility of the $\pi^N$ and $R^\ell$ transformations}
\label{sec:compatibility}

Our transformations $\pi^N$ and $R^\ell$ behave nicely for the sup
norm on $C_b(P(E),TV)$, see~\eqref{eq:compat:infty}.  More generally
we shall consider ``duality pairs'' of metric spaces:
\begin{defin}
We say that a pair $(\FF,\GG)$ of normed vectorial spaces 
are ``in duality'' if
\begin{equation} \label{eq:dualiteFG} \forall \, f \in \GG,
  \,\, \forall \, \phi \in \FF \qquad |\langle f, \phi \rangle
  |\le \|f \|_\GG \, \| \phi \|_\FF
\end{equation} 
where $\|\cdot \|_\FF$ denotes the norm on
$\FF$ and $\| \cdot \|_\GG$ denotes the norm on $\GG$.
\end{defin}

The ``compatibility'' of the transformation $R^\ell$ for any such pair
follows from the multilinearity: if $\FF$ and $\GG$ are in duality,
$\FF \subset C_b(E)$ and $P_\GG$ is endowed with the metric associated
to $\|\cdot \|_\GG$, then for any $\varphi = \varphi_1 \times \dots
\times \varphi_\ell \in \FF^{\otimes \ell}$, the polynomial function
$R^\ell_\varphi$ is of class $C^{1,\eta}(P_\GG,\R)$ for any $\eta \in (0,1]$. Indeed, defining for $f_\alpha \in P_{\GG_1}$,
$\alpha = 1, 2$, 
\[ 
\GG \to \R, \quad h \mapsto DR^\ell_\varphi [f_1] (h) :=
\sum_{i=1}^\ell \left( \prod_{j \not= i } \int_E \varphi \, df_1 \right)
\langle \varphi_i, h \rangle,
\]
we have
\begin{eqnarray*}
  R^\ell_\varphi(f_2) - R^\ell_\varphi(f_1) &=& \sum_{i=1}^\ell
  \left( \prod_{1 \le k < i } \int_E \varphi_k \, df_2 \right) \, \langle
  \varphi_i, f_2 - f_1 \rangle \, \left( \prod_{i < k \le \ell }
  \int_E \varphi_k \, df_1 \right), 
\end{eqnarray*}
and
\begin{multline}\nonumber
  R^\ell_\varphi(f_2) - R^\ell_\varphi(f_1) - DR^\ell_\varphi
  [f_1] (f_2 - f_1) =
  \\
  = \sum_{1\le j < i \le \ell} \left( \prod_{1 \le k < j } \int_E
    \varphi_k df_2 \right) \, \langle \varphi_j, f_2 - f_1
  \rangle \, \left( \prod_{j < k <i} \int_E\varphi_k df_1 \right)
  \, \langle \varphi_i, f_2 - f_1 \rangle \Bigl( \prod_{i < k \le
    \ell } \int_E\varphi_k df_1 \Bigr).
\end{multline}
We deduce then for instance $R^\ell_\varphi \in C^{1,1}(\GG;\R)$ since
\begin{eqnarray} \nonumber
  && |R^\ell_\varphi(f_2) - R^\ell_\varphi(f_1) | \le \ell \, \| \varphi
  \|_{\FF \otimes (L^\infty) ^{\ell-1}} \, \| f_2 - f_1 \|_\GG , \quad
  | DR^\ell_\varphi [f_1] (h) | \le \ell \, \| \varphi \|_{\FF \otimes
    (L^\infty) ^{\ell-1}} \, \| h \|_\GG, \\ \label{eq:Ck1polyk}
  && \qquad |R^\ell_\varphi(f_2) - R^\ell_\varphi(f_1) -
  DR^\ell_\varphi [f_1] (f_2 - f_1)| \le \frac{\ell (\ell -1)}2
  \| \varphi \|_{\FF^2 \otimes (L^\infty) ^{\ell-2}} \, \| f_2 - f_1 \|^2_\GG ,  
\end{eqnarray}
where we have defined 
\begin{eqnarray*}
  \| \varphi \|_{\FF^k \otimes (L^\infty) ^{\ell-k}} 
  &:=& \max_{i_1, ...,  i_k \mbox{ {\tiny distincts in }} [|1,\ell|]} 
  \| \varphi_{i_1} \|_{\FF} \, \dots \| \varphi_{i_k} \|_{\FF} \! 
  \prod_{j \neq (i_1,\dots,i_k)} \| \varphi_j \|_{L^\infty(E)}.
 \end{eqnarray*}
 
\begin{rems}\label{rem:RphiC11} 
   \begin{itemize}
   \item It is easily seen in this computation that the tensorial
     structure of $\varphi$ is not necessary. In fact it is likely
     that this assumption could be relaxed all along our proof. We do
     not pursue this line of research.
   \item The assumption $\FF \subset C_b(E)$ could also be
     relaxed. For instance, when $\FF := \hbox{\rm Lip}_0(E)$ is the
     space of Lipschitz function which vanishes in some fixed point
     $x_0 \in E$, $\GG$ is its dual space, and $P_\GG := \{Êf
     \in P_1(E); \,\, \langle f, \hbox{\rm dist}_E(\cdot ,x_0)
     \rangle \le a \}$ for some fixed $a > 0$, we have $R^\ell_\varphi
     \in C^{1,1}(P_1(E);\R)$ with
     $$
     [R^\ell_\varphi]_{C^{0,1}} \le \ell \, a^{\ell-1} \, \|Ê\varphi
     \|_{\FF^{\otimes\ell}}, \quad [R^\ell_\varphi]_{C^{1,1}} \le
     \frac{\ell (\ell -1)}2 \, a^{\ell-1} \, \|Ê\varphi
     \|_{\FF^{\otimes\ell}},
     $$
     or equivalently $R^\ell_\varphi \in C^{1,1}_\Lambda(P_1(E);\R)$ with
     $\Lambda(f) := \|Êf \|_{M^1_1}^{\ell-1}$.
\end{itemize}
\end{rems}

\smallskip In the other way round, for the projection $\pi^N$ it is
clear that if the empirical measure $V \mapsto \mu^N_V$ belongs to
$C^{k,\eta}(E^N,\GG)$ for some norm space $\GG \in M^1(E)$, then by
composition one has
\begin{equation}\label{eq:compat-pi-F}
  \left\| \pi^N (\Phi) \right\|_{C^{k,\eta}(E ^N;\R)} 
  \le C_\pi \, \left\| \Phi \right\|_{C^{k,\eta}(\GG)}.
\end{equation}
However the regularity of the empirical measure depends on the metric
$\GG$.

\begin{ex}\label{expleTVbis} In the case $\FF = (C_b(E),L^\infty)$ and
  $\GG = (M^1(E),TV)$, \eqref{eq:compat-pi-F} is trivial with
  $k=\eta=0$.
\end{ex}

\begin{ex}\label{expleWpbis}  When $\FF = \mbox{Lip}_0(E)$ (Lipschitz
  function vanishing at some given point $v_0$) endowed with the norm
  $\|\phi \|_{Lip}$ and $P_\GG(E)$ (constructed in
  Example~\ref{expleWp}) is endowed with the Wasserstein distance
  $W_1$ with linear cost, one has \eqref{eq:compat-pi-F} with $k=0$,
  $\eta=1$:
  \[
  \left| \Phi\left(\mu^N _X \right) - \Phi\left(\mu^N _Y \right)
  \right| \le \| \Phi \|_{C^{0,1}(P_\GG)} \, W_1\left(\mu^N _X ,\hat
    \mu^N _Y \right) \le \| \Phi \|_{C^{0,1}(P_\GG)} \, \|X-Y
  \|_{\ell^1},
  \]
  where we use (\ref{Wqellq}), which proves that
  \[
  \| \pi^N (\Phi) \|_{C^{0,1} (E^N)} \le \| \Phi \|_{C^{0,1}(P_\GG)},
  \]
  when $E^N$ is endowed with the $\ell^1$ distance defined in
  (\ref{Wqellq}).
\end{ex}

 \subsection{Examples of distances on measures}
 \label{subsec:ExpleMetrics} 
 Let us list some well-known distances on $P(\R^d)$ (or on subsets of
 $P(\R^d)$) useful for the sequel. These distances are all topologically
 equivalent to the weak topology $\sigma (P(E),C_b(E))$ (on the sets
 $\BB P_{k,a}(E)$ for $k$ large enough and for any $a \in (0,\infty)$)
 and they are all uniformly topologically
 equivalent (see  \cite{TV,coursCT} and section~\ref{subsect:ComparisonDistance}). 

\begin{ex}[Dual-H\"older (or Zolotarev's) distances]\label{expleZolotarev} 
  Denote by $\mbox{dist}_E$ a distance on $E$ and let us fix $v_0 \in
  E$ (e.g. $v_0=0$ when $E = \R^d$ in the sequel). Denote by
  $C^{0,s}_0(E)$, $s \in (0,1)$ (resp. $\mbox{Lip}_0(E)$) the set of
  $s$-H\"older functions (resp. Lipschitz functions) on $E$ vanishing
  at one arbitrary point $v_0 \in E$ endowed with the norm
  $$
  [\varphi ]_{s} := \sup_{x,y \in E} {|\varphi(y) - \varphi(x)| \over
    \hbox{{\em dist}}_E(x,y)^s}, s \in (0,1], 
  \qquad  [\varphi ]_{\mbox{{\scriptsize {\em Lip}}}} := [\varphi ]_1.
  $$
  We then define the dual norm: take $m_\GG := 1$, ${\bf m}_\GG := 0$
  and $P_\GG(E)$ endowed with
 \begin{equation}\label{def;[]*s} 
   \forall \, f,g \in P_\GG,
  \qquad [g-f]^*_s := \sup_{\varphi \in C^{0,s}_0(E)} {\langle
    g - f, \varphi \rangle \over [\varphi]_s }. 
  \end{equation}
\end{ex}

\begin{ex}[Wasserstein distances]\label{expleWp} 
  For $q \in (0,\infty)$, define $W_q$ on
  $$
  P_\GG (E) = P_q(E):= \{ f \in P(E); \,\, \left \langle f,
    \mbox{dist}(\cdot,v_0)^q \right \rangle < \infty \}
  $$
  by
  \begin{eqnarray*}
    \forall \, f,g \in P_q(E), \quad  W_q(f,g) 
    :=  \inf_{\Pi \in \Pi(f,g)}  \int_{E\times E}  \hbox{dist}_E(x,y)^q \, \Pi(dx,dy),
\end{eqnarray*}
where $\Pi(f,g)$ denote the set of probability measures $\Pi \in
P(E \times E)$ with marginals $f$ and $g$ ($\Pi(A,E) = f
(A)$, $\Pi(E,A) = g (A)$ for any Borel set $A \subset E$).  Note
that for $V,Y \in E^N$ and any $q \in [1,\infty)$, one has
\begin{equation}\label{Wqellq} 
  W_q\left(\mu^N_V,\mu^N_Y\right) =
  d_{\ell^q(E^N/\SN)}
  (V, Y) := \min_{\sigma \in \SN} \left( {1 \over N} \sum_{i=1}^N
    \hbox{dist}_E(v_i,y_{\sigma(i)})^q \right)^{1/q}, 
\end{equation} 
and that
\begin{equation}\label{W1KR} 
  \forall \, f, \, g \in
  P_1(E) , \quad W_1 (f,g) 
  = [f-g]^*_1 = \sup_{\phi \in \mbox{\tiny Lip}_0(E)}\, \left \langle f-g,
  \phi \right \rangle.  
\end{equation} 
We refer to \cite{VillaniTOT} and the references therein for more details on the Wasserstein distances.
\end{ex}

\begin{ex}[Fourier-based norms] \label{expleFourier} 
For $E=\R^d$, $m_{\GG_1} := |v|$,
  ${\bf m}_{\GG_1} :=0$, let us define 
  \[ 
  \forall \, f \in \mathcal{T} P_{\GG_1}, \quad \| f \|_{\GG_1} = |f|_s
  := \sup_{\xi \in \R^d} \frac{|\hat f(\xi)|}{|\xi|^s}, \quad s \in (0,1].
  \]
Similarly, for $E=\R^d$, $m_{\GG_2} := |v|^2$,
  ${\bf m}_{\GG_2} :=v$, we define 
  \[ 
  \forall \, f \in \mathcal{T} P_{\GG_2}, \quad \| f \|_{\GG_2} = |f|_s
  := \sup_{\xi \in \R^d} \frac{|\hat f(\xi)|}{|\xi|^s}, \quad s \in (1,2].
  \]
 
\end{ex}

\begin{ex}[More Fourier-based norms] \label{expleFourierGal} 
More generally, for $E=\R^d$, $k \in \N^*$, we define $m_{\GG} := |v|^k$,
 ${\bf m}_{\GG} := (v^j)_{j \in \N^d, \, |j|Ê\le k-1}$ where for
$j = (j_1, ..., j_d)  \in \N^d$ we set $v^j = (v_1^{j_1}, ..., v_d^{j_d})$ and 
$|j| = j_1 + ... + j_d$, and
  \[ 
  \forall \, f \in \mathcal{T} P_{\GG}, \quad \| f \|_{\GG} = |f|_s
  := \sup_{\xi \in \R^d} \frac{|\hat f(\xi)|}{|\xi|^s}, \quad s \in (0,k].
  \]
In fact, we may extend the above norm to $M^1_k(\R^d)$ in the following way. 
We first define for $f \in M^1_{k-1}(\R^d)$ and $j \in \N^d, \, |j|Ê\le k-1$, 
$$
M_j[f] := \int_{\R^d} v^j\, f(dv).
$$
For a fixed (once for all) function $\chi$ in the Schwartz space $\SS(\R^d)$ such that 
$\chi \equiv 1$ on the set $\{Êv \in \R^d, |v| \le 1 \}$, so that in particular $\int_{\R^d} \FF^{-1} (\chi) (v) \, dv = \chi(0) = 1$,
we define $\MM_k[f]$ be its Fourier transform 
\[ \hat{\mathcal{M}_k}[f](\xi) 
:= \chi(\xi)  \, \left(
  \sum_{|j| \le k-1} M_j[f] \, \frac{\xi^j}{j!} \, i^{|j|} \right), \quad 
\]
which is some smooth version of the Taylor expansion of $\hat f$ at $\xi=0$. Then we may define the seminorms 
$$
| f |_k := \sup_{\xi \in \R^d} \left( |\xi|^{-k} \, \left| \hat f (\xi) - \hat{\mathcal{M}_k}[f](\xi) \right|  \right) 
$$
and 
$$
||| f |||_k := |f|_k + \sum_{j \in \N^d, \, |j|Ê\le k-1} |M_j[f]|.
$$
\end{ex}

\begin{ex}[Negative Sobolev norms]\label{expleH-s} 
  For any $s \in (d/2,d/2+1/2)$ take $E=\R^d$, $m_{\GG_1} := |v|$,
  ${\bf m}_{\GG_1} :=0$ and
  \[ 
  \forall \, f \in \mathcal{T} P_{\GG_1}, \quad \| f \|_{\GG_1}
  = \| f \|_{\dot H^{-s} (\R^d)} := \left\|
    \frac{\hat f(\xi)}{|\xi|^s} \right\|_{L^2}.
  \]
  Similarly, for any $s \in [d/2+1/2,d/2+1)$ take $E=\R^d$, $m_{\GG_2}
  := |v|^2$, ${\bf m}_{\GG_2} :=v$ and
  \[ 
  \forall \, f \in \mathcal{T} P_{\GG_2}, \quad \| f \|_{\GG_2}
  = \| f \|_{\dot H^{-s} (\R^d)} := \left\|
    \frac{\hat f(\xi)}{|\xi|^s} \right\|_{L^2}.
  \]
  \end{ex}

 \subsection{Comparison of distances when $E=\R^d$}
\label{subsect:ComparisonDistance}
 \begin{lem}\label{lem:ComparDistances} Let $f,g \in P(\R^d)$, then
\begin{equation}\label{estim:W1Wq}
  \forall \, q,k \in (1,\infty) \quad 
  W_1(f,g) \le W_q(f,g) \le  M_{k+1}^{1-\alpha} \, W_1(f,g)^\alpha,
\end{equation}
with $\alpha := 1 - (q-1)/k$, 
\begin{equation}\label{estim:dsWs}
  \forall \, s \in (0,1], \quad | f - g |_s \le W_s(f,g) \le W^s_1(f,g),
\end{equation}
\begin{equation}\label{estim:H-kd1}
  \forall \, s \in (d/2,d/2+1), \quad \| f-g \|_{\dot H^{-s}}
  \le  C \, |f-g |_1^{2s-d} , \quad C=C(d,s) >0,
\end{equation}
\begin{equation}\label{estim:*1s} 
\forall \, s > 0, \, \, k > 0 \quad [f-g]^*_1 \le C \,
  M_{k+1}^{\alpha_1} \, | f-g |_s^{\gamma_1},  
  \quad C=C(d,s,k) >0,
\end{equation}
with 
$$
\alpha_1 := \frac{d}{d+k+k(d+s-1)}, \quad 
\gamma_1 := \frac{k}{d+k+k(d+s-1)}, 
$$
\begin{equation}\label{estim:W1H-k}
  \forall \, s \ge 1, \, k > 0, \quad 
  [f-g]_1 ^* \le C \,  M_{k+1}^{\alpha_2} \,  \| f-g
  \|_{\dot H^{-s}}^{\gamma_2}, 
 \quad C=C(d,s,k) >0,
\end{equation}
with 
\[
\alpha_2 := \frac{d/2}{d/2+k+ k(s-1)}, 
\quad \gamma_2 := \frac{k}{d/2+k+ k(s-1)}
\]
and
$$
M_k := \max \left\{ \int_{\R^d} (1+|x|^k) \, f(dx) \, ; \ \int_{\R^d}
  (1+|x|^k) \, g(dx) \right\}.
$$
\end{lem}

 \noindent
 {\sl Proof of Lemma~\ref{lem:ComparDistances}. }  
We split the proof in several steps. For the proof of \eqref{estim:W1Wq} we refer to \cite{TV,coursCT}. 

 \smallskip\noindent{\sl Proof of \eqref{estim:dsWs}. } Let $\pi \in
 \Pi(f,g)$. We write
\begin{eqnarray*}
| \hat f (\xi) -  \hat g (\xi) | 
&=&
\left|  \int_{\R^d \times \R^d} (e^{-i \, v \cdot \xi} - e^{-i \, w \cdot \xi}) \, \pi (dv,dw) \right|
\\
&\le&
 \int_{\R^d \times \R^d} |e^{-i \, v \cdot \xi} - e^{-i \, w \cdot \xi}| \, \pi (dv,dw) 
\\
&\le&
 C_s \, \int_{\R^d \times \R^d} |v-w|^s \, |\xi|^s \, \pi (dv,dw),
\end{eqnarray*}
which yields \eqref{estim:dsWs} by taking the supremum in $\xi \in
\R^d$ and the infimum in $\pi \in \Pi(f,g)$.  
 
 \smallskip\noindent{\sl Proof of \eqref{estim:H-kd1}. } Consider $R >
 0$  and the ball $B_R=\{ x \in \R^d \ ; \ |x|\le R\}$, and write
\begin{eqnarray*}
\| f-g \|_{\dot H^{-s}}^2
&=& \int_{B_R} { | \hat f (\xi) -  \hat g (\xi) |^2  \over |\xi|^{2s}} \, d\xi 
+  \int_{B_R^c} { | \hat f (\xi) -  \hat g (\xi) |^2  \over |\xi|^{2s}} \, d\xi
\\
&\le& |f- g|^2_1 \int_{B_R} { d\xi  \over |\xi|^{2(s-1)}}
+  4 \, \int_{B_R^c} {d\xi  \over |\xi|^{2s}}
\\
&\le& C(d) \, R^{d - 2(s-1)} \, |f - g|^2_1  + 4 \, R^{d - 2s}.
\end{eqnarray*}
Then \eqref{estim:H-kd1} follows by choosing (the optimal) $R := |f - g|_1^{-1}$. 
  
\smallskip\noindent{\sl Proof of \eqref{estim:*1s}. } We introduce a
truncation function $\chi_R(x) = \chi(x/R)$, $R>0$, where $\chi \in
C^\infty(\R^d)$, $[\chi]_1 \le 1$, $0 \le \chi \le 1$, $\chi \equiv 1$ on
$B(0,1)$, supp$\, \chi \subset B(0,2)$, and a mollifer function
$\omega_\eps(x) = \eps^{-d} \, \omega(x/\eps)$, $\eps>0$ where for instance $\omega(x) =
(2\pi)^{-d/2} \, \exp(-|x|^2/2)$ (and thus $\hat \omega_\eps
(\xi) = \hat \omega  (\eps \, \xi) = \exp (- \eps^2 \, |\xi|^2/2)$). Fix $\varphi \in W^{1,\infty}(\R^d)$
such that $[\varphi]_1 \le 1$, $\varphi(0) = 0$, define $\varphi_R := \varphi \, \chi_R$,
$\varphi_{R,\eps} = \varphi_R \ast \omega_\eps$ and write
 $$
 \int \varphi \, (df-dg) = \int \varphi_{R,\eps} \, (df-dg) +
 \int \left(\varphi_R - \varphi_{R,\eps}\right) \, (df-dg) + 
 \int \left(\varphi - \varphi_R \right) \, (df-dg).
 $$ 
For the last term, we have 
\begin{eqnarray}
\label{estim:*1s-1}
  \left|  \int  (\varphi_R- \varphi) \, (df-dg) \right| 
  &\le&
   \int  (1- \chi_R) Ê\, \varphi \, (df+dg)
  \\  \nonumber
 &\le& \int_{B_R^c} [Ê\varphi ]_1 \, {|x|^{k+1} \over R^k} \, (df+dg) \le 
  \frac{M_{k+1}[f+g]}{R^k},
\end{eqnarray}
where $M_{k+1}[f+g]$ denotes the $(k+1)$-th moment of $f+g$.  
In order to deal with the second term, we observe that 
$$
|\nabla \varphi_R |Ê\le \chi (x/R) + |\varphi|Ê\, |\nabla (\chi_R)|
\le \chi (x/R) + {|x| \over R} Ê\, |\nabla \chi|(x/R),
$$
so that for any $q \in [1,\infty]$ there holds $\|Ê\nabla \varphi_R \|_{L^q}Ê \le C \, R^{d/q}$, for some constant depending only on $\chi$, $d$.
Next, using that
$$
\| \varphi_R - \varphi_{R,\eps} \|_\infty \le \| \nabla
\varphi_R \|_\infty \int_{\R^d} \omega_\eps (x) \, |x| \, dx \le C \,
\eps,
$$
we find
\begin{equation}\label{estim:*1s-2}
\left| \int  \left(\varphi_R - \varphi_{R,\eps} \right) \, (df  -
  dg) \right| \le C \, \eps.
\end{equation}
For the first term,  using Parseval's identity,
\begin{eqnarray}\nonumber
  \left| \int \varphi_{R,\eps} \, (f-g) \right|
  &=& {1 \over 2\pi} \left| \int \hat \varphi_R \, \hat \omega_\eps \,
    \overline{(\hat f - \hat g)} \, d\xi \right| 
  \\ \nonumber
  &\le& {1 \over 2\pi} \, \| \nabla \varphi_R \|_{L^1} \, 
  \left\| \frac{\hat f - \hat g}{|\xi|^s} \right\|_{L^\infty}  \,  \int |\xi|^{s-1} 
  \, \exp (- \eps^2 \, |\xi|^2/2) \, d\xi 
  \\  \nonumber
  &\le& C \, R^{d} \, \left( \int (1+|y|) \, \chi(y) \,
    dy \right) \, |f- g|_s  \, \eps^{-(d+s-1)}  \, 
  \left( \int
    |z|^{s-1} \, e^{-\frac{|z|^2}{2}} \, dz \right) 
    \\ \label{estim:*1s-3}
  &\le& C \, R^{d} \, \eps^{-(d+s-1)}  \, |f - g|_s. 
\end{eqnarray}
 Gathering \eqref{estim:*1s-1}, \eqref{estim:*1s-2} and
 \eqref{estim:*1s-3}, we get 
 \[ [f - g]_1 ^* \le C \, \left( \eps + \frac{M_{k+1}[f+g]}{R^k} +
   R^{d} \, \eps^{-(d+s-1)} \, |f - g|_s \right).
 \]
 This yields \eqref{estim:*1s} by optimizing the parameters $\eps$ 
 and $R$.
 
\smallskip\noindent{\sl Step 4. Proof of \eqref{estim:W1H-k}.} 
We start with the same decomposition as before:
$$
\int \varphi \, (df-dg) = \int \varphi_{R,\eps} \, (df-dg) +
\int \left(\varphi_R - \varphi_{R,\eps}\right) \, (df-dg) + \int
\left(\varphi - \varphi_R \right) \, (df-dg).
$$ 

The first term is controled by
\[
\left| \int \varphi_{R,\eps} \, (df-dg) \right| 
= \left| \int \hat \varphi_{R,\eps} \,
|\xi|^s \, \frac{(\hat f - \hat g)}{|\xi|^s} \right| \le 
\| \varphi_{R,\eps} \|_{\dot H^s} \, \| f - g \|_{\dot H^{-s}}
\]
with 
\begin{eqnarray}\nonumber
  \|  \varphi_{R,\eps} \|_{\dot H^s} 
  &=&  \left(  \int |\xi|^2 \, |\widehat{\varphi \, \chi_R}|^2 \, |\xi|^{2(s-1)} \, |\hat\omega_\eps|^2 \, d\xi \right)^{1/2}
  \\  \nonumber  
  &\le& \|Ê\nabla (\varphi \, \chi_R) \|_{L^2} \, \| |\xi|^{s-1} \, \hat  \omega_\eps (\xi) \|_{L^\infty}
 \\ \label{RachevW1-2}
   &\le& \|Ê\nabla (\varphi \, \chi_R) \|_{L^2} \, \eps^{1-s} \| |z|^{s-1} \, \hat  \omega (z) \|_{L^\infty}
\le  C \, R^{d/2} \,  \eps^{-(s-1)}.
\end{eqnarray} 
The second term and the last term are  controled as before  by $C \, \eps$ and $M_{k+1} \, R^{-k}$ respectively.
Summing we obtain 
\[ 
[f-g]_1 ^* \le C \, \left( \eps + \frac{M_{k+1}[f+g]}{R^{k}} +
  R^{d/2} \, \eps^{-(s-1)} \, |f-g|_{\dot H^{-s}} \right).
\]
This yields (\ref{estim:W1H-k}) by optimizing the parameters $\eps$ and
$R$. \qed

\subsection{On the law of large numbers for measures}
\label{subsec:LawLN}

\subsubsection{Remark on the meaning of the $\WW$ function.} 
For a
given function $D : P(E) \times P(E) \to \R_+$ continuous for the weak
topology, and such that $D(f,g) = 0$ if and only $f=g$ ($D$
stands for a distance on $P(E)$ or a function of a distance), we define
$$
\WW_{D} ^N (f) := \int_{E^{N}} D \left(\mu^N_V; f \right) \, 
f^{\otimes N} (dV).
$$
On the one hand, from the definition of $\pi^N_P f^{\otimes N}$ in section~\ref{framework}, 
we have
$$
\WW_{D} ^N (f) = \WW^N _D\left(f^{\otimes N}; f\right) = 
\WW_{D} ^\infty\left(\pi^N_P f^{\otimes N}; f\right),
$$
with
\begin{eqnarray*}
  &&\forall \,f^N \in P(E^N) \qquad  
  \WW^{N} _D (f^N; f):= \int_{E^N} D (\mu^N_V,f) \, f^{N} (dV) \\
  &&\forall \, \pi \in P(P(E)) \qquad 
  \WW_{D} ^\infty (\pi; f):= \int_{P(E)} D (\rho,f) \, \pi (d\rho) .
\end{eqnarray*}
  
  \Black
 
\subsubsection{Comparison of $\WW^N_D$ functions for different distances
  $D$.}
Let us begin with an elementary result.

 \begin{lem}\label{lem:Ap-element}
   If $D_1 \le C \, D_2^r$ for some constants $C > 0$ and $r < 1$, we have
   for some constant $C' = C'(C,r)$ and for any $f \in P(E)$
 \begin{equation}\label{Omegad1d2}
 \WW_{D_1} ^N(f) \le C' \, \left( \WW_{D_2} ^N(f) \right)^r.
 \end{equation}
 \end{lem}

 \noindent
 {\sl Proof of Lemma~\ref{lem:Ap-element}. }  For any $\eps > 0$, we have
 from Young's inequality
 \begin{eqnarray*}
   \WW_{D_1} ^N(f) 
   &\le& C \, \int_{E^N} D_2\left(\mu^N_V;f \right)^r \, f^{\otimes N} (dV) 
   \\
   &\le&
   C \int_{\R^{Nd}} \left[ (1-r) \, \eps^{r/(1-r)} +
     \frac{r}{\eps} \, D_2\left(\mu^N_V ; f\right) \right] \, 
   f^{\otimes N} (dV) 
   \\
   &\le&
   C(r) \, \left[  \eps^{r/(1-r)} + \frac{\WW_{D_2} ^N(f)}{\eps} \right].
 \end{eqnarray*}
This yields (\ref{Omegad1d2}) by optimizing $\eps > 0$. \qed

\begin{lem}\label{lem:Rachev&W1} 
  We have the following rates for the $\WW$ function:

  \begin{itemize}\item 
    For any $f \in P_2(\R^d)$, any $s \in (d/2,d/2+1)$ and any $N
    \ge 1$ there holds
    \begin{equation}\label{estim:RachevHdotk}
      \WW^N _{\| \cdot \|^2_{\dot H^{-s}}} (f) = 
      \int_{\R^{Nd}} \left\| \mu^N_V - f \right\|^2_{\dot H^{-s}} \, 
      f^{\otimes N} (dV) 
      \le \mbox{{\em Cst}}(d,M_2)  \, N^{- 1}.  
    \end{equation}

  \item For any $\eta > 0$ there exists $k \ge 1$ such that for any
    $f \in P_{k}(\R^d)$ and any $N \ge 1$ there holds
    \begin{equation}\label{estim:RachevW1}
      \WW^N _{W_1} (f) \le 
      \mbox{{\em Cst}}(\eta,k,M_k) \, N^{ - 1/ (d+\eta)}.
    \end{equation}

  \item For any $\eta > 0$ there exists $k \ge 2$ such that for any
    $f \in P_{k}(\R^d)$ and any $N \ge 1$ there holds
    \begin{equation}\label{estim:Rachev1}
      \WW^N _{W_2^2} (f) \le 
      \mbox{{\em Cst}}(\eta,k,M_k) \, N^{ - 1/ (d+\eta)}.
    \end{equation}
  \end{itemize}

 \end{lem}
 
 \begin{rem}\label{rem:Rachev} 
   Estimate \eqref{estim:Rachev1} has to be compared with the 
   following classical estimate established in \cite{BookRachev}: for
   any $f \in P_{d+5}(\R^d)$ and any $N \ge 1$ there holds
  \begin{equation}\label{estim:Rachev2}
    \WW^N _{W_2^2} (f) \le \mbox{{\em Cst}}(d,M_{d+5})\, N^{- {2 \over d+4}}.
 \end{equation}
It is worth mentioning that \eqref{estim:Rachev1} improve \eqref{estim:Rachev2} when $d \le 3$ and $k$ is large enough
(so that $\eta < 2 - d/2$). 
 \end{rem}

 \noindent
 {\sl Proof of Lemma~\ref{lem:Rachev&W1}. } 
 We split the proof into two steps. 

 \smallskip\noindent{\sl Proof of \eqref{estim:RachevHdotk}. } 
 Let us fix $f \in P_2(\R^d)$. First,
 writing
 \[
 \left(\hat\mu^N_V - \hat f \right)(\xi) = \frac{1}{N} \, \sum_{j=1}^N
 \left(e^{-i \, v_j \cdot \xi } - \hat f(\xi) \right) ,
 \]
 we have
 \begin{eqnarray*}
   \WW^N_{\| \cdot \|_{\dot H^{-s}}^2} (f) 
   &=&  \int_{\R^{Nd}} \left( \int_{\R^d} 
     \frac{\left|\hat\mu^N_V - \hat f \right|^2}{|\xi|^{2 \, s}} \,
     d\xi\right) \, f^{\otimes N} (dV) 
   \\
   &=& \frac{1}{N^2} \, \sum_{j_1,j_2=1}^N \int_{\R^{(N+1) d}} 
   \frac{\left(e^{-i \, v_{j_1} \cdot \xi} - \hat f(\xi)\right) \, 
     \left(\overline{ e^{-i \, v_{j_2} \cdot \xi} - \hat
         f(\xi)}\right)}
   {|\xi|^{2 \, s}} \, d\xi  \, f^{\otimes N} (dV).
 \end{eqnarray*}
 We use 
 $$
  \int_{\R^d} (e^{-i \, v_j \cdot \xi}- \hat f(\xi) ) \, f(dv_j) = 0,
  \quad j=1, \dots, d,
  $$
 and
 \begin{eqnarray*}
 \int_{\R^d} \left| e^{-i \, v \cdot \xi } - \hat f(\xi) \right|^2 \, f(dv) 
 &=& \int_{\R^d} \left[ 1 - e^{-i \, v \cdot \xi } \,  
   \overline{\hat f(\xi)} - e^{i \, v \cdot \xi } \, 
   \hat f(\xi) + | \hat f(\xi) |^2 \right]  \, f(dv) \\
 &=& 1 - |\hat f(\xi) |^2, 
 \end{eqnarray*}
 to deduce
 \begin{eqnarray*}
   \WW^N_{\| \cdot \|_{\dot H^{-s}}^2} (f) 
   &=&
   \frac{1}{N^2} \, \sum_{j=1}^N \int_{\R^{(N+1) d}} 
   \frac{ \left| e^{-i \, v_j \cdot \xi } - \hat
       f(\xi)\right|^2}{|\xi|^{2 \, s}} \, d\xi  \, f^{\otimes N} (dV) 
     \\
     &=& \frac{1}{N} \, \int_{\R^{2 \,  d}} \frac{\left| e^{-i \, v \cdot
           \xi } - \hat f(\xi) 
       \right|^2}{|\xi|^{2 \, s}} \, d\xi  \, f(dv) 
     \\
     &=& \frac{1}{N} \, 
     \int_{\R^{d}} \frac{(1 - |\hat f(\xi) |^2)}{|\xi|^{2 \, s}} \, d\xi.
 \end{eqnarray*}
 
 Finally, observing that $\hat f(\xi) = 1 +i \, \langle f , v \rangle
 \cdot \xi + \OO(M_2 \, |\xi|^2)$, and therefore
 \begin{eqnarray*}
   |\hat f(\xi)|^2  
   &=& \left(1 +i \,  \langle f , v \rangle \cdot \xi + \OO(M_2 \,
   |\xi|^2)\right) \, 
   \left(1  - i \,  \langle f , v \rangle \cdot \xi 
     + \overline{\OO(M_2 \, |\xi|^2)}\right) \\
   &=& 1 + \OO\left(M_2 \, |\xi|^2\right), 
 \end{eqnarray*}
 we obtain
 \begin{eqnarray*}
   \WW^N_{\| \cdot \|_{\dot H^{-s}}^2} (f) 
   &=&
   \frac{1}{N} \, \left( \int_{|\xi| \le 1} \frac{(1 - |\hat f(\xi)
       |^2)}{|\xi|^{2 \, s}} 
     \, d\xi +  \int_{|\xi| \ge 1} \frac{(1 - |\hat f(\xi)
       |^2)}{|\xi|^{2 \, s}} 
     \, d\xi \right)
   \\
   &=&
   \frac{1}{N} \, \left( \int_{|\xi| \le 1} \frac{M_2}{|\xi|^{2 \,
         (s-1) }} \, 
     d\xi +  \int_{|\xi| \ge 1} \frac{1}{|\xi|^{2 \, s}} \, d\xi \right),
 \end{eqnarray*}
 from which (\ref{estim:RachevHdotk}) follows. 
 
 \smallskip\noindent{\sl Proof of \eqref{estim:RachevW1} and \eqref{estim:Rachev1}. } 
 By gathering \eqref{estim:RachevHdotk}, (\ref{estim:W1H-k}) in
 Lemma~\ref{lem:Ap-element} and Lemma~\ref{lem:ComparDistances}, we
 straightforwardly deduce
 \begin{eqnarray*}
   \WW^N_{W_1} (f)  
   &=& \int_{R^{Nd}} [ \mu^N_V - f]^*_1 \, f^{\otimes N}(dV) \\
   &\le& C_{k,s,d}(M_{k+1}) \int_{R^{Nd}} 
   \left(  \|\mu^N_V - f \|_{\dot H^{-s}}^2 \right)^{\gamma_2/2} \, 
   f^{\otimes N}(dV)   \\
   &\le&C_{k,s,d}(M_{k+1}) \, N^{ - {\gamma_2/2}},
 \end{eqnarray*}
from which we deduce  (\ref{estim:RachevW1}) because in the limit case $k=\infty$ we have $\gamma_2/2 = 1/(2s)$ and we may choose $s$ as
close from $d/2$ as we wish. 

\smallskip\noindent{\sl Proof of \eqref{estim:Rachev1}. } Estimate
\eqref{estim:Rachev1} follows from \eqref{estim:RachevW1} with the
help of Lemma~\ref{lem:ComparDistances} and \eqref{estim:W1Wq}.  \qed

\section{The abstract theorem}
\label{sec:abstract-theo}

\subsection{Assumptions for the abstract theorem}
\label{sec:hyp-evol}

Assume that 
\begin{itemize}

\item {\bf (A1) } and {\bf (A2) } hold, so that in particular the
  semigroups $S^N_t$, $T^N_t$, $S^{N \! L}_t$ and $T^\infty_t$ are
  well defined as well as the generators $G^N$ and $G^\infty$.

\end{itemize}

\fbox{
\begin{minipage}{0.95\textwidth}
\begin{itemize}
\item[{\bf (A3)}] {\bf Convergence of the generators.} In the
  probability metrized set $P_{\GG_1}$ introduced in {\bf (A2)}
  (associated to the weight function $m_{\GG_1}$ and constraint
  function ${\bf m}_{\GG_1}$) we define 
  $$
  {\bf R_{\GG_1}} := \{Ê{\bf r} \in \R^D; \,\, \exists \, f \in P_{\GG_1} \, \hbox{s.t.} \, {\bf m}_{\GG_1} (f) = {\bf r} \}.
  $$
  Then  for the weight function $\Lambda_1 (f) = \langle f,m_1 \rangle$ and for some function $\eps_2(N)$ going to $0$ as $N$ goes to  
  infinity, we assume that the generators $G^N$ and $G^\infty$ satisfy
  \begin{eqnarray}
  \nonumber
  && \forall \,  \Phi \in \bigcap_{{\bf r} \in   {\bf R_{\GG_1}}} C_{\Lambda_1}^{1,\eta}(P_{\GG_1,{\bf r}};\R)
  \\
   \label{eq:estimcvg1}
  &&
   \left\|  \left( M^N _{m_1} \right)^{-1} \!\!
      \left( G^N \, \pi_N - \pi_N \, G^\infty \right)  \,  
      \Phi \right\|_{L^\infty(\mathbb{E}_N)} \le \eps_2(N) \, 
   \sup_{{\bf r} \in   {\bf R_{\GG_1}}} [ \Phi ]_{C^{1,\eta} _{\Lambda_1}(P_{\GG_1,{\bf r}})},
  \end{eqnarray}
 where $M^N _{m_1}$ was defined in~\eqref{eq:defMm}. 
\end{itemize}
\end{minipage}
}

\fbox{
\begin{minipage}{0.95\textwidth}
\begin{itemize}  
\item[{\bf (A4)}] {\bf Differential stability of the limiting semigroup.} 
   We assume that the flow $S^{N \! L}_t$ is $C^{1,\eta}_{\Lambda_2}(P_{\GG_1, {\bf r}},P_{\GG_2})$ for any ${\bf r} \in   {\bf R_{\GG_1}}$ in the sense that 
   there exists $C_T ^\infty >0$ such that
  \begin{equation}
   \sup_{{\bf r} \in   {\bf R_{\GG_1}}}  \int_0^T  \left([ S^{N \! L}_t]_{C^{1,\eta}_{\Lambda_2}(P_{\GG_1, {\bf r} },P_{\GG_2})} 
    + [ S^{N \!
      L}_t]^{1+\eta'}_{C^{0,\eta''}_{\Lambda_2}(P_{\GG_1,  {\bf r} },P_{\GG_2})} \right) 
    \, dt  \le C_T ^\infty,
  \end{equation}
  where $\eta \in (0,1)$ is the same as in {\bf (A3)}, $(\eta',\eta'') = (\eta,1)$ or $(\eta',\eta'') = (1,(1+\eta)/2)$,  $\Lambda_2 =
  \Lambda_1^{1/(1+\eta')}$ and $P_{\GG_2}Ê= \{Ê\rho, \, \rho \in P_{\GG_1} \}Ê$ but it is endowed with the norm associated to a normed space $\GG_2 \supset \GG_1$.
\end{itemize}
\end{minipage}
}

\fbox{
\begin{minipage}{0.95\textwidth}
\begin{itemize}   
\item[{\bf (A5)}] {\bf Weak stability of the limiting semigroup.} 
We assume that, for some probabilistic space $P_{\GG_3}(E)$  (associated to a weight function $m_{\GG_3}$, a constraint function ${\bf m}_{\GG_3}$ and some metric structure $\hbox{dist}_{\GG_3}$) and that for any $a,T > 0$ there exists a concave and continuous function $\Theta_{a,T} : \R_+ \to \R_+$ such that $\Theta_{a,T}(0) = 0$, we have
 \begin{eqnarray}
  \nonumber
  &&  \forall  \, f_1, f_2 \in \BB P_{\GG_3,a}(E)
  \\
   \label{eq:stab-weak-S}
  &&
   \sup_{[0,T)} \mbox{dist}_{\GG_3} 
    \left( S^{N \! L}_t(f_1), S^{N \! L}_t(f_2) \right) \le \Theta_{a,T} \left( 
    \mbox{dist}_{\GG_3} (f_1,f_2) \right).
  \end{eqnarray}
\end{itemize}
\end{minipage}
}

\subsection{Statement of the result}
\label{sec:abstract-result}

\begin{theo}[Fluctuation estimate]\label{theo:abstract}
  Consider a family of $N$-particle initial conditions $f^N_0 \in
  P_{sym}(E^N)$, $N \ge 1$, and the associated solution $f_t^N = S^N_t
  f_0^N$. Consider a $1$-particle initial condition $f_0 \in P(E)$ and
  the associated solution $f_t = S^{N \! L}_t f_0$. Assume that {\bf
    (A1)-(A2)-(A3)-(A4)-(A5)} hold for some spaces $P_{\GG_k}$,
  $\GG_k$ and $\FF_k$, $k=1,2,3$ with $\FF_k \subset C_b(E)$, and
  where $\FF_k$ and $\GG_k$ are in duality.

  Then there is an explicit absolute constant $C \in (0,\infty)$ such
  that for any $N, \ell \in \N^*$, with $N \ge 2 \ell$, and for any
  \[ \varphi = \varphi_1 \otimes \varphi_2 \otimes \dots \otimes \,
  \varphi_\ell \in (\FF_1 \cap \FF_2 \cap \FF_3)^{\otimes \ell} \]
  we have
  \begin{eqnarray}
  \label{eq:cvgabstract1}
  &&\quad \sup_{[0,T)}\left| \left \langle \left( S^N_t(f_0 ^N) - \left(
          S^{N \! L}_t(f_0) \right)^{\otimes N} \right), \varphi 
    \right\rangle \right| 
  \\ \nonumber 
  &&\quad 
  \le C \, \Bigg[ \ell^2 \, \frac{\|\varphi\|_\infty}{N} 
  + C^N_{T,m_1} \, C_T^\infty \, \varepsilon_2(N) \, \ell^2 \, 
  \|\varphi\|_{\FF_2 ^2 \otimes (L^\infty)^{\ell-2}} 
  \\ \nonumber
  &&\qquad\qquad\qquad 
  + \ell \, \, \|\varphi\|_{\FF_3 \otimes (L^\infty)^{\ell-1}} \, 
  \Theta_{C^N_{0,m_3},T}  \left(  \WW_{\!\! \hbox{ \small {\em dist}}_{\GG_3}} \!\! 
    \left(\pi^N_P f^N_0,\delta_{f_0}\right) \right) \Bigg],
  \end{eqnarray} 
  where $ \WW_{\!\! \hbox{ \small {\em dist}}_{\GG_3}}$ stands for the Monge-Kantorovich distance
  in $P(P_{\GG_3}(E))$, see example~\ref{expleWp}, which means for
  such particular probabilities
  \begin{equation}\label{eq:defOmegaN} 
    \WW_{\!\! \hbox{ \small {\em dist}}_{\GG_3}} \!\!  \left(\pi^N_P f^N_0,\delta_{f_0}\right) = 
    \int_{E^N} \hbox{{\em  dist}}_{\GG_3}(\mu^N_V,f_0) \, f^N_0 (dV).
   \end{equation}
\end{theo}





\begin{rem}
  1) Our goal here is to treat the $N$-particles system as a
  perturbation (in a very degenerated sense) of the limiting problem,
  and to minimize assumptions on the many-particle systems in order to
  avoid complications of many dimensions dynamics.
  
  2) In the applications the worst decay rate in the right-hand side
  of \eqref{eq:cvgabstract1} is always the last one, which deals with
  the chaoticity of the initial data.
  
  3) It is worth mentioning that in the case  when $f^N_0 = f_0^{\otimes N}$ we have 
   $$
   \WW_{\!\! \hbox{ \small {\em dist}}_{\GG_3}} \!\! 
    \left(\pi^N_P f^N_0,\delta_{f_0}\right) 
    = \WW_{\!\! \hbox{ \small {\em dist}}_{\GG_3}}^N \! (f_0)
    $$
    and the decay rate is obtained thanks to the law of large numbers for measures 
    presented in section~\ref{subsec:LawLN}. For more general initial datum
    we refer to section~\ref{sec:extensions}. 
\end{rem}

\subsection{Proof of Theorems~\ref{theo:abstract}} 
For a given function $\varphi \in (\FF_1 \cap \FF_2 \cap \FF_3)^{\otimes
  \ell}$, we break up the term to be estimated into four parts:
\begin{eqnarray*}
  &&\left| \left \langle \left( S^N_t(f_0^{N}) - \left( S _t ^\infty
          (f_0)
        \right)^{\otimes N} \right),  \varphi  
      \otimes 1^{\otimes N-\ell} \right\rangle \right| \le
  \\
  &&\le \left| \left\langle S^N_t(f_0^N), 
      \varphi  \otimes
      1^{\otimes N-\ell} \right\rangle -
    \left \langle S^N_t(f_0^N), 
      R^\ell_\varphi \circ \mu^N_V \right\rangle  \right| \qquad\qquad (=:  \TT_ 1 )
  \\
  &&+ \left| \left\langle f_0^N, T^N_t ( R^\ell_\varphi \circ
      \mu^N_V) \right\rangle
    - \left\langle f_0^N, 
      (T_t ^\infty R^\ell_\varphi ) \circ \mu^N_V) \right\rangle  \right| \qquad\qquad (=:  \TT_ 2 )
  \\  
  &&+ \left| \left\langle f_0^N, 
      (T_t ^\infty R^\ell_\varphi ) \circ \mu^N_V) \right\rangle
    -  \left\langle (S _t ^\infty (f_0))^{\otimes \ell} ,
      \varphi \right\rangle \right|  \qquad\qquad\qquad (=:  \TT_ 3 ).
\end{eqnarray*} 

We deal separately with each part step by step: 
\begin{itemize} 
\item $\TT_1$ is controled by a purely combinatorial arguments
  introduced in \cite{Grunbaum}. In some sense it is the price we have
  to pay when we use the injection $\pi^N$;
\item $\TT_2$ is controled thanks to the consistency estimate {\bf
    (A3)} on the generators, the differential stability assumption
  {\bf (A4)} on the limiting semigroup and the moments propagation
  {\bf (A1)}; 
\item $\TT_3$ is controled in terms of the chaoticity of the initial
  data thanks to the weak stability assumption {\bf (A5)} on the
  limiting semigroup and {\bf (A1)}-(iii).
\end{itemize}

\smallskip

\noindent {\bf Step~1: Estimate of the first term $\TT_1$.} 
Let us prove that for any $t \ge 0$ and any $N \ge 2 \ell$ there holds
\begin{equation}\label{estim:T1}
  \TT_1 := \left| \left\langle S^N_t(f_0^N), 
      \varphi \otimes
      1^{\otimes N-\ell} \right\rangle - \left \langle S^N_t(f_0^N),
      R^\ell_\varphi \circ \mu^N_V \right\rangle \right| \le \frac{2 \,
    \ell^2 \, \| \varphi \|_{L^\infty(E^\ell)}}{N}.
\end{equation} 
Since $S^N_t(f_0^N)$ is a symmetric probability measure, estimate
(\ref{estim:T1}) is a direct consequence of the following lemma

\begin{lem}\label{lem:symmetrization}
For any $\varphi \in C_b(E^\ell)$ we have
\begin{equation}\label{estim:symmetrization1}
  \forall \, N \ge 2 \ell, \qquad  
  \left| \left( \varphi 
      \otimes {\bf 1}^{\otimes N-\ell} \right)_{sym} -
    \pi_N R^\ell_\varphi \right| 
  \le \frac{2 \, \ell^2 \, \| \varphi \|_{L^\infty(E^\ell)}}{N}
\end{equation}
where for a function $\phi \in C_b(E^N)$, we define its symmetrized
version $\phi_{sym}$ as:
\begin{equation}\label{def:sym}
  \phi_{\mbox{{\tiny sym}}} = 
  \frac{1}{|\SSS_N|} \, \sum_{\sigma \in \SSS_N} \phi_\sigma.
\end{equation}

As a consequence for any symmetric measure $f^N \in P(E^N)$ we have 
\begin{equation} \label{eq:mNRphi} \left| \langle f^N,
    R^\ell_\varphi(\mu^N _V) \rangle - \langle f^N, \varphi \rangle
  \right| \le {2 \, \ell^2 \, \| \varphi \|_{L^\infty(E^\ell)} \over
    N}.
\end{equation}
\end{lem}


\smallskip\noindent {\bf Proof of Lemma~\ref{lem:symmetrization}.} For a
given $\ell \le N/2$ we introduce 
\[ 
A_{N,\ell} := \left\{ (i_1, ..., i_\ell) \in [|1,N|]^\ell \, : \
  \forall \, k \not= k', \ i_k \not= i_{k'} \ \right\} \qquad
\mbox{and} \qquad B_{N,\ell} := A_{N,\ell}^c.
\]
Since there are $N (N-1) \dots (N-\ell+1)$ ways of choosing $\ell$
distinct indices among $[|1,N|]$ we get
\begin{eqnarray*}
  {\left|B_{N,\ell}\right|  \over N^\ell} &=& 
  1 - \left(1 - {1 \over N}\right) \, ... \, 
  \left(1 - {\ell-1 \over N}\right) 
  = 1 - \exp \left( \sum_{i = 0}^{\ell-1} 
    \ln \left(1 - \frac{i}N \right) \right) \\
  &\le& 1 - \exp \left( - 2 \sum_{i = 0}^{\ell-1} \frac{i}N \right) 
  \le {\ell^2 \over  N},
\end{eqnarray*}
where we have used 
\[
\forall \, x \in [0,1/2], \quad \ln ( 1 - x) \ge - 2 \, x \qquad \mbox{and}
 \qquad \forall \, x \in \R, \quad e^{-x} \ge 1 - x.
\]
Then we compute
\begin{eqnarray*} &&R^\ell_\varphi(\mu^N_V) =
  {1 \over N^\ell} \sum_{i_1, ..., i_\ell = 1}^N 
 \varphi (v_{i_1}, \dots, v_{i_\ell}) \\
  &&= {1 \over N^\ell} \sum_{(i_1, ..., i_\ell) \in A_{N,\ell}} 
  \varphi (v_{i_1}, \dots, v_{i_\ell})
  +  {1 \over N^\ell} \sum_{(i_1, ..., i_\ell) \in B_{N,\ell}}  
   \varphi (v_{i_1}, \dots, v_{i_\ell}) \\
  &&= {1 \over N^\ell} \, {1 \over (N-\ell)!} \sum_{\sigma \in \SSS_N}
  \varphi (v_{\sigma(1)}, \dots, v_{\sigma(\ell)})
  +  {\mathcal O} \left(  {\ell^2 \over N} \, \|\varphi \|_{L^\infty} \right) \\
  &&= {1 \over N!} \sum_{\sigma \in \SSS_N} 
  \varphi (v_{\sigma(1)}, \dots , v_{\sigma(\ell)})
  +  {\mathcal O} \left(  {2 \, \ell^2 \over N} \, 
    \| \varphi \|_{L^\infty} \right) \\
\end{eqnarray*}
and the proof of (\ref{estim:symmetrization1}) is complete. Next for any
$f^N \in P(E^N)$ we have
$$
\left\langle f^N, \varphi \right\rangle = \left\langle f^N,
  \left(\varphi \otimes {\bf
      1}^{\otimes N-\ell} \right)_{sym} \right\rangle,
 $$
 and \eqref{eq:mNRphi} trivially follows from
 (\ref{estim:symmetrization1}).  \qed

\medskip\noindent {\bf Step~2: Estimate of the second term $\TT_2$. }
Let us prove that for any $t \in [0,T)$ and any $N \ge 2 \ell$ there
holds 
\begin{eqnarray}\label{estim:T2}
  \quad \TT_ 2 &:= & \left| \left\langle
      f_0^N, T^N_t \left( R^\ell_\varphi \circ \mu^N_V \right) \right\rangle -
    \left\langle f_0^N, \left(T_t ^\infty R^\ell_\varphi \right) \circ
        \mu^N_V \right\rangle \right| \\ \nonumber &\le&
  C^N_{T,m_2} \, C_T ^\infty \, \| \varphi \|_{\infty, \FF_2^2 \otimes
    (L^\infty)^{\ell-2}} \, \ell^2 \, \eps(N).
\end{eqnarray} 
We start from the following identity
$$
T^N_t \pi_N - \pi_N T^\infty _t = - \int_0 ^t \frac{d}{ds} \left( T^N
  _{t-s} \, \pi_N \, T^\infty _s \right) \, ds = \int_0 ^t T^N _{t-s}
\, \left[ G^N \pi_N - \pi_N G^\infty \right] \, T^\infty _s \, ds.
$$
From assumptions {\bf (A1)} and {\bf (A3)}, we have for any $t \in
[0,T)$
\begin{eqnarray}\nonumber
  && \left|  \left\langle f_0^N, 
      T^N_t \left( R^\ell_\varphi \circ \mu^N_V \right) \right\rangle 
    - \left\langle f_0^N, \left(T_t ^\infty R^\ell_\varphi \right) \circ 
      \mu^N_V  \right\rangle \right|  \\
  \nonumber
  &&\qquad \le \int_0 ^T \left| \left\langle M^N_{m_1} \, 
      S^N _{t-s}\left(f_0 ^N\right), 
      \left(M^N _{m_1}\right)^{-1}  \, 
      \left[ G^N \pi_N - \pi_N G^\infty \right] \, 
      \left(T^\infty_s R^\ell_\varphi\right) \right\rangle \right| 
  \, ds \\ \nonumber
  && \qquad \le
  \left( \sup_{0 \le t < T} \left\langle f^N _t, M^N _{m_1} \right 
    \rangle \right) \,  
  \left( \int_0^T \left\| \left(M^N _{m_1}\right)^{-1} \, 
      \left[ G^N \pi_N - \pi_N G^\infty \right] \, 
      \left(T^\infty_s R^\ell_\varphi\right)
    \right\|_{L^\infty(\mathbb{E}_N)} \, ds \right) \\
  \label{eq:trotter}
  &&\qquad \qquad \qquad 
  \le  \eps(N) \, C_{T,m}^N \, \sup_{{\bf r}Ê\in {\bf R}_{\GG_1}}Ê\int_0^T 
  \left[ T^\infty_s R^\ell_\varphi\right]_{C^{1,\eta}
    _{\Lambda_1}(P_{\GG_1},{\bf r} )} \, ds. 
\end{eqnarray}

Now, let us fix ${\bf r}Ê\in {\bf R}_{\GG_1}$. Since $T^\infty_t(R^\ell_\varphi) = R^\ell_\varphi \circ
S^{N \!  L}_t$ with $S^{N \! L}_t \in C^{1,\eta}_{\Lambda_2}
(P_{\GG_1,{\bf r}};P_{\GG_2})$ thanks to assumption {\bf (A4)}, and
$R^\ell_\varphi \in C^{1,1}(P_{\GG_2};\R)$ because $\varphi \in
\FF_2^{\otimes\ell}$ (see subsection~\ref{sec:compatibility}), we
obtain with the help of Lemma~\ref{lem:DL} that
$T^\infty_t(R^\ell_\varphi) \in C^{1,\eta}
_{\Lambda_2^{1+\theta'}}(P_{\GG_1,{\bf r}};\R)$ with
\[
\left[ T^\infty_s\left(R^\ell_\varphi\right) \right]
_{C^{1,\eta} _{\Lambda_2^{1+\theta'}} (P_{\GG_1,{\bf r}})}
  \le \left( \left[ S^{N \! L}_t\right]_{C^{1,\eta}_{\Lambda_2}(P_{\GG_1, {\bf r}},P_{\GG_2})} 
  + \left[ S^{N \! L}_t\right]^{1+\theta'}_{C^{0,1}_{\Lambda_2}(P_{\GG_1, {\bf r}},P_{\GG_2})} \right)
  \, \left\| R^\ell_\varphi \right\|_{C^{1, \eta}(P_{\GG_2})}.
\]
With the help of $\Lambda_2 = \Lambda_1^{1/(1+\theta')}$, \eqref{eq:Ck1polyk} and
assumption {\bf (A4)}, we hence deduce
\begin{equation}
  \label{estim:TtPhi}
  \int_0^T [ T^\infty_s(R^\ell_\varphi) ]_{C^{1,\eta} _{\Lambda_1} (P_{\GG_1,{\bf r}})} \, ds
  \le C^\infty_T \,\ell^2 \,  \left(\| \varphi \|_{\FF_2 \otimes
      (L^\infty)^{\ell-2}} 
    + \| \varphi \|_{\FF_2 \otimes (L^\infty)^{\ell-1}} \right) 
\end{equation}
Then we go back to the computation \eqref{eq:trotter}, and plugging
(\ref{estim:TtPhi}) we deduce (\ref{estim:T2}).  

\medskip\noindent {\bf Step~3: Estimate of the third term $\TT_3$. }
Let us prove that for any $t \ge 0$, $N \ge \ell$ 
\begin{multline}\label{estim:T3}
  \TT_3 := \left| \left\langle f_0^N, \left(T_t ^\infty
        R^\ell_\varphi \right) \circ \mu^N_V \right\rangle -
    \left\langle \Big(S _t ^\infty (f_0)\Big)^{\otimes \ell} ,
      \varphi \right\rangle \right|  \le \\
  \le  \left[R_\varphi \right]_{C^{0,1}}  \,   \Theta_{C^N_{0,m_3},T}  \left( \WW_{1,P_{\GG_3}} \left(\pi^N_P f^N_0,\delta_{f_0}\right) \right).
\end{multline}

We shall proceed as in Step~2, using that:
\begin{itemize}
\item $\hbox{supp} \, \pi^N_P f_0^N \subset \KK := \{ f \in
  P_{\GG_3}; \,\, M_{m_3}(f) \le C^N_{0,m_3} \}$ thanks to
  assumption {\bf (A1)},
\item $S^{N \! L}_t$ satisfies some H\"older like estimate uniformly
  on $\KK$ and $[0,T)$ thanks to assumption {\bf (A5)},
\item $R^\ell_\varphi \in C^{0,1}(P_{\GG_3},\R)$ because $\varphi \in
  \FF_3^{\otimes\ell}$.
\end{itemize}
We deduce thanks to the Jensen inequality (for a concave function)
\begin{eqnarray*}
  \TT_3 
  &=& \left| \left\langle f^N_0, 
      R^\ell_\varphi \left(S^{N\!L}_t (\mu^N_V)\right) 
    \right\rangle -  \left \langle f_0 ^N, R^\ell_\varphi
      \left(S^{N\!L}_t(f_0 )\right) \right \rangle \right| 
  \\
  &=& \left| \left\langle f^N_0, R^\ell_\varphi \left(S^{N\!L}_t (\mu^N_V)\right) 
      - R^\ell_\varphi \left(S^{N\!L}_t(f_0 )\right) \right\rangle \right| 
  \\
  &\le&  \big[R_\varphi\big]_{C^{0,1}(P_{\GG_3})} \, \left\langle f^N_0, \, 
   \mbox{dist}_{\GG_3} 
    \left( S^{N \! L}_t(f_0), S^{N \! L}_t( \mu^N_V) \right) \right\rangle 
 \\
  &\le&  \big[R_\varphi\big]_{C^{0,1}(P_{\GG_3})} \, \left\langle f^N_0, \, 
  \Theta_{a,T} \left(  \mbox{dist}_{\GG_3} (f_0 , \mu^N_V) \right) \right\rangle 
 \\
  &\le&  \big[R_\varphi\big]_{C^{0,1}(P_{\GG_3})} \,   \Theta_{a,T}  \left( \left\langle f^N_0, \, 
  \mbox{dist}_{\GG_3} (f_0 , \mu^N_V)  \right\rangle \right).
\end{eqnarray*}
Now, by definition of the optimal transport Wasserstein distance we
have
$$
\forall \, \mu_1,\,\mu_2 \in P\left(P_{\GG_3}\right), \quad
\WW_{1,P_{\GG_3}} \left(\mu_1,\mu_2\right) = \inf_{\pi \in
  \Pi \left(\mu_1, \mu_2\right)} \int_{P_{\GG_3} \times P_{\GG_3}}
\hbox{dist}_{\GG_3} (\mu,\mu') \, \pi(d\mu,d\mu'),
$$
where $\Pi(\mu_1,\mu_2)$ denotes the probabilty measures on the
product space $P_{\GG_3} \times P_{\GG_3}$ with first and second
marginals $\mu_1,\mu_2$. In the case when $\mu_2 = \delta_{f_0}$ then
$\Pi (\mu_1 , \delta_{f_0}) = \{ \mu_1 \otimes \delta_{f_0} \}$ has
only one element, and therefore
\begin{eqnarray*}
  \WW_{1,P_{\GG_3}}  \left(\pi^N_P f_0^N , \delta_{f_0}\right) 
  &=& \inf_{\pi \in \Pi (\pi^N_P f_0^N , \delta_{f_0})} 
  \int_{P_{\GG_3} \times P_{\GG_3}} 
  \hbox{dist}_{\GG_3} (f,g) \, \pi(df,dg) \\
  &=&  \int \!\! \int_{P_{\GG_3} \times P_{\GG_3}}  
  \hbox{dist}_{\GG_3} (f,g) \, \pi^N_P f_0^N (dg)
  \, \delta_{f_0}(df) \\
  &=&   \int_{P_{\GG_3}}  \hbox{dist}_{\GG_3} (f_0,g) \, 
  \pi^N_P f_0^N (dg)  \\
  &=&  \int_{E^{N}}  \hbox{dist}_{\GG_3} (\mu^N_V,f_0)   \, f_0^N (dV).
\end{eqnarray*}
We then easily conclude. \qed



\section{(True) maxwellian molecules}
\label{sec:BddBoltzmann}
\setcounter{equation}{0}
\setcounter{theo}{0}


\subsection{The model}
\label{sec:modelEBbounded}

Let us consider $E = \R^d$, $d \ge 2$, and a $N$-particles system
undergoing space homogeneous random Boltzmann type collisions
according to a collision kernel $B = \Gamma(z) \, b (\cos \theta)$
(see Subsection~\ref{sec:introEB}).  More precisely, given a
pre-collisional system of velocity variables $V = (v_1, ..., v_N)
\in E^N = (\R^d)^N$, the stochastic process is:
\begin{itemize}
\item[(i)] for any $i'\neq j'$, draw 
  a random time $ T_{\Gamma(|v_{i'}- v_{j'}|)}$ of collision
  accordingly to an exponential law of parameter $\Gamma(|v_{i'}-
  v_{j'}|)$, and then choose the collision time $T_1$ and the
  colliding couple $(v_i,v_j)$ (which is a.s. well-defined) in such a
  way that
  $$
  T_1 = T_{\Gamma(|v_i - v_j|)} := \min_{1 \le i' \neq j' \le N} T_{\Gamma(|v_{i'}-
    v_{j'}|)};
  $$
\item[(ii)] then draw $\sigma \in S^{d-1}$ according to the law
  $b(\cos \theta_{ij})$, where 
  $\cos \theta_{ij} = \sigma \cdot (v_j-v_i)/|v_j-v_i|$;
\item[(iii)] the new state after collision at time $T_1$ becomes
  $$
  V^*_{ij} = (v_1, \dots, v^*_i, \dots, v^*_j, \dots , v_N),
  $$
  where only velocities labelled $i$ and $j$ have changed, according
  to the rotation
  \begin{equation}\label{vprimvprim*}
    \quad\quad   v^*_i = {v_i + v_j \over 2} 
    + {|v_i - v_j| \, \sigma \over 2}, \quad
    v^*_j= {v_i + v_j \over 2} - {|v_i - v_j| \, \sigma \over 2}.
  \end{equation}
\end{itemize}

\smallskip The associated Markov process $(\VV_t)$ on the velocity
variables on $(\R^d)^N$ is then built by iterating the above
construction. After scaling the time (changing $t \to t/N$ in
order that the number of interactions is of order $\OO(1)$ on finite
time interval, see \cite{spohn}) we denote by $f^N_t$ the law of
$\VV_t$, $S^N_t$ the associated semigroup, $G^N$ and $T^N_t$
respectively the dual generator and dual semigroup, as in the previous
abstract construction. The so-called {\em Master equation} on the law
$f^N_t$ is given in dual form by
\begin{equation}\label{eq:BoltzBddKolmo}
  \partial_t \langle f^N_t,\varphi \rangle = 
  \langle f^N_t, G^N \varphi \rangle 
\end{equation}
with 
\begin{equation}\label{defBoltzBddGN}
  (G^N\varphi) (V) = {1 \over N} \, 
  \sum_{i,j= 1}^N \Gamma\left(|v_i-v_j|\right)
  \,   
  \int_{\mathbb{S}^{d-1}} b(\cos\theta_{ij}) \, \left[\varphi^*_{ij} -
    \varphi\right] \, d\sigma
\end{equation}
where $\varphi^*_{ij}= \varphi(V^*_{ij})$ and $\varphi = \varphi(V)
\in C_b(\R^{Nd})$.

This collision process is invariant under velocities permutations and
satisfies the microscopic conservations of momentum and energy at any
collision time
$$
\sum_{j=1} ^N v^*_j = \sum_{j=1} ^N v_j \quad \mbox{ and } \quad
|V^*|^2 = \sum_{j=1} ^N |v^* _j|^2 = \sum_{j=1} ^N |v_j|^2 = |V|^2.
$$
As a consequence, for any symmetric initial law $f_0^N \in
P_{\mbox{{\scriptsize sym}}}(\R^{Nd})$ the law $f_t^N$ at later times
is also a symmetric probability, and it conserves momentum and energy:
\[
 \forall \, \alpha = 1, \dots, d, \quad 
  \int_{\R^{dN}} \left( \sum_{j=1}^N v_{j,\alpha} \right) \, f_t^N (dV) =
  \int_{\R^{dN}} \left( \sum_{j=1}^N v_{j,\alpha} \right)  \, f_0^N
  (dV),
\]
where $(v_{j,\alpha})_{1 \le \alpha \le d}$ denote the components of
$v_j \in \R^d$, and
\begin{equation}\label{eq:preconservationE}
  \forall \, \phi : \R_+ \to \R_+, \quad  
  \int_{\R^{dN}}\phi( |V|^2 ) \, f_t^N (dV) = 
  \int_{\R^{dN}} \phi( |V|^2 )  \, f_0^N (dV)
\end{equation}
(equality between possibly infinite non-negative quantities).

\medskip The (expected) limiting nonlinear homogeneous Boltzmann
equation is defined by \eqref{el}, \eqref{eq:collop},
\eqref{eq:rel:vit}. The equation generates a nonlinear semigroup
$S^{N\! L}_t (f_0) := f_t$ for any $f_0 \in P_2(\R^d)$ (probabilities
with bounded second moment): for the Maxwell case we refer to
\cite{T1,TV}, for the hard spheres case we refer to \cite{MW99} in a
$L^1$ setting, and for the generalization of these $L^1$ solutions to
$P_2(\R^d)$ we refer to \cite{EM}, \cite{Fo-Mo} and
\cite{Lu-Mouhot}. For these solutions, one has the
conservation of momentum and energy
$$
\forall \, t \ge 0, \quad \int_{\R^{d}} v \, f_t (dv) = \int_{\R^{d}}v
\, f_0 (dv), \qquad \int_{\R^{d}} |v|^2 \, f_t (dv) =
\int_{\R^{d}}|v|^2 \, f_0 (dv).
$$

Without restriction, by using the change of variable $\sigma \mapsto -
\sigma$, from now on we restrict the angular domain to $\theta \in
[-\pi/2,\pi/2]$ for the limiting equation as well as the $N$-particle
system. Therefore we assume $\mbox{supp} \, b \subset [0,1]$. We still
denote by $b$ the symmetrized version of $b$ by a slight abuse of
notation.

\subsection{Statement of the result}
\label{sec:resultEBbounded} In this section we consider the case of
the {\sl Maxwell molecules kernel}. More precisely we shall assume
\begin{equation}\label{Maxwelltrue}
  \Gamma \equiv 1, \quad b \in
  L^\infty_{loc}([0,1)), 
 \quad \forall \, \eta > 1/2, \quad  
C_\eta(b) := \int_{\mathbb{S}^{d-1}}  b(\cos \theta)\, |1-\cos \theta|^{\eta/2} \, 
  d\sigma < \infty. 
\end{equation}
Indeed for any positive real function $\psi$ and any given vector $u
\in \R^d$ we have
\[
 \int_{S^{d-1}} \psi (\hat
u \cdot \sigma) \, d\sigma
= |\mathbb{S}^{d-2}| \, \int_0^\pi \psi (\cos\theta) \, \sin^{d-2} \, \theta \,
d\theta 
\]
and the true Maxwell angular collision kernel {\bf (tMM}) defined in
Subsection~\ref{sec:introEB} satisfies (in dimension $d=3$) 
$b(z) \sim K \, (1-z)^{-5/4}$ as $z \to 1$, which hence fulfills
(\ref{Maxwelltrue}). These assumptions also trivially include the
Grad's cutoff Maxwell molecules {\bf (GMM)} introduced in
Subsection~\ref{sec:introEB}.

Our fluctuations estimate result then states as follows:

\begin{theo}[Boltzmann equation for Grad's cut-off/true Maxwell
  molecules]\label{theo:tMM}
  Let us consider an initial distribution $f_0 \in P(\R^d)$ with
  compact support. 
  Let us consider a hierarchy of $N$-particle distributions $f^N _t =
  S^N _t(f_0 ^{\otimes N})$ issued from the tensorized initial data
  $f^N _0 = f_0 ^{\otimes N}$.
   Then for any
  \[ 
  \varphi = \varphi_1 \otimes \varphi_2 \otimes \dots \otimes \,
  \varphi_\ell \in \FF^{\otimes\ell}, 
  \] 
  with 
  $$
  \FF := \left\{ \varphi : \R^d \to \R; \,\, \| \varphi \|_\FF :=
    \int_{\R^d} (1 + |\xi|^4) \, |\hat \varphi (\xi)|Ê\, d\xi <
    \infty \right\},
  $$
  we have 
  \begin{eqnarray} \label{eq:cvgBddBE} 
  &&\sup_{t \ge 0} \left| \left
        \langle \left( S^N_t(f_0 ^N) - \left( S^{N\! L}_t(f_0)
          \right)^{\otimes N} \right), \varphi \right\rangle \right|
  \\ \nonumber
  &&\le 
 C \, \Bigg[ \ell^2 \, \frac{\|\varphi\|_\infty}{N} 
  + C^N_{T,4} \, {C_{\eta,\infty}^\infty \over N^{1-\eta}} \, \ell^2 \, 
  \|\varphi\|_{\FF ^2 \otimes (L^\infty)^{\ell-2}} 
  \\ \nonumber
  &&\qquad\qquad\qquad 
  + \ell \, \, \|\varphi\|_{W^{1,\infty} \otimes (L^\infty)^{\ell-1}} \, 
  \WW^N_{W_2} ( f_0)   \Bigg],
  \\ \nonumber
  &&
    \le \ell^2 \, {C_\eta \over N^{{2 \over d} - \eta}} \,
    \| \varphi \|_{\FF^{\otimes\ell}}
    \end{eqnarray} 
    for some constants $C_{\eta,\infty}^\infty, \, C_\eta,  \in (0,\infty)$ which may blow up when
    $\eta \to 0$ (and depend on $b$ and $\EE := M_2(f_0)$).
\end{theo}




In order to prove Theorem~\ref{theo:tMM}, we have to establish
the assumptions {\bf (A1)-(A2)-(A3)-(A4)-(A5)} of
Theorem~\ref{theo:abstract} with $T=\infty$.

\subsection{Proof of (A1)} \label{sec:MaxA1} The operators $S^N_t$,
$T^N_t$ and $G^N$ are well defined on $L^2(\mathbb{S}^{d N
  -1}(\sqrt{\EE}))$ for any $\EE > 0$ thanks to the facts 
\[
\left\langle G^N \varphi, \varphi \right\rangle_{L^2(\mathbb{S}^{d N
  -1}(\sqrt{\EE}))} = - {1 \over N} \, \sum_{i,j= 1}^N \Gamma\left(|v_i-v_j|\right)
  \,   
  \int_{\mathbb{S}^{d-1}} b(\cos\theta_{ij}) \, \left[\varphi^*_{ij} -
    \varphi\right]^2 \, d\sigma
\]
and $(S^N_t)^* = S^N _t = T^N _t$ on this space, for any $\EE >
0$ (see for instance~\cite[Chapter~9, section~2]{Kato}).

Then $G^N$ is well defined on the domain $C^1 _b (\R^{dN})$, and
it is closable on this space by using the fact that for any $\varphi_n
\in \mbox{Dom}(G^N)$, $\varphi_n \to 0$, $G^N \varphi_n \to \psi$, one
has for any $\EE >0$ that the restrictions of $\varphi_n$, $G^N
\varphi_n$ to the sphere of energy $\EE$ belongs to $L^2$ on this
compact space, and then, using the preceeding discussion, we deduce 
that $\psi \equiv 0$ on this subspace. Since we can argue in this way
for any $\EE >0$, we deduce that $\psi \equiv 0$, which is the
definition of being closable. 

Then it remains to prove bounds on the polynomial moments of the
$N$-particles system. We shall prove the following more general lemma:

\begin{lem}\label{lem:momentsN}
  Consider the collision kernel $B = |v-v_*|^\gamma \, b(\theta)$ with
  $\gamma = 0$ or $1$ and $b \in L^1([0,1), (1-z)^2)$. This covers the
  three cases {\bf (HS)}, {\bf (tMM)} and {\bf (GMM)}.

Assume that the initial datum of the $N$-particle system satisfies:
\[
\hbox{{\em supp}} \, f^N_0 \subset \left\{V \in \R^{Nd}; \,\, M^N_2(V)
  \le \EE \right\}, \quad M_2 ^N = \frac1N \, \sum_{j=1} ^N |v_j|^2
\]
and
\[
\left \langle f^N _0, M_k ^N \right \rangle \le C_{0,k} <\infty, \quad 
M_k ^N = \frac1N \, \sum_{j=1} ^N |v_j|^k, \quad k \ge 2.
\]
Then we have
\[
\sup_{t \ge 0} \left \langle f^N _t, M_k ^N \right\rangle \le 
\max \left\{ C_{0,k}; \, \bar a_k \right\} 
\]
where $\bar a_k \in (0,\infty)$ depends on $k$ and $\EE$.
\end{lem}

\noindent
{\it Proof of Lemma~\ref{lem:momentsN}.}  
From (\ref{eq:pNUnifEnergy}) we have
\begin{equation}\label{eq:MaxwellSupportE}
\hbox{supp} \, f^N_t \subset 
\left\{V \in \R^{Nd}; \,\, M^N_2(V) \le \EE \right\},
\end{equation}
by taking $\phi (z) := {\bf 1}_{z > N \, \EE}$ in
\eqref{eq:preconservationE}.

Next, we write the differential equality on $k$-th moment 
$$
\frac{d}{dt} \left\langle f^N _t , \frac1N \, \sum_{j=1} ^N
  |v_j|^{k}\right\rangle = \frac{1}{N^2} \, \sum_{j_1 \neq j_2} ^N
\left\langle f^N _t, \left|v_{j_1}-v_{j_2}\right|^\gamma \,
  \KK\left(v_{j_1},v_{j_2}\right) \right\rangle,
$$
with 
$$
\KK\left(v_{j_1},v_{j_2}\right) = {1 \over 2} 
\int_{\mathbb{S}^{d-1}} b(\theta_{j_1 j_2}) \, \left[ |v_{j_1} ^*|^{k}
  + |v_{j_2} ^*|^{k} - |v_{j_1}|^{k} - |v_{j_2}|^{k} \right] \, d\sigma.
$$
From the so-called Povner's Lemma proved in \cite[Lemma 2.2]{MW99}
(valid for singular collision kernel as in our case), we have
$$
\KK (v_{j_1},v_{j_2}) \le C_1 \, (|v_{j_1}|^{k-1} \, |v_{j_2}| +
|v_{j_1}| \, |v_{j_2}|^{k-1}) - C_1 \, (|v_{j_1}|^{k} + |v_{j_2}|^{k})
$$
for
some constants $C_1, C_2 \in (0,\infty)$ depending only on $k$ and
$b$.

By using the inequalities $|v_{j_1}-v_{j_2}| \ge |v_{j_1}| -
|v_{j_2}|$ and $|v_{j_1}-v_{j_2}| \ge |v_{j_2}| - |v_{j_1}|$ in order
to estimate the last term when $\gamma =1$, we then deduce
\begin{multline*}
  |v_{j_1}-v_{j_2}| \, \KK (v_{j_1},v_{j_2}) \le \\
  \le C_3 \, [ (1+|v_{j_1}|^{k+\gamma-1}) \, (1+ |v_{j_2}|^2) + (1+ |v_{j_1}|^2)
  \, (1+|v_{j_2}|^{k+\gamma-1})] - C_1 \, (|v_{j_1}|^{k+\gamma} +
  |v_{j_2}|^{k+\gamma}),
\end{multline*}
for a constant $C_3$ depending on $C_1$ and $C_2$. 

Using (symmetry hypothesis) that
$$
\forall \, k \ge 0, \quad \langle f^N_t , |v_1|^k \rangle = \left
  \langle f^N_t, M^N_k \right \rangle,
$$
and  \eqref{eq:MaxwellSupportE}  we get 
\begin{multline*} 
  {d \over dt} \left\langle f^N_t , |v_1|^{k} \right\rangle \le 2 \,
  C_3 \left\langle f^N_{t} , (1 + M^N_{k+\gamma-1}) \, (1 + M^N_{2})
  \right\rangle - 2 \, C_2 \left\langle f^N_{t} , M_{k+\gamma} \right\rangle \\
  \le C_1 \, (1+\EE) \, \left(1 + \left\langle f^N_{t} ,
      |v_1|^{k+\gamma-1} \right\rangle \right) - 2 \, C_2 \left\langle
    f^N_{t} , |v_1|^{k+\gamma} \right\rangle.
\end{multline*}
Using finally H\"older's inequality 
\[
\left\langle f^N_1,|v|^{k-\gamma+1} \right\rangle \le \left\langle
f^N_1,|v|^{k+\gamma} \right\rangle^{(k-\gamma+1)/(k+\gamma)}
\]
we conclude that $y(t) = \langle f^N_t , |v_1|^{k} \rangle$ satisfies a
differential inequality of the following kind
$$
y' \le - K_1 \, y^{\theta_1} +  K_2 \, y^{\theta_2} + K_3
$$
with $\theta_1 \ge 1$ and $\theta_2 < \theta_1$, which concludes the
proof of the lemma. \qed

\medskip
Lemma~\ref{lem:momentsN} proves {\bf (A1)}-(i) with $m_e (v) = |v|^2$,
 {\bf (A1)}-(ii) with $m_1(v) = |v|^4$ (and we do not need {\bf (A1)}-(iii) in this
 particular case: we may take $m_3 \equiv 0$). 


\subsection{Proof of (A2).}
Let us define
\[
P_{\GG_1} := \left\{ f \in P_4(\R^d) \, ; \ \langle f, m_e
  \rangle \le \EE \right\}
\quad 
\mbox{endowed with the distance induced by} \  | \cdot |_2.
\] 

Let us prove {\bf (A2)}-(i)-(ii) and {\bf (A2)}-(iii) with:
\begin{equation}\label{A2trueMaxwell}
\forall \, f_0 \in P_{\GG_1}, \quad  
\forall \, t \in (0,1] \qquad \left|
  S^{N\!L}_t f_0 - f_0 - t \, Q(f_0,f_0) \right|_2 \le C \, t^2,
\end{equation}
(one can prove more generally that the application $\SS(f_0) :
[0,\tau) \to P_{\GG_1}$, $t \mapsto \SS(f_0)(t) := S^{N \!
  L}_t(f_0)$ is $C^{1,1}([0,\tau);P_{\GG_1})$, with
$\SS(f_0)'(0) = Q(f_0,f_0)$ and $\tau >0$).


 

\medskip
Let us recall the following result proved in \cite{TV}.  
We provide its proof for the the sake of completeness and because we will
need to modify it in order to obtain similar results in the next sections.

\begin{lem} \label{lem:contraction}
For any $f_0,g_0 \in P_2(\R^d)$, the associated solutions $f_t$ and
$g_t$ to the Boltzmann equation for Maxwellian collision kernel satisfy
\begin{equation}\label{estim:C01TrueMax}
        \sup_{t \ge 0} \left|f_t-g_t\right|_2 \le 
        \left|f_0-g_0\right|_2.
\end{equation}

\end{lem}

\noindent
{\it Proof of Lemma~\ref{lem:contraction}.} 
We recall Bobylev's identity for Maxwellian collision kernel
({\it cf.} \cite{Bobylev-88})
$$
\FF\left(Q^+(f,g)\right) (\xi) = \hat Q^+ (F,G) (\xi) =: {1 \over 2}
\int_{S^{d-1}} b\left(\sigma \cdot \hat\xi\right) \, [F^+ \, G^- + F^- \, G^+ ]\,
d\sigma,
$$
with $F = \hat f$, $G = \hat g$, $F^\pm= F(\xi^\pm)$, $G^\pm=
G(\xi^\pm)$, $\hat \xi = \xi / |\xi|$ and
$$
\xi^+ = {1\over 2} (\xi +  |\xi| \, \sigma),
 \quad
\xi^- = {1\over 2} (\xi -  |\xi| \, \sigma).
$$
Denoting by $D = \hat g - \hat f$, $S = \hat g + \hat f$, the
following equation holds
\begin{equation}\label{eq:BoltzMaxD}
  \partial_t D = \hat Q (S,D) =  \int_{S^2} b
  \left(\sigma \cdot \hat \xi \right) \, 
  \left[ \frac{D^+ \, S^-}{2} + \frac{D^- \, S^+}{2} - D \right] \, d\sigma.
\end{equation}
We perform the following cutoff on the angular collision kernel:
\[
\int_{\mathbb{S}^{d-1}} b_K \left( \sigma \cdot \hat \xi \right) \,
d\sigma = K, \quad b_K = b \, {\bf 1}_{|\theta|\ge \delta(K)}
\]
for some well-chosen $\delta(K)$, and we shall relax this assumption
in the end (using uniqueness of measure solutions of \cite{TV}).
Using that $\| S \|_\infty \le 2$, we deduce in distributional sense
$$
\frac{d}{dt} {|D| \over |\xi|^2} + K \, {|D| \over |\xi|^2} \le 
\left( \sup_{\xi \in \R^d} {|D| \over |\xi|^2} \right) \, \left( \sup_{\xi \in \R^d}
\int_{S^{d-1}} b_K\left(\sigma \cdot \hat \xi \right) \, \left(
  \left|\hat\xi^+\right|^2 + \left|\hat\xi^-\right|^2 \right) \,
d\sigma \right)
$$
with
$$
\left|\hat\xi^+\right| = 
{1 \over \sqrt{2}} \, \left(1 + \sigma \cdot \hat\xi\right)^{1/2},
\qquad \left|\hat\xi^-\right| = {1 \over \sqrt{2}} \, \left(1 - \sigma \cdot
\hat\xi\right)^{1/2}.
$$
By using $|\hat\xi^+|^2 + |\hat\xi^-|^2 =1$, we deduce 
\[
\frac{d}{dt} {|D| \over |\xi|^2} + K \, {|D| \over |\xi|^2} \le 
K \, \left( \sup_{\xi \in \R^d} {|D| \over |\xi|^2} \right)
\]
from which we deduce 
\[
\left( \sup_{\xi \in \R^d} \frac{|D_t(\xi)|}{|\xi|^2} \right) \le 
\left( \sup_{\xi \in \R^d} \frac{|D_0(\xi)|}{|\xi|^2} \right) 
\]
for any value of the cutoff parameter $K$. Therefore by relaxing the
cutoff $K \to \infty$, we deduce \eqref{estim:C01TrueMax}.
\qed  


 
 \medskip
Hence we deduce that $S^{N\!L}_t$ is $C^{0,1}(P_{\GG_1},P_{\GG_1})$ and {\bf (A2)-(i)} is proved. 

\begin{lem}
  \label{lem:estimQtoscani1}
  For any $f,g \in P_2$ with same momentum and finite second moment, we have
  \begin{equation}\label{eq:Qtoscani1}
    \left| Q(f,f) \right|_1 \le C \, \left( \int_{\R^d} (1+|v|) \,
      df(v) \right)^2
  \end{equation}
  \begin{equation}\label{eq:Qtoscani3}
    \left| Q(f,f) \, v \right|_1 \le C \, \left( \int_{\R^d} (1+|v|) \,
      df(v) \right)^2
  \end{equation}
  \begin{equation}\label{eq:Qtoscani2}
    \left| Q(f,f) \right|_2 \le C \, \left( \int_{\R^d} \left(1+|v|^2\right) \,
      df(v) \right)^2
  \end{equation}
  and 
  \begin{equation}\label{eq:Qbiltoscani1}
    \left|  Q(f+g,f-g) \right|_1 \le C \, \left( \int_{\R^d} (1+|v|) \,
      \left( df(v) + dg(v) \right) \right) \, 
    \left( | f-g |_1 + \int_{\R^d} (1+|v|) \, |df-dg|(v) \right).
  \end{equation}
  \begin{equation}\label{eq:Qbiltoscani3}
    \left| Q(f+g,f-g) \, v  \right|_1 \le C \, \left( \int_{\R^d} (1+|v|) \,
      \left( df(v) + dg(v) \right) \right) \, 
    \Big( | f-g |_1 + \big|  (f-g) \, v \big|_1 \Big).
  \end{equation}
   \begin{equation}\label{eq:Qbiltoscani2}
    \left| Q(f+g,f-g) \right|_2 \le C \, \left( \int_{\R^d} (1+|v|) \,
      \left( df(v) + dg(v) \right) \right) \, 
    \Big( | f-g |_2 + \big|  (f-g) \, v  \big|_1 \Big).
  \end{equation}
\end{lem}

\noindent
{\it Proof of Lemma~\ref{lem:contraction}.} 
We prove the second inequalities
\eqref{eq:Qbiltoscani1}-\eqref{eq:Qbiltoscani2}. The first
inequalities \eqref{eq:Qtoscani1}-\eqref{eq:Qtoscani2} are then a
trivial consequence by using
\[
Q(f,f) = Q(f,f) - Q(M,M) = Q(f-M,f+M)
\]
where $M$ is the maxwellian distribution with same momentum and energy
as $f$, and then applying \eqref{eq:Qbiltoscani1} or
\eqref{eq:Qbiltoscani2} with $f-M$ and $f+M$.

We write in Fourier:
\[
\mathcal{F}\left( Q(f+g,f-g) \right) = \hat Q(D,S) = \frac12 \,
\int_{\mathbb{S}^{d-1}} b(\sigma \cdot \hat \xi) \, \left( S(\xi^+) \,
  D(\xi^-) + S(\xi^-) \, D(\xi^+) - 2 \, D(\xi) \right)
\]
where $\hat Q$ is the Fourier form the symmetrized collision operator
$Q$, which yields
\[
\frac{\left| \hat Q(D,S) \right|}{|\xi|^2} \le \TT_1 + \TT_2 + \TT_3
\]
with 
\[
\TT_1 \le \int_{\mathbb{S}^{d-1}} b(\sigma \cdot \hat \xi) \, 
\left|S(\xi^+)\right| \, \frac{\left| D(\xi^-) \right|}{|\xi^-|^2} \,
\frac{\left| \xi^- \right|^2}{|\xi|^2} \, d\sigma 
\le C \, | D |_2
\]
and 
\[
\TT_2 \le \int_{\mathbb{S}^{d-1}} b(\sigma \cdot \hat \xi) \,
\frac{\left|D(\xi^+)\right|}{|\xi^+|} \, \frac{\left| S(\xi^-) -2
  \right|}{|\xi^-|} \, \frac{\left| \xi^- \right|}{|\xi|} \, d\sigma
\le C \, | D |_1 \, \left( \int_{\R^d} (1+|v|) \, ( df(v) + dg(v) )
\right) 
 \] 
and 
\begin{multline*}
\TT_3 \le 2 \, \int_{\mathbb{S}^{d-1}} b(\sigma \cdot \hat \xi) \,
\frac{\left| D(\xi^+) - D(\xi) \right|}{|\xi|} \, d\sigma \\
\le \int_{\mathbb{S}^{d-1}} b(\sigma \cdot \hat \xi) \,
\frac{|\xi^-|}{|\xi|} \, \int_0 ^1 \frac{|\nabla D (\theta \xi +
  (1-\theta) \xi^+)|}{|\theta \xi +
  (1-\theta) \xi^+|} \, d\theta \, d\sigma 
\le C \, | (f - g)  \, v |_1
\end{multline*}
which concludes the proof of  \eqref{eq:Qbiltoscani1}. The proof of \eqref{eq:Qbiltoscani3} and  \eqref{eq:Qbiltoscani2} are similar. 
\qed

\medskip

The proof of \textbf{(A2)-(ii)} is a consequence of \eqref{eq:Qbiltoscani2}, while the 
proof of \textbf{(A2)-(iii)}  can be done by repeated use of
Lemma~\ref{lem:estimQtoscani1}. Indeed, we start from
\[
f_t - f_0 = \int_0 ^t Q(f,f) \, ds
\]
from which we deduce thanks to \eqref{eq:Qtoscani1}, \eqref{eq:Qtoscani3}, \eqref{eq:Qtoscani2} 
\[
\left| f_t - f_0 \right|_1 + \left| (f_t - f_0) \, v \right|_1 + \left| f_t - f_0 \right|_2 \le C \, t. 
\]
Then write
\[
f_t - f_0 - t \, Q(f_0,f_0) = \int_0 ^t \left( Q(f_s,f_s) - Q(f_0,f_0)
\right) \, ds = \int_0 ^t Q(f_s-f_0,f_s+f_0) \, ds, 
\]
from which we deduce thanks to \eqref{eq:Qtoscani3}, \eqref{eq:Qtoscani2} 
\[
\left| f_t - f_0 - t \, Q(f_0,f_0) \right|_2 \le C \, \int_0 ^t \left(
  \left| f_s - f_0 \right|_2 + \left| v \, (f_s - f_0) \right|_1
\right) \, ds \le C \, \int_0 ^t s \, ds \le C \, t^2.
\]

\subsection{Proof of  (A3) with $\GG_1$. } 
Define $\Lambda_1(\rho) := \langle \rho , 1 + |v|^4 \rangle^{1-\eta}$ for any $\rho \in P_{\GG_1}$ and some $\eta \in (0,1)$. 
Let us prove that for any $\Phi \in C^{1,\eta}_{\Lambda_1}\left(P_{\GG_1,{\bf r}};\R\right)$, ${\bf r}Ê\in {\bf R}_{\GG_1}$, 
\begin{equation}
  \label{eq:H1VMix}
   \left\| M^N_4(V)^{\eta-1} \, \left( G^N \, \pi^N - \pi^N \, G^\infty \right)  \,
    \Phi \right\|_{L^\infty(\mathbb{E}_N)} \le
  {C_{1} \, \EE \over N} \, [ \Phi ]_{C^{1,\eta}_{\Lambda_1} \left(P_{\GG_1,{\bf r}};\R\right)},
\end{equation}
for some constant $C_1$ and where $\mathbb{E}_N = \{ V \in (\R^d)^N \ , \
|V|^2 \le \EE \}$.

First, consider velocities $v,v_*, w,w_* \in \R^d$ such that
$$
w = {v+v_* \over 2} + {|v-v_*| \over 2} \, \sigma, \quad w_* = {v+v_*
  \over 2} - {|v-v_*| \over 2} \, \sigma, \quad \sigma \in \mathbb{S}^{d-1}.
$$
Then $\delta_v + \delta_{v_*} - \delta_w - \delta_{w_*} \in \mathcal{I} P_{\GG_1}$.
Performing a Taylor expansion up to order two and one, we have 
\begin{multline*}
e^{i v \cdot \xi} + e^{i v_* \cdot
          \xi} - e^{i w \cdot \xi} - e^{i w_* \cdot \xi} = \\
 = i \, (w-v) \xi \, e^{i v \cdot \xi} + \OO(|w-v|^2 \, |\xi|^2) +  i \, (w_*-v_*) \xi \, e^{i v_* \cdot \xi} + \OO(|w_*-v_*|^2 \, |\xi|^2) 
 \\
 = i \, (w-v) \, e^{i v \cdot \xi} + \OO(|w-v|^2 \, |\xi|^2) +  i \, (w_*-v_*) \, (e^{i v_ \cdot \xi} + \OO ( |v-v_*|Ê\, |\xi|) + \OO(|w_*-v_*|^2 \, |\xi|^2) 
 \\
 = \OO ( |\xi|^2 \, |v-v_*|^2  \cos \theta) 
 \end{multline*}   
thanks to the impulsion conservation and the fact that 
 $$
 |w-v| + |w_* - v_*| \le  \sqrt{2} \, |v-v_*| \, (1-\cos \theta).
 $$
We then deduce 
\begin{multline*}
    \left| \delta_v + \delta_{v_*}- \delta_w - \delta_{w_*} \right|_2 
    = \sup_{\xi \in \R^d} \frac{ \left| e^{i v \cdot \xi} + e^{i v_* \cdot \xi} 
          - e^{i w \cdot \xi} - e^{i w_* \cdot  \xi}\right| }{|\xi|^2}
          \le  C \, |v-v_*| \, (1-\cos \theta).
\end{multline*}
Consider $V \in \mathbb{E}_N$ and define 
$$
{\bf r}_V := \left(\langle \mu^N_V,z_1 \rangle, \dots, \langle
  \mu^N_V,z_d \rangle, \langle \mu^N_V, |z|^2 \rangle\right) \in {\bf
  R}_\EE.
$$
Then for a given $\Phi \in C_\Lambda^{1,\eta}(P_{\GG_1,{\bf r}_V};\R)$,
we set $\phi := D\Phi[ \mu_V ^N]$, $u_{ij} = (v_i-v_j)$ and we compute:
\begin{eqnarray*}
&&  G^N(\Phi \circ  \mu ^N _V) 
=  \frac{1}{2N} \, \sum_{i,j=1} ^N 
  \int_{\mathbb{S}^{d-1}} \left[ \Phi\left(\mu^N_{V^*_{ij}}\right) 
    - \Phi\left(\mu^N _V\right) \right] \, b(\theta_{ij}) \, d\sigma 
\\
&&\qquad= \frac{1}{2N} \, \sum_{i,j=1} ^N
  \int_{\mathbb{S}^{d-1}} \left\langle \mu^N _{V^*_{ij}}
    - \mu^N _V , \phi \right\rangle \, b(\theta_{ij}) \, d\sigma 
\\
&&\qquad+  \frac{[\Phi]_{ C_\Lambda^{1,\eta}(P_{\GG_1,{\bf r}_V}) }}{2N} \, 
  \sum_{i,j=1} ^N  \int_{\mathbb{S}^{d-1}} 
  \!\!\! \ \left[ M_{4}\left(\mu^N _{V^*_{ij}} +  \mu^N _{V}\right) \! \right]^{1-\eta} \, 
  \OO \! \left( \left| \mu^N _{V^*_{ij}} 
      - \mu^N _V \right|_2^{1+\eta} \right)  \, d\sigma \\
&&\qquad =:  I_1(V) + I _2(V).
\end{eqnarray*}
For the first term $I_1(V)$, thanks to Lemma~\ref{lem:H0}, we have 
\begin{multline*}
  I_1(V) = \frac{1}{2N^2} \, \sum_{i,j=1} ^N
  \int_{\mathbb{S}^{d-1}} b(\theta_{ij}) \, \left[ \phi(V_i ^*) +
    \phi(V_j^*) - \phi(V_i) - \phi(V_j) \right] \, d\sigma \\
  = \frac{1}{2N^2} \, \int_v \int_w 
  \int_{\mathbb{S}^{d-1}} b(\theta) \, \left[ \phi(v^*) +
    \phi(w^*) - \phi(v) - \phi(w) \right] \, \mu^N _V(dv) \, 
  \mu^N_V (dw) \, d\sigma \\
  = \left\langle Q(\mu^N _V, \mu^N _V), \phi \right\rangle = 
  \left(G^\infty\Phi\right)(\mu^N _V).
\end{multline*}
For the second term $I_2(V)$, using that 
\begin{eqnarray*}
  M_{4}\left(\mu^N_{V^*_{ij}}\right) \le C \, M^N_{4}(V),
\end{eqnarray*}
we deduce 
\begin{eqnarray*}
  |I_2(V)|  
  &\le& \frac{C}{N^{2+\eta}} \, [M^N_{4}(V)]^{1-\eta} \,   [\Phi]_{C^{1,\eta}_\Lambda(P_{\GG_1,{\bf r}_V})}
  \, \sum_{i,j= 1}^N \OO \left\{ \int_{\mathbb{S}^{d-1}}
    b \left(\cos \theta_{ij}\right)  \, 
    \left(1+|v_i|^2+|v_j|^2\right) \, 
    \left(1 - \sigma \cdot \hat u_{ij}\right) \, d\sigma  \right\}  
\\
  & \le & C \, {  [M^N _{4}(V) ]^{1-\eta}} \, 
  [\Phi]_{C^{1,\eta}_\Lambda(P_{\GG_1,{\bf r}_V})}  \,  {1 \over N^\eta} \, (1+\EE),
 \end{eqnarray*}
which concludes the proof.

\subsection{Proof of (A4) in Fourier-based norms $|\cdot|_2$ and
  $|\cdot|_4$}
 
For $f_0, g_0 \in P_4(\R^d)$, let us define the associated solutions $f_t$ 
and $g_t$ to the nonlinear Boltzmann equation as well as $h_t := \LL^{N\!  L}_t[f_0](g_0-f_0)$
 the solution to the  linearized Boltzmann equation around $f_t$. More precisely, 
 we define
\[
\left\{
\begin{array}{l}
\partial_t f_t = Q(f_t,f_t), \qquad f_{|t=0} = f_0 \vspace{0.3cm} \\
\partial_t g_t = Q(g_t,g_t), \qquad g_{|t=0} = g_0 \vspace{0.3cm} \\
\partial_t h_t = Q(f_t,h_t) + Q(h_t,f_t), \qquad h_{|t=0} = h_0 := g_0 - f_0.
\end{array}
\right.
\]
 
\begin{lem}\label{lem:a4maxwellfourier} There exists $\lambda \in (0,\infty)$ 
and for any $\eta \in (0,1)$ there exists $C_\eta$ such that for any 
$f_0, g_0 \in P_{\GG_1,{\bf r}}$, ${\bf r}Ê\in {\bf R_{\GG_1}}$, we have
\begin{eqnarray}\label{ineq;tMMC01}
  \left| f_t - g_t \right|_2 \le C_\eta \, e^{- (1-\eta) \, \lambda \, t } \, \max ( M_4(f_0), M_4(g_0))^{(1-\eta)/2} \, \left| f_0 - g_0 \right|_2^\eta, 
\end{eqnarray}
and 
 \begin{eqnarray}\label{ineq;tMMC10}
   \left| h_t \right|_2 \le C_\eta \, e^{- (1-\eta) \, \lambda \, t } \, \max ( M_4(f_0), M_4(g_0))^{(1-\eta)/2} \, \left| f_0 - g_0 \right|_2^\eta. 
\end{eqnarray}
As a consequence, the operator $\LL^{N\!  L}_t$ defined by  $\LL^{N\!  L}_t[f_0](h_0) := h_t$ is such that   $\LL^{N\!  L}_t[f_0] \in \MM^1(\GG_{1,{\bf r}},\GG_{1,{\bf r}})$. 
\end{lem}

\noindent
{\it Proof of Lemma~\ref{lem:a4maxwellfourier}.} {\sl Step 1. Estimate in   $|\cdot|_4$. } We proceed in the spirit of \cite{T1,CGT}. 
With the notation $\MM = \MM_4$, $\hat \MM = \hat \MM_4$, introduced in Example~\ref{expleFourierGal}, 
let us define $d:=f-g$,  $s:=f+g$ and $\tilde d := d - \MM[d]$, and then $D := \FF(d)$, $S := \FF(s)$ and 
\[
\tilde D := \FF(\tilde d) = D - \hat{\mathcal{M}}[d] . 
\]
The equation satisfies by $\tilde D$ is 
\begin{eqnarray}\label{eq:BoltzTildeD}
  \partial_t \tilde D 
  &=&  \hat Q (D,S) - \partial_t \hat \MM[d]
  \\ \nonumber
  &=&  \hat Q (\tilde D,S) + \{ÊÊ\hat Q(\hat\MM[d],S) - \hat\MM[Q(d,s)] \}.
\end{eqnarray}
We infer from \cite{T1} that for any $j \in \N^d$, there exists some absolute coefficients $(a_{i,j})$ 
depending on the collision kernel $b$ such that 
$$
\int_{\mathbb{S}^{d-1}} b(\cos\theta) \left[(v^i)' + (v^i)'_* - (v^i) -
(v^i)_* \right] \, d\sigma = \sum_{j, \, |j| \le |i|} a_{i,j} \, \left(v^j\right)
\, \left(v^{i-j}\right)_*.
$$
We deduce that 
\[
\forall \, |i|\le 3, \quad \nabla^i _\xi  
\mathcal{\hat M}[Q(d,s)]_{\big|\xi=0} = M_i[Q(d,s)] = \sum_{j, \, |j| \le |i|}
a_{i,j} \, M_j[d] \, M_{i-j}[s]
\]
together with 
\begin{multline*}
\forall \, |i|\le 3, \quad \nabla^i _\xi \hat Q(\mathcal{\hat
  M}[d],S)_{\big|\xi=0} = M_i[Q(\mathcal{M}[d],s)] \\ 
= \sum_{j, \, |j| \le |i|}
a_{i,j} \, M_j[\mathcal{M}[d]] \, M_{i-j}[s] =  \sum_{j, \, |j| \le |i|}
a_{i,j} \, M_j[d] \, M_{i-j}[s]
\end{multline*}
since  $M_i[\mathcal{M}[d]] = M_i[d]$ for any $|i| \le 3$
by construction. As a consequence, we get 
\begin{eqnarray}\label{eq:MMQds-QMMdS}
\forall \, \xi \in \R^d,\qquad 
\left| \hat{\mathcal{M}}[Q(d,s)] -
  \hat Q( \hat{\mathcal{M}}[d], S ) \right| 
\le C \, |\xi|^4 \,
\Bigl( \sum_{|i| \le 3} \left|M_i[f_t-g_t]\right| \Bigr).
\end{eqnarray}
On the other hand, from \cite[Theorem 8.1]{T1} and its corollary, we know that there exists some
constants $C,\lambda \in (0,\infty)$ such that 
\begin{eqnarray}\label{eq:MomentsTanaka}
\forall \, t \ge 0 \qquad 
\Bigl( \sum_{|i| \le 3} \left|M_i[f_t-g_t]\right| \Bigr) \le 
C \, e^{-\lambda \, t} \Bigl( \sum_{|i| \le 3} \left|M_i[f_0-g_0]\right| \Bigr). 
\end{eqnarray}
Gathering \eqref{eq:BoltzTildeD}, \eqref{eq:MMQds-QMMdS} and \eqref{eq:MomentsTanaka} and
performing the same cutoff on the angular collision kernel as in the proof of Lemma~\ref{lem:contraction}, 
we have 
\begin{multline*}
  \frac{d}{dt} {|\tilde D(\xi)| \over |\xi|^4} + K \, {|\tilde D(\xi)|
    \over |\xi|^4} \le \left( \sup_{\xi \in \R^d} {|\tilde D(\xi)|
      \over |\xi|^4} \right) \, \left( \sup_{\xi \in \R^d}
    \int_{S^{d-1}} b_K\left(\sigma \cdot \hat \xi \right) \, \left(
      \left|\hat\xi^+\right|^4 + \left|\hat\xi^-\right|^4 \right) \,
    d\sigma \right) \\
  + C \, e^{-\lambda \, t} \left( \sum_{|i| \le 3} \left|M_i[f_0-g_0]\right| \right).
\end{multline*}
Let us compute (the supremum has been droped by the spherical invariance)
\[
\lambda_K := \int_{S^{d-1}} b_K\left(\sigma \cdot \hat \xi \right) \, \left(
  \left|\hat\xi^+\right|^4 + \left|\hat\xi^-\right|^4 \right) \,
d\sigma = \int_{S^{d-1}} b_K\left(\sigma \cdot \hat \xi \right) \, \left(
 \frac{1+ \left( \sigma \cdot \hat \xi \right)^2}{2} \right) \,
d\sigma,
\]
so that 
\begin{multline*}
  \lambda_K - K = - \left( \int_{S^{d-1}} b_K\left(\sigma \cdot \hat \xi
    \right) \, \left( \frac{1- \left( \sigma \cdot \hat \xi
        \right)^2}{2} \right) \, d\sigma \right) \\
\xrightarrow[K \to
  \infty]{} - \int_{S^{d-1}} b\left(\sigma \cdot \hat \xi \right) \,
  \left( \frac{1- \left( \sigma \cdot \hat \xi \right)^2}{2} \right)
  \, d\sigma := - \bar\lambda \in (-\infty,0). 
\end{multline*}
Thanks to Gronwall lemma, we get 
\begin{multline*}
\left( \sup_{\xi \in \R^d} \frac{|\tilde D_t(\xi)|}{|\xi|^4} \right) \le
e^{(\lambda_K-K) \, t} \, \left( \sup_{\xi \in \R^d} \frac{|\tilde D_0(\xi)|}{|\xi|^4}
\right)
 \\ 
 + C_3 \, \left( \sum_{|i| \le 3} \left|M_i[f_0-g_0]\right| \right) \, 
\left( \frac{e^{-\lambda \, t}}{K-\lambda_K-\lambda} -  \frac{e^{(\lambda_K-K) \, t}}{K-\lambda_K-\lambda} \right).
\end{multline*}
Therefore, passing to the limit $K \to \infty$ and choosing (without any restriction) $\lambda \in (0,\bar\lambda)$, 
we obtain 
\[
\sup_{\xi \in \R^d} \frac{|\tilde D_t(\xi)|}{|\xi|^4}  \le
C \, e^{ - \lambda \, t} \, \left( \sup_{\xi \in \R^d} \frac{|\tilde D_0(\xi)|}{|\xi|^4}
+ \sum_{|i| \le 3} \left|M_i[f_0-g_0]\right| \right)
\]
or equivalently (and with the notations of  Example~\ref{expleFourierGal}), 
\[
|||Êd_t |||_4  \le C \, e^{ - \lambda \, t} \, |||Êd_0 |||_4 . 
\]

\smallskip\noindent
 {\sl Step 2. } From the preceding step and a trivial interpolation argument, we have 
\begin{multline*}
\left| f - g \right|_2 \le \left| f - g - \mathcal{M}[f-g] \right|_2 +
C \, \left( \sum_{|i| \le 3} \left|M_i[f-g]\right| \right)
\\
\le 
\left\| f - g - \mathcal{M}[f-g] \right\|_{M^1} ^{1/2}
\, \left| f - g - \mathcal{M}[f-g] \right|_4 ^{1/2} + 
C \, \left( \sum_{|i| \le 3} \left|M_i[f-g]\right| \right)
\\
\le C \, (1 +  M_4(f_0) + M_4(g_0)) \, e^{ - (\lambda/2) \, t}.
\end{multline*}
 We conclude by writing 
$$
\left| f - g \right|_2 \le  \left| f - g \right|_2^\eta \, \left| f - g \right|_2^{1-\eta},
$$
using Lemma~\ref{estim:C01TrueMax} for the first term and the previous decay estimates
for the second term.

\smallskip\noindent
 {\sl Step 3. } The same computations imply at least formally the same estimate on $h_t$ 
 as stated in Lemma~\ref{lem:a4maxwellfourier}. 
  Now, we proceed to the rigorous statement.  
 We consider a truncated (cut-off) model associated to $b_\eps= b \,
 {\bf 1}_{\cos \theta \le 1-\eps} \in L^1(S^{d-1})$ and an initial
 datum $h^\eps_0 \in L^1_2(\R^3)$ such that $\langle h^\eps_0, 1
 \rangle = 0$.  By standard argument there exists a unique solution
 $h^\eps_t \in C([0,\infty);L^1_2(\R^3))$ to the associated linearized
 Boltzmann equation arround $f_t$, and this one fulfills
 (\ref{ineq;tMMC10}).  Letting $\eps$ goes to $0$, we get that for
 some $h_t \in C([0,\infty);\SS'(\R^3)-w)$ there holds $h^\eps_t \wto
 h_t$ in $\SS'(\R^3)$ and $h_t$ is a weak solution to the linearized
 Boltzmann  equation (corresponding now to $b$) and
 which satisfies again (\ref{ineq;tMMC10}).  In other words, we
 have built a linear flow $\LL^{N \! L}_t [f_0] : \GG_1 \to \GG_1$ si
 that for any $h_0 \in \GG_1$ the distribution $\LL^{N \! L}_t [f_0]
 (h_0) := h_t$ is the solution to linearized Boltzmann equation.
\qed

\begin{lem}\label{lem:a4maxwellfourier2}
There exists $\lambda \in (0,\infty)$  and for any $\eta \in (0,1)$ there exists $C_\eta$ such that for any 
$f_0, g_0 \in P_{\GG_1,{\bf r}}$, ${\bf r}Ê\in {\bf R_{\GG_1}}$, we have
  $$
  \left| \omega_t \right|_{4} \le C \, e^{- (1-\eta) \, \lambda \,
    t} \, \big|g_0 - f_0\big|_2 ^{1+\eta} 
  $$
  with 
  \[
  \omega_t := g_t - f_t - h_t = S^{N \! L}_t(g_0) - S^{N \! L}_t(f_0)-
  \LL^{N \! L}_t[f_0] (g_0 - f_0)
  \]
  (as proved below, $\omega_t$ always has vanishing moments up to
  order $3$). 
\end{lem}

\noindent
{\it Proof of Lemma~\ref{lem:a4maxwellfourier2}.} 
The following arguments can be fully justified as in \cite{TV} by
truncating $b$ and passing to the limit.
Consider the error term
$$
\omega := g-f-h, \qquad \Omega := \hat \omega.
$$
Again we perform the angular cutoff of Lemma \ref{lem:contraction}
with cutoff parameter $K$, then the evolution equation (in the Fourier
side) satisfied by $\Omega$ is
\begin{equation}\label{eq:BIevold}
  \partial_t \Omega = \hat Q(\Omega,S) + \hat Q^+(H,D).
\end{equation}

Let us prove that 
\[
\forall \, |i|\le 3, \ \forall \, t \ge 0, \quad M_i[\omega_t] :=
\int_{\R^d} v^i \, d\omega_t(v) = 0.
\]
Consider the equation on $\omega_t$: 
\[
\partial_t \omega_t = Q(\omega,(f+g)) + Q(h,(f-g))
\]
and the fact that, for maxwell molecules, the $i$-th moment of
$Q(f_1,f_2)$ is a sum of terms given by product of moments of $f_1$
and $f_2$ whose orders sum to $|i|$. Hence using that for some
absolute coefficients $(a_{i,j})$ we have
$$
\int_{\mathbb{S}^{d-1}} b(\cos\theta) \left[(v^i)' + (v^i)'_* - (v^i) -
(v^i)_* \right] \, d\sigma = \sum_{j, \, |j| \le i} a_{i,j} \, \left(v^j\right)
\, \left(v^{i-j}\right)_*,
$$
we deduce
\[
\forall \, |i|\le 3, \quad 
\frac{d}{dt} M_i[\omega_t]  = \sum_{j \le i} 
  a_{i,j} \, M_j[\omega_t] \, M_{i-j}[f_t+g_t]
 + \sum_{j \le i} a_{i,j} \, M_j[h_t] \, M_{i-j}[f_t-g_t]
\]
\Black
and since
\[
\forall \, |i|\le 1, \quad M_i[h_t]  = M_i[f_t-g_t] = 0
\]
we deduce 
\[
\forall \, |i|\le 3, \quad 
\frac{d}{dt} M_i[\omega_t]  = \sum_{j \le i} 
  a_{i,j} \, M_j[\omega_t] \, M_{i-j}[f_t+g_t] 
\]
and from the initial data $\omega_0 = 0$ this concludes the proof.

We now consider the equation in Fourier form
\[
\partial_t \Omega = \hat Q(\Omega,S) + \hat Q^+(H,D)
\]
and we deduce in distributional sense
$$
\left( \frac{d}{dt} {|\Omega (\xi)| \over |\xi|^4} + K \,
  {|\Omega (\xi)| \over |\xi|^4}\right) \le \TT_1 + \TT_2,
$$
where
\begin{eqnarray*}
  \TT_1
  &:=& \sup_{\xi \in \R^3} \int_{\mathbb{S}^{d-1}} 
  {b\left(\sigma \cdot \hat\xi\right) \over |\xi|^4} \, 
  \left( \left| \frac{\Omega (\xi^+) \, S (\xi^-)}{2} \right|  
    + \left|\frac{\Omega (\xi^-) \, S
        (\xi^+)}{2} \right| \right) \, d\sigma
  \\
  &\le& \sup_{\xi \in \R^3} \int_{\mathbb{S}^{d-1}} 
  b\left(\sigma \cdot \hat\xi\right) \, 
  \left( {\left|\Omega (\xi^+)\right| \over \left|\xi^+\right|^4} \, 
    {\left|\xi^+ \right|^4 \over \left|\xi\right|^4}
    + {\left|\Omega (\xi^-)\right| \over \left|\xi^-\right|^4} \, 
    {\left|\xi^-\right|^4 \over \left|\xi\right|^4} \right) \, d\sigma
  \\
  &\le&  \left( \sup_{\xi \in \R^3} 
    {\left|\Omega (\xi)\right|^2 \over \left|\xi\right|^2} \right) \,  
  \left( \sup_{\xi \in \R^3} 
    \int_{\mathbb{S}^{d-1}} b\left(\sigma \cdot \hat\xi\right) \, 
    \left( \left|\hat\xi^+\right|^4 + \left|\hat\xi^-\right|^4 \right)
    \, d\sigma \right) \\
  &\le&  \lambda_K \, \left( \sup_{\xi \in \R^3} 
    {\left|\Omega (\xi)\right| \over \left|\xi\right|^4} \right),
\end{eqnarray*}
where $\lambda_K$ was defined in Lemma \ref{lem:contraction}, and
\begin{eqnarray*}
  \TT_2 &:=& {1\over 2} \,  \sup_{\xi \in \R^3} \int_{\mathbb{S}^{d-1}} 
  {b\left(\sigma \cdot \hat\xi\right) \over |\xi|^4} \, 
  \left| H (\xi^+) \, D (\xi^-) +  H (\xi^-) \, 
    D (\xi^+) \right|  \,  d\sigma
  \\ 
  &\le& {1 \over 2} \, \sup_{\xi \in \R^3} \int_{\mathbb{S}^{d-1}} 
  b\left(\sigma \cdot \hat\xi\right)
  \, \left(  {| H (\xi^+) | \over |\xi^+|^2} \, {| D (\xi^-) | \over |\xi^-|^2}
    \, {| \xi^- |^2 \over |\xi|^2}  + {| D (\xi^+) |^2 \over |\xi^+|^2} \, 
    {| H (\xi^-) |^2 \over |\xi^-|^2}
    \, {| \xi^- |^2 \over |\xi|^2}  \right) \, d\sigma
  \\
  &\le& \left| h_t \right|_2 \, \left| d_t \right|_2  \, 
  \int_{\mathbb{S}^{d-1}} b\left(\sigma \cdot
    \hat\xi_0\right) \, 
  \left(1 - \sigma \cdot \hat\xi_0\right) \, d\sigma \\
  &=& C \, e^{-(1-\eta) \, \lambda \, t} \, 
  \left| h_0 \right|_2 \, \left| d_0 \right|_2 ^\eta \le 
  C \, e^{-(1-\eta) \, \lambda \, t} \,  \left| d_0 \right|_2 ^{1+\eta}
\end{eqnarray*}
by using the estimates of Lemma~\ref{lem:contraction}.

Hence we obtain
$$
\left( \frac{d}{dt} {|\Omega (\xi)| \over |\xi|^4} + K \, {|\Omega
    (\xi)| \over |\xi|^4}\right) \le \lambda_K \, \left( \sup_{\xi \in \R^3}
  {\left|\Omega (\xi)\right| \over \left|\xi\right|^4} \right) + C \,
e^{-(1-\eta) \, \lambda \, t} \, \left| d_0 \right|_2 ^{1+\eta}.
$$
We then from the Gronwall inequality and relaxing the cutoff parameter
$K$ as in Lemma \ref{lem:a4maxwellfourier} (assuming
without restriction $(1-\eta) \lambda \le \bar \lambda$
\[
\left( \sup_{\xi \in \R^3}
  {\left|\Omega_t (\xi)\right| \over \left|\xi\right|^4} \right) \le C
\, e^{-(1-\eta) \, \lambda \, t} \, \left| g_0 - f_0 \right|_2 ^{1+\eta}.
\]
\qed 






\subsection{Proof of (A5) in Wasserstein distance} We know from \cite{T1} that 
$$
\sup_{t \ge 0}ÊW_2( S^{NL}_t f_0, S^{NL}_t g_0) \le W_2 (f_0,g_0). 
$$
As a consequence, from $W_1 \le W_2 \le C_k \, W_1^{\eta(k)}$ and $W_1 = [ \cdot ]^*_1$, we 
deduce that {\bf (A5)} holds with $\Theta(x) = x^{\eta(k)}$, $\FF_3 = Lip$ and $P_{\GG_3} = P_2$
endowed with $\|Ê\cdot \|_{\GG_3} = [Ê\cdot ]^*_1$. 




\subsection{Proof of (A5) in negative Sobolev norms} 

It is also possible to prove easily, in the cutoff case,
that the weak stability holds in Sobolev space on finite time:
\begin{lem}\label{lem:stab-max-sob}
   For any $T \ge 0$ and $f_t$, $g_t$ solutions of the Boltzmann
   equation with maxwellian kernel and initial data $f_0$ and $g_0$,
   that there exists $s \in (d/2,d/2+1)$ and a constant $C_k$ such that
   $$
   \sup_{[0,T]} \left\| f_t - g_t \right\|_{H^{-s}} \le C_{T,s}
   \, 
   \left\| f_0 - g_0 \right\|_{H^{-s}} .
   $$
 \end{lem}
 Remark that this proves {\bf (A5)} with
 $\GG_3=P(\R^d)$, 
 endowed with the $\dot{H}^{-s}$ norm.
This is useful in order to obtain the optimal rate $1/\sqrt{N}$ of the
law of large numbers. 

Let us only sketch the proof. We integrate \eqref{eq:BoltzMaxD}
against $D/(1+|\xi|^{2})^s$:
\[
\frac{d}{dt} \| D \|_{\dot{H}^{-k}} ^2 = { 1 \over 2 } \, \int_\xi
\int_{\mathbb{S}^{d-1}} b\left(\sigma \cdot \hat \xi\right) \, \frac{\left[ D^-
    S^+ D + D^+ S^- \, D - 2 \, 
    |D|^2\right]}{(1+|\xi|^{2})^s} \, d\sigma \, d\xi
\]
and we use Young's inequality together with the bounds
\[
\left\|S^+\right\|_\infty, \quad \left\| S^{-} \right\|_\infty 
\le \|f+g\|_{M^1}\le 2
\]
to conclude. 

\section{Hard spheres}
\label{sec:hardspheres}
\setcounter{equation}{0}
\setcounter{theo}{0}


\subsection{The model}
\label{sec:modelHS}
The limit equation was introduced in Subsection~\ref{sec:introEB} and
the stochastic model has discussed
Subsection~\ref{sec:modelEBbounded}. We consider here the case of the
Master equation \eqref{eq:BoltzBddKolmo}, \eqref{defBoltzBddGN} and
the limit nonlinear homogeneous Boltzmann equation \eqref{el},
\eqref{eq:collop}, \eqref{eq:rel:vit} with $\Gamma(z  ) = |z|$.

\subsection{Statement of the result}
\label{sec:statementHS}

We have the following theorem: 

\begin{theo}[Hard spheres homogeneous Boltzmann equation]\label{theo:HS}
  Let us consider $f_0 \in P(\R^d)$ with compact support. We denote by
  $\EE$ and $A$ two positive constants such that
  $$
  M_2(f_0) \le \EE \qquad\hbox{and}\qquad \hbox{supp} \, f_0 \subset
  \{Êv \in \R^d, \,\, |v| \le A \}.
  $$
  Let us also consider the hierarchy of $N$-particle distributions
  $f^N _t = S^N _t(f_0 ^{N})$ issued
  from the tensorial initial datum $f^N _0 = f_0 ^{\otimes N}$, $N \ge
  1$. Let us finally fix some $\delta \in (0,1)$.

  Then there are some constants $k_1 \ge 2$, $C >0$, only depending on
  $\delta$ and $\mathcal{E}$, some constant $C_\eta$ depending on
  $\eta$, and some constants $a >0$ depending on the collision kernel
  such that for any
  \[ \ell \in \N^*, \quad \varphi = \varphi_1 \otimes \varphi_2 \otimes \dots \otimes \,
  \varphi_\ell \in W^{1,\infty}(\R^d)^{\otimes \ell}, \] 
  we have
  \begin{multline} \label{eq:cvgHS} 
    \forall \, N \ge 2 \, \ell, \quad 
    \sup_{[0,T]}\left| \left \langle
        \left( S^N_t(f_0 ^N) - \left( S^{N \! L}_t(f_0)
          \right)^{\otimes N} \right), \varphi \right\rangle \right|
    \\ \le C \, \ell \, \|\varphi\|_{ W^{1,\infty}(\R^d)^{\otimes \ell} } \, \Bigg[ \frac{ \|f_0\|_{M^1_1}}{N} +
    \frac{ \|f_0\|_{M^1_{k_1}}}{N^{1-\delta}} + e^{a \, A} \,
    \eps_3(N) \Bigg]
  \end{multline}
  and the dominant error is the last term given by
  \[
  \eps_3(N) =
  e^{- \left[\log\left(\frac{C_\eta}{N^{\frac{1}{d}-\eta}}\right)_-\right]^{1/2}}
  \]
  for any $\eta >0$ small (the constant $C_\eta$ blows up when $\eta
  \to 0$).  As a consequence we deduce propagation of chaos with explicit
  rate. 
  
 \end{theo}

 In order to prove Theorem~\ref{theo:HS}, we shall prove each
 assumption {\bf (A1)-(A2)-(A3)-(A4)-(A5)} of
 Theorem~\ref{theo:abstract} step by step. Its application then
 exactly gives Theorem~\ref{theo:HS}.  We fix $\FF_1=\FF_2=C_b(\R^d)$
 and $\FF_3 = \mbox{Lip}(\R^d)$.


\subsection{Proof of (A1)} In the proof of
Lemma~\ref{lem:momentsN} 
we have already proved that
$$
\forall \, t \ge 0, \quad \hbox{supp} \, f^N_t \subset  \mathbb{E}_N :=
\left\{V \in \R^{Nd}; \,\, M^N_{2}(V) \le \EE\right\},
$$
which is precisely {\bf (A1)-(i)} with $m_2(v) := |v|^2$, as well as
that for any $k \ge 2$,
\[
\sup_{t \ge 0} \left \langle f^N _t, M_k ^N \right\rangle \le C^N_k 
\]
where $C^N_k$ depends on $k$, $\EE$, on the collision kernel and on
the initial value $\langle f^N _0, M_k ^N \rangle$ which is uniformly
bounded in $N$ (in function of $k$ and $A$ for instance). That is
precisely {\bf (A1)-(ii)} with $m_1(v) := |v|^{k_1}$ where $k_1 \ge 2$
will be chosen (large enough) in section~\ref{sec:HSA3bis}.  As for
\textbf{(A1)-(iii)}, we remark that for a given N-particle velocity $V
= (V_1, ..., V_N) \in \R^{dN}$, we have
 $$
 V \in \hbox{supp} \, f_0^{\otimes N} \ \Longleftrightarrow \ \forall
 \, i=1, \dots, N, \ V_i \in \hbox{supp} \, f_0 \ \Rightarrow \ \forall
 \, i=1,\dots, N, \ m_3(V_i) \le m_3(A),
 $$
with $m_3(v) := e^{a \, |v|}$.  We conclude that 
$$
\hbox{supp} \, f_0^{\otimes N} \subset \left\{ V \in \R^{dN}; \
  M^N_{m_3}(V) \le m_3(A) \right\},
$$
and {\bf (A1)-(iii)} holds. 

\subsection{Proof of (A2)} \label{sec:HSA2} 
We define 
$$
P_{\GG_1} := \left\{ f \in P_{k_1}(\R^d); \ M_2(f) \le \EE \right\}
$$
that we endow with the total variation norm $\| \cdot \|_{\GG_1} := \|
\cdot \|_{M^1}$.
 
\smallskip\noindent $\bullet$ The proof of assertion (i), that is for
any $a \ge a_{k_1} > 0$ and for any $t > 0$, the application $f_0
\mapsto S^{N\!L}_t f_0$ maps $\BB P_{\GG_1,a} := \{ f \in
P_{\GG_1}, \,\, M_{k_1} (f) \le a \}$ continuously into itself, is
postponed to section~\ref{sec:HSA4} where it is proved in
\eqref{estim:dt} a H\"older continuity of the flow.
 
\smallskip\noindent
$\bullet$ For any $f,g \in P_{\GG_1}$ we have 
\begin{eqnarray} \nonumber \| Q(g,g) - Q(f,f) \|_{M^1} &=&\| \tilde
  Q(g-f,g+f) \|_{M^1} \\ \nonumber &\le& 2 \int_{\R^d} \int_{\R^d}
  \int_{\mathbb{S}^{d-1}} b(\theta) \, |v-v_*| \, |f - g | \,
  |f_*+g_*| \, d\sigma \, dv_* \, dv \\ \label{eq:HSA2ii} &\le& 2 \,
  (1+\EE) \, \| b \|_{L^1} \, \| (f - g) \, \langle v \rangle \|_{M^1}.
\end{eqnarray}
We deduce that 
$$
\| Q(g,g) - Q(f,f) \|_{M^1} \le  2 \, (1+\EE)^{3/2} \, \| b \|_{L^1} \,
\| f - g \|_{M^1}^{1/2}
$$
which yields $Q \in C^{0,1/2}(P_{\GG_1}; \GG_1)$. 

\smallskip\noindent
 $\bullet$ Finally, we have for any $f_0 \in P_{\GG_1}$
\begin{multline*}
\|f_t  - f_0 \|_{M^1_1} = \left\| \int_0^t Q(f_s,f_s) \, ds \right\|_{M^1_1} 
\le\int_0^t  \left\| Q(f_s,f_s) \right\|_{M^1_1}  \, ds  \le 3 \, t \, \| b \|_{L^1} \, (1+\EE)^2,
\end{multline*}
and then, using \eqref{eq:HSA2ii}, 
\begin{multline*}
\big\|f_t  - f_0 - t \, Q(f_0,f_0) \big\|_{M^1} 
= \left\| \int_0^t ( Q(f_s,f_s) - Q(f_0,f_0) ) \, ds \right\|_{M^1} \\
\le \int_0^t \big\| Q(f_s,f_s) - Q(f_0,f_0) \big\|_{M^1} \, ds \\
\le \int_0^t  2 \, (1+\EE) \, \| b \|_{L^1} \,  \| f_s - f_0 \|_{M^1_1} \, ds
  \le 3 \, t^2 \, \| b \|_{L^1}^2 \, (1+\EE)^3, 
\end{multline*}
which precisely means that $t \mapsto S^{N \! L}_t (f_0)$ is
$C^{1,1}$ at $t=0^+$ in the space $P_{\GG_1}$. A similar argument
would yield the $C^{1,1}([0,\tau],P_{\GG_1})$ regularity on a small
time interval $[0,\tau]$.

\subsection{Proof of (A3)} \label{sec:HSA3} For any $k \in \N^*$, any
energy bound $\EE$ (which has been fixed once for all in the statement
of the main result), any energy $r_{d+1} \in [0,\EE]$ and any mean
velocity $ (r_1, ..., r_d) \in B_{\R^d}(0,r_{d+1})$ we define
$$
P_{\GG_1,{\bf r}} := \left\{Êf \in P_{k}(\R^d); \,\, \langle f, v_j \rangle = r_j, \, j = 1, ..., d, 
\,\, M_2(f) = r_{d+1} \right\}, \quad {\bf r}:= (r_1, ..., r_{d+1}). 
$$
We also denote by ${\bf R}_\EE$ the set of all admissible vectors
${\bf r}Ê\in \R^{d+1}$ constructed as above.  We claim that there
exists $C = C_{k_1,\EE}$ such that for any $\eta \in (0,1)$ and any
function
$$
\Phi \in \bigcap_{{\bf r } \in {\bf R}_\EE}
C_\Lambda^{1,\eta}(P_{\GG_1,{\bf r}};\R),
$$
with $\Lambda(f) := M_{k_1}(f)$, we have
\begin{equation}\label{eq:HSA3}
\forall \, V \in \mathbb{E}_N, \quad
\left| G^N (\Phi \circ \mu^N_V) - (G^\infty \Phi) (\mu^N_V) \right| 
\le C \, 
\left( \sup_{ {\bf r} \in {\bf R}_\EE} [ \Phi
  ]_{C^{1,\eta}_\Lambda(P_{\GG_1,{\bf r}})} \right) \, {M^N_{k_1}(V)\over N^\eta},
\end{equation}
which is precisely assumption {\bf (A3)} with $\eps(N) = C/
N^{\eta}$. 

Consider $V \in \mathbb{E}_N$ and define 
$$
{\bf r}_V := \left(\langle \mu^N_V,z_1 \rangle, \dots, \langle
  \mu^N_V,z_d \rangle, \langle \mu^N_V, |z|^2 \rangle\right) \in {\bf
  R}_\EE.
$$
Then for a given $\Phi \in C_\Lambda^{1,\eta}(P_{\GG_1,{\bf r}_V};\R)$,
we set $\phi := D\Phi[ \mu_V ^N]$ and we compute:
\begin{eqnarray*}
&&  G^N(\Phi \circ  \mu ^N _V) 
=  \frac{1}{2N} \, \sum_{i,j=1} ^N |V_i - V_j| \, 
  \int_{\mathbb{S}^{d-1}} \left[ \Phi\left(\mu^N_{V^*_{ij}}\right) 
    - \Phi\left(\mu^N _V\right) \right] \, b(\theta_{ij}) \, d\sigma 
\\
&&= \frac{1}{2N} \, \sum_{i,j=1} ^N |V_i - V_j| \,
  \int_{\mathbb{S}^{d-1}} \left\langle \mu^N _{V^*_{ij}}
    - \mu^N _V , \phi \right\rangle \, b(\theta_{ij}) \, d\sigma 
\\
&&+  \frac{[\Phi]_{ C_\Lambda^{1,\eta}(P_{\GG_1,{\bf r}_V}) }}{2N} \, 
  \sum_{i,j=1} ^N |V_i - V_j|  \int_{\mathbb{S}^{d-1}} 
  \!\!\! \ \max\left\{ M_{k_1}\left(\mu^N _{V^*_{ij}}\right) ;\, 
    M_{k_1}\left( \mu^N _{V}\right) \! \right\} \, 
  \OO \! \left( \left\| \mu^N _{V^*_{ij}} 
      - \mu^N _V \right\|_{M^1}^{1+\eta} \right)  \, d\sigma \\
&& =:  I_1(V) + I _2(V).
\end{eqnarray*}
For the first term $I_1(V)$, thanks to Lemma~\ref{lem:H0}, we have 
\begin{multline*}
  I_1(V) = \frac{1}{2N^2} \, \sum_{i,j=1} ^N |V_i - V_j| \,
  \int_{\mathbb{S}^{d-1}} b(\theta_{ij}) \, \left[ \phi(V_i ^*) +
    \phi(V_j^*) - \phi(V_i) - \phi(V_j) \right] \, d\sigma \\
  = \frac{1}{2N^2} \, \int_v \int_w |v - w| \,
  \int_{\mathbb{S}^{d-1}} b(\theta) \, \left[ \phi(v^*) +
    \phi(w^*) - \phi(v) - \phi(w) \right] \, \mu^N _V(dv) \, 
  \mu^N_V (dw) \, d\sigma \\
  = \left\langle Q(\mu^N _V, \mu^N _V), \phi \right\rangle = 
  \left(G^\infty\Phi\right)(\mu^N _V).
\end{multline*}
For the second term $I_2(V)$, using that 
\begin{eqnarray*}
  M_{k_1}\left(\mu^N_{V^*_{ij}}\right) 
  &:=& M^N_{k_1}\left(V^*_{ij}\right) = {1 \over N} \, 
  \left( \left(\sum_{n \not= i,j} |V_n|^{k_1} \right) + |V^*_{i}|^{k_1} +
    |V^*_{j}|^{k_1} \right) 
  \\
  &\le& {1 \over N} \, 
  \left( \left( \sum_{n \not= i,j} |V_n|^{k_1} \right) 
    + 2 \, \left(|V_i|^2 + |V_j|^2\right)^{k_1/2} \right) 
  \\
  &\le& \frac{2^{1+k_1/2}}{N} \, 
  \left( \sum_{n} |V_n|^{k_1} \right) = C \, M^N_{k_1}(V),
\end{eqnarray*}
we deduce 
\begin{eqnarray*}
  |I_2(V)|  
  &\le& C \, {  M^N _{k_1}(V)} \, 
  [\Phi]_{C^{1,\eta}_\Lambda(P_{\GG_1,{\bf r}_V})}  \, 
  \left( \frac{1}{2N} \, \sum_{i,j=1} ^N |V_i - V_j| \,
  \left( \frac4N \right)^{1+\eta} \right) \\
  & \le & C \, {  M^N _{k_1}(V)} \, 
 [\Phi]_{C^{1,\eta}_\Lambda(P_{\GG_1,{\bf r}_V})}  \,  \left( {1 \over N^\eta} \, 
 \frac{1}{N^2} \, \sum_{i,j=1} ^N \left[ \langle V_i \rangle + 
    \langle V_j \rangle \right] \right) \\
  & \le & C \, {  M^N _{k_1}(V)} \, 
  [\Phi]_{C^{1,\eta}_\Lambda(P_{\GG_1,{\bf r}_V})}  \,  {1 \over N^\eta} \, (1+\EE).
 \end{eqnarray*}

We conclude that \eqref{eq:HSA3} holds by gathering these two estimate.

\subsection{Proof of (A4) on a finite time interval $[0,T]$} 
\label{sec:HSA4}

\begin{lem}\label{lem:expansionHS} 
  For any given energy $\EE > 0$ and any $\delta > 0$ there exists
  some constants $k_1 \ge 2$ (depending on $\EE$ and $\delta$) and $C$
  (depending on $\EE$), such that for any $f_0, g_0 \in
  P_{\GG_1}(\R^d)$, we have
  \begin{eqnarray} \label{estim:dt} && \left\| S^{N \!  L}_t(g_0) -
      S^{N \! L}_t(f_0) \right\|_{M^1} \le e^{C \, (1+t)} \, \sqrt{
      \max \left\{ M_{k_1}(f_0),M_{k_1}(g_0) \right\} } \, \left\| f_0
      - g_0 \right\|_{M^1}^{1-\delta}, \\ \label{estim:ht} && \left\|
      \mathcal{L}S^{N \! L}_t[f_0]\left(f_0-g_0\right) \right\|_{M^1}
    \le e^{C \, (1+t)} \, \sqrt{ M_{k_1}(f_0) } \, \left\| f_0 - g_0
    \right\|_{M^1}^{1-\delta}, \\ \label{estim:dt-ht} && \left\| S^{N
        \!  L}_t(g_0) - S^{N \! L}_t(f_0) - \mathcal{L}S^{N \!
        L}_t[f_0](g_0-f_0) \right\|_{M^1} \\ \nonumber 
    &&\qquad\qquad\qquad\qquad\qquad\quad
    \le e^{C \, (1+t)} \, \sqrt{ \max \left\{
        M_{k_1}(f_0),M_{k_1}(g_0) \right\} } \, \left\| f_0 - g_0
    \right\|_{M^1}^{2-\delta}.
  \end{eqnarray}
\end{lem}

\noindent
{\bf Proof of Lemma~\ref{lem:expansionHS}.} We proceed in several
steps. Let us define 
$$
\forall \, f \in M^1(\R^d), \quad \| f \|_{M^1_k} := \int_{\R^d} \langle v \rangle^k \, |f|(dv), 
\quad
\| f \|_{M^1_{k,\ell}} := \int_{\R^d} \langle v \rangle^k \, (1 + \log \langle v \rangle)^\ell \, |f|(dv).
$$

\smallskip\noindent{\sl Step 1. The strategy.}
Let us define 
\[
\left\{
\begin{array}{l}
\partial_t f_t = Q(f_t,f_t), \qquad f_{|t=0} = f_0 \vspace{0.3cm} \\
\partial_t g_t = Q(g_t,g_t), \qquad g_{|t=0} = g_0 \vspace{0.3cm} \\
\partial_t h_t = Q(f_t,h_t) + Q(h_t,f_t), \qquad h_{|t=0} = h_0 := g_0 - f_0.
\end{array}
\right.
\]
Existence and uniqueness for $f_t$, $g_t$ and $h_t$ is a consequence
of the following important stability argument that we use several
times. This estimate is due to DiBlasio~\cite{DiB74} in a $L^1$
framework, and it has been recently extended to a measure framework in
\cite[Lemma 3.2]{EM}. Let us recall the argument for $h$. We first
write
\begin{eqnarray}\label{eq:htM12}
{d \over dt} \int \langle v \rangle^2 \, |h_t|(dv) 
\!\!&\le& \!\!
\int\!\! \! \int\!\! \! \int |h_t|(dv) \, f_{t} (dv_*)\, |u| \,
b(\theta) \, \Big[ \langle v' \rangle^2 \!+\! \langle v'_* \rangle^2\!
-\! \langle v \rangle^2 \!-\!  \langle v_* \rangle^2 \Big] \, d\sigma 
\\ \nonumber
&&  + 2 \int\!\!\! \int\!\! \! \int  |h_t|(dv) \, f_{t} (dv_*)\,  |u| \, b(\theta) \, \langle v_* \rangle^2  \, d\sigma \, dv \, dv_*
\end{eqnarray}
(this formal computation can be justified by a regularization
proceedure, we refer to \cite{EM} for instance). Because the first
term vanishes, we conclude with
\begin{equation}\label{ineq:htHS}
{d \over dt} \| h_t \|_{M^1_2} \le C \, \| f \|_{M^1_3} \, \| h_t \|_{M^1_2}.
\end{equation}
When $\|Êf_t \|_{M^1_3} \in L^1(0,T)$, we may integrate this
differential inequality and we deduce that $h$ is unique. 

More precisely, we have established
\begin{equation}\label{eq:unifh}
\sup_{[0,T]} \| h_t \|_{M^1_2} \le  \left\| g_0 - f_0 \right\|_{M^1 _2} \, 
\exp \left( C \, \int_0 ^T \left\| f_s \right\|_{M^1 _3} \, ds \right), 
\end{equation}
and similar arguments imply
\begin{equation}\label{eq:unifd}
\sup_{[0,T]} \| f_t - g_t \|_{M^1_2} \le  \left\| g_0 - f_0 \right\|_{M^1 _2} \, 
\exp \left( C \, \int_0 ^T \left\| f_s +g_s \right\|_{M^1 _3} \, ds \right).
\end{equation}
It is worth mentioning that one cannot prove $\| f_t \|_{M^1_3} \in
L^1(0,T)$ under the sole assumption $\| f_0 \|_{M^1_{2}} < \infty$
because since it would contradict the non-uniqueness result
of~\cite{LuW02}. However, thanks to Povzner's inequality, one may show
that $\| f_t \|_{M^1_3} \in L^1(0,T)$ whenever $\| f_0 \|_{M^1
  _{2,1}}$ is finite, with the definition
\[
\| f_0 \|_{M^1_{k,\ell}} := \int \langle v \rangle^k \, \log \left( \langle v
  \rangle \right)^\ell \, df_0(v) < + \infty
\]
(see \eqref{fM13L1t} below or \cite{MW99,Lu99}), which will be the key
step for establishing \eqref{estim:dt} and \eqref{estim:ht}.

\smallskip Now, our goal is to estimate (in terms of $\| g_0 - f_0
\|_{M^1}$) the $M^1$ norm of 
\[ 
\zeta_t := f_t -g_t - h_t.
\]
The measure $\zeta_t$ satisfies the following evolution equation:
\[
\partial_t \zeta_t = Q(f_t,f_t) - Q(g_t,g_t) - 
Q(h_t,f_t) - Q(f_t,h_t), \quad \zeta_0 = 0.
\]
We can rewrite this equation as
\[
\partial_t \zeta_t = Q(\zeta_t, f_t+g_t) + Q(h_t,f_t-g_t).
\]
The same arguments as in \eqref{eq:htM12}-\eqref{ineq:htHS} yield the
following differential inequality
\[ 
\frac{d}{dt} \left\|\zeta_t\right\|_{M^1 _2} \le 
C \, \left\|\zeta_t\right\|_{M^1 _2} \, 
 \left\|f_t+g_t\right\|_{M^1_3}  + \left\| \tilde Q(h_t,f_t-g_t)
\right\|_{M^1_2}, \quad  \left\|\zeta_0\right\|_{M^1 _2} = 0.
\]
We deduce 
\begin{equation}\label{estim:uwv}
 \sup_{t \in [0,T]} \left\|\zeta_t\right\|_{M^1 _2} \le  
  \left( \int_0 ^T \left\| \tilde Q(h_s,f_s-g_s)
    \right\|_{M^1_2} \, ds \right) \, 
  \exp\left(C \, \int_0 ^T \left\|f_s+g_s\right\|_{M^1_3} \, ds\right).
\end{equation}
Since
\begin{eqnarray}\label{eq:w} 
  \int_0^T  \left\| \tilde Q(h_s,f_s-g_s)
  \right\|_{M^1_2} \, ds &\le& 
  C \, \left( \sup_{[0,T]} \left\| h_t \right\|_{M^1_2} \right) \,
  \left( \int_0^T \| f_s -g_s \|_{M^1_3} \, ds \right) \\ \nonumber 
  && + \, C \, \left( \sup_{[0,T]} \left\| g_t - f_t
  \right\|_{M^1_2} \right) \, \left( \int_0^T \| h_s \|_{M^1_3} \, ds \right),
\end{eqnarray}
we deduce from~\eqref{eq:unifh} and \eqref{eq:unifd}
\begin{multline*}
 \sup_{[0,T]} \left\| \zeta_t \right\|_{M^1_2} \le C
  \, \left\| g_0 - f_0 \right\|_{M^1 _2} \, \exp \left( \int_0 ^T
    \left\| f_s \right\|_{M^1 _3} + \left\| g_s \right\|_{M^1 _3}
  \right) 
  \\
   \left\{ \left( \int_0 ^T \left\| g_s - f_s \right\|_{M^1
        _3} \right) \, \exp \left( \int_0 ^T \left\| f_s \right\|_{M^1
        _3} \right) + \left( \int_0 ^T \left\| h_s \right\|_{M^1 _3}
    \right) \, \exp \left( \int_0 ^T \left\| f_s \right\|_{M^1 _3} +
      \left\| g_s \right\|_{M^1 _3} \right) \right\}.
\end{multline*}
Hence the problem is reduced to time integral controls over
$\|f_t\|_{M^1_3}$, $\|g_t\|_{M^1_3}$, $\|f_t-g_t\|_{M^1_3}$ and $\|
h_t\|_{M^1_3}$.

\smallskip\noindent{\sl Step 2. Time integral control of $f$ and $g$
  in $M^1 _3$.}  In this step we prove
\begin{equation}\label{fM13L1t} 
  \int_0^T \| f_t
  \|_{M^1_{3,\ell-1}}  \, dt \le C_\EE \, T + C' \, \| f_0
  \|_{M^1_{2,\ell}} \quad \ell = 1, 2, 
\end{equation} 
for $f$, where $C_\EE$ is some energy dependent constant, and $C'$ is
a numerical constant. The same estimate obviously holds for $g$.
These estimates are a consequence of the accurate version of the
Povzner inequality as one can find in \cite{MW99,Lu99}. Indeed it has
been proved in \cite[Lemma 2.2]{MW99} that for any convex function
$\Psi : \R^d \to \R$, $\Psi (v) = \psi (|v|^2)$, the solution $f_t$ to
the hard spheres Boltzmann equation satisfies
$$
{d \over dt} \int_{\R^d} \Psi(v) \, f_t(dv) = \int_{\R^d}\!
\int_{\R^d} f_t(dv) \, f_t (dv_*) \, |v-v_*| \, K_\Psi(v,v_*)
$$
with $K_\Psi = G_\Psi - H_\Psi$, where the term $G_\Psi$ ``behaves
mildly'' and the term $H_\Psi$ is given by (see~\cite[formula
(2.7)]{MW99})
$$
H_\Psi (v,v_*) = 2\pi \, 
\int_0^{\pi/2} \Big[ \psi (|v|^2 \cos^2 \theta + |v_*|^2\sin^2 \theta
  ) - \cos^2 \theta \, \psi (|v|^2) - \sin^2 \theta  \, \psi (|v_*|^2
  ) \Big] \, d\theta,
$$
(note that $H_\Psi \ge 0$ since its integrand is nonnegative from the
convexity of $\psi$).  More precisely, in the cases that we are
interested with, namely $\Psi(v) = \psi_{2,\ell}(|v|^2)$ with
$\psi_{k,\ell}(r) = r^{k/2} \, (\log \, r)^\ell$ and $\ell = 1, 2$, it has
been established that (with obvious notations) 
\begin{equation}\label{eq:BddG}
  \forall \, v,v_* \in \R^d, \quad 
  \left|G_{\psi_{2,\ell}}(v,v_*)\right| \le 
  C_\ell \, \langle v \rangle \, (\log (1+
  \langle v \rangle^2))^\ell \, \langle v_* \rangle \, (\log (1 +
  \langle v_* \rangle^2))^\ell.
\end{equation}
On the other hand, in the case $\ell = 1$ we easily compute (with the
notation $x := \cos^2 \theta$ and $u = |v_*|/|v|$)
\begin{multline*} 
\forall \, x,\, u \in \R,
\ x \in [1/4,3/4], \, u \in [0,1/2], \\
\psi_{2,1}
\Big(|v|^2 \cos^2 \theta + |v_*|^2\sin^2 \theta \Big) - \cos^2 \theta \,
\psi_{2,1} \left(|v|^2\right) - \sin^2 \theta \, \psi_{2,1} 
\left(|v_*|^2 \right)=
\\
 = |v|^2 \, \Big[ (1-x) \, \psi_{2,1} \left(u^2\right) 
   + x \, \psi_{2,1} (1) - \psi_{2,1} \left((1-x) \, u^2 
     + x\right) \Big] \ge \kappa_0 \, |v|^2,
\end{multline*} 
for some numerical constant $\kappa_0 > 0$ (note that this lower bound
only depends on the strict convexity of the real function
$\psi_{2,1}$). We straightforwardly deduce that there exists a
numerical constant $\kappa_1 > 0$ such that
$$
H_{2,1} (v,v_*) \ge \kappa_1 \, |v|^2 \, {\bf 1}_{|v| \ge 2 \, |v_*|}.
$$

Similarly, in the case $\ell = 2$, we have
\begin{multline*}
\forall \, x,\, u \in \R, \ x \in [1/4,3/4], \, u \in [0,1/2], \\
\psi_{2,2} \Big(|v|^2 \cos^2 \theta + |v_*|^2\sin^2 \theta \Big) - \cos^2
\theta \, \psi_{2,2} \left(|v|^2\right) - \sin^2 \theta  \, \psi_{2,2} 
\left(|v_*|^2\right)= \\
= 2 \, |v|^2  \, \log |v|^2 \,  \left\{ (1-x) \, \psi_{2,1} \left(u^2\right) 
+ x \, \psi_{2,1} (1) - \psi_{2,1} \left((1-x) \, u^2 + x\right)
\right\} \\
+  |v|^2 \,  \Big[ (1-x) \, \psi_{2,2} \left(u^2\right) 
+ x \, \psi_{2,2} (1) - \psi_{2,2} \left((1-x) \, u^2 + x\right) \Big] 
\ge 2\, \kappa_0 \,  |v|^2 \, \log |v|^2,
\end{multline*}
and then 
$$
H_{2,2} (v,v_*) \ge 2 \, \kappa_1 \, |v|^2 \, \log |v|^2 \, {\bf
  1}_{|v| \ge 2 \, |v_*|}.
$$

Putting together the estimates obtained on $G_{2,\ell}$ and
$H_{2,\ell}$ we deduce the Povzner inequality
\begin{equation}\label{estim:Povnzer} 
|v-v_*| \, K_{2,\ell} \le C \, (\langle v
\rangle^2 + \langle v_* \rangle^2 ) - \kappa \, |v|^3 \, (\log \langle
v \rangle)^{\ell-1}, 
\end{equation} 
and we finally obtain the differential inequality
\[
\frac{d}{dt} \|f_t\|_{M^1_{2,\ell}} \le 2\, C \, (1 + \EE) - \kappa \, M_{3,\ell-1},
\]
from which \eqref{fM13L1t} follows. 


\smallskip\noindent{\sl Step 3. Exponential time integral control of
  $f$ and $g$ in $M^1 _3$ (proof of \eqref{estim:dt} and
  \eqref{estim:ht}).}
Let us first prove that 
\begin{equation}\label{estim:expfg} 
e^{C \, \int_0 ^T \left\| f_s \right\|_{M^1 _3} \, ds} \le
\sqrt{M_k(f_0)}, \qquad 
e^{C \, \int_0 ^T \left\| g_s \right\|_{M^1 _3} \, ds} \le
\sqrt{M_k(g_0)}
\end{equation}
for any $k \ge k_\EE$, with $k_\EE$ big enough. 

This is a straightforward consequence of the previous step and the
following interpolation argument. For any given probability measure $f
\in P_k(\R^d)$ with $M_2(f) \le \EE$, we have for any $a > 2$
\begin{eqnarray*}
\| f \|_{M^1_{2,1}} &=& \int_{\R^d} \langle v \rangle^2 \, \left(1+
\log(\langle v \rangle^2)\right) \, \left({\bf 1}_{\langle v \rangle^2 \le a} +
{\bf 1}_{\langle v \rangle^2 \ge a} \right) \, f(dv)
\\
&\le& (1+\EE) \, (1+ \log a) + \frac{1}{a} \, \int_{\R^d} \langle v
\rangle^4 \, \left(1+ \log(\langle v \rangle^2)\right) \, f(dv)
\\
&\le& (1+\EE) \, (1+ \log a) + \frac{1}{a} \, \| f \|_{M^1_5}.
\end{eqnarray*}
By choosing $a:= \| p \|_{M^1_5}$, we get
\begin{equation}\label{estim:M121M15} 
\|p\|_{M^1_{2,1}} \le 2 \, (1+\EE) \,
\left(1+ \log \| p \|_{M^1_5}\right).
\end{equation} 
On the other hand, the following elementary H\"older inequality holds
\begin{equation}\label{eq:Holderkell} 
\forall \, k,k' \in \N, \ k' \le k, \ \forall \, f \in M^1_k, \quad 
\| f \|_{M^1_{k'}} \le
\| f \|_{M^1} ^{1-k'/k} \, \| p\|_{M^1_k} ^{j/k}.
\end{equation} 
Then estimate~\eqref{estim:expfg} follows from \eqref{fM13L1t},
\eqref{estim:M121M15} and \eqref{eq:Holderkell} with $k' = 2$ and $k =
k_\EE \ge 5$ large enough in such a way that
\[
C' \, 2 \, (1+\EE) \, \frac{5}{k} \le \frac12,
\] 
where $C'$ is the constant which appears in~\eqref{fM13L1t}.  We then
deduce \eqref{estim:dt} from \eqref{eq:unifd}, and (similarly)
\eqref{estim:ht} from \eqref{eq:unifh}.


\smallskip\noindent{\sl Step 4. Time integral control on $d$ and $h$.} 
Let us prove 
\begin{multline}\label{estim:h&dM131L1T}
\int_0^T \left(
\left\| d_t \right\|_{M^1_{3}} + \left\| h_t \right\|_{M^1_{3}} \right) \, dt
\\ \le 
e^{C_\EE \, (1+T)} \, e^{C' \, \left( \| f_0 \|_{M^1 _{2,1}} + \| g_0
  \|_{M^1_{2,1}} \right)} \, 
\left( \left\| f_0 \right\|_{M^1 _{2,2}} + \left\| g_0
  \right\|_{M^1_{2,2}} \right) \, \left\| g_0 - f_0 \right\|_{M^1_2} 
+ \left\| g_0 - f_0 \right\|_{M^1 _{2,1}},
\end{multline}
for some energy dependent constant $C_\EE$ and some numerical constant
$C'$. Performing similar computations to those leading to
\eqref{eq:htM12}, we obtain
\begin{eqnarray*} {d \over dt} \| h_t \|_{M^1_{2,1}} &\le& \int\!\!
  \int\!\! |h_t|(dv) \, f_{t} (dv_*)\, |u| \, K_{2,1}(v,v_*) \\
  \nonumber && + 2 \, C \, \int\!\! \int\!\! \int |h_t|(dv) \, f_{t}
  (dv_*)\, |u| \, \langle v_* \rangle^2 \, \left(1 + \log \langle v_*
  \rangle^2 \right).
\end{eqnarray*}
Thanks to the Povzner inequality \eqref{estim:Povnzer}, we deduce for some numerical constants $C,\kappa > 0$
\[ 
{d \over dt} \left\| h_t \right\|_{M^1_{2,1}} \le 
C \, \left\| h_t \right\|_{M^1_{2}} \,
\left\| f_t \right\|_{M^1_{3,1}} - \kappa \, \left\| h_t \right\|_{M^1_{3}} .
\]
Integrating that differential inequality yields 
\[
\left\| h_T \right\|_{M^1_{2,1}} + \kappa \, 
\int_0^T \left\| h_t \right\|_{M^1_{3}} \, dt \le 
C \, \left( \sup_{[0,T]} \left\| h_t \right\|_{M^1_{2}} \right) \, 
\left( \int_0^T \left\| f_t \right\|_{M^1_{3,1}} \, dt \right)
+ \left\| h_0 \right\|_{M^1_{2,1}}.
\]
Using the previous pointwise control on $\| h_t \|_{M^1_{2}}$
and~\eqref{fM13L1t} (with $\ell = 2$) we deduce
\begin{equation}\label{eq:hM131L1t}
\int_0^T \left\| h_t \right\|_{M^1_{3}} \, dt \le 
e^{C_\EE \, (1+T)} \, e^{C_1 \, \| f_0 \|_{M^1_{2,1}}} \, 
\left\| f_0 \right\|_{M^1_{2,2}} \, \left\| g_0 - f_0 \right\|_{M^1_2}
+ \left\| g_0 - f_0 \right\|_{M^1 _{2,1}}.
\end{equation}
Arguing similarly for $d_t$, we deduce~\eqref{estim:h&dM131L1T}.


\smallskip\noindent{\sl Step 5. Conclusion. }  By gathering the 
estimates
\eqref{estim:uwv}-\eqref{eq:w}-\eqref{estim:expfg}-\eqref{estim:h&dM131L1T},
we obtain
\begin{multline*}
\sup_{[0,T]} \| \zeta_t \|_{M^1_2} \le e^{C_\EE \, (1+T)} \, e^{C' \,
  \left( \left\| f_0 \right\|_{M^1 _{2,1}} + \left\| g_0 \right\|_{M^1
    _{2,1}} \right) } \\ \left( \left\| f_0 \right\|_{M^1 _{2,2}} +
\left\| g_0 \right\|_{M^1 _{2,2}} \right) \, \left\| g_0 - f_0
\right\|_{M^1_2} \, \left\| g_0 - f_0 \right\|_{M^1 _{2,1}},
\end{multline*}
for some energy dependent constant $C_\EE$ and some numerical constant
$C'$. Then arguing as in the end of step 3, using \eqref{estim:M121M15} and
\eqref{eq:Holderkell} with $k$ large enough, we get
$$
\sup_{[0,T]} \| \zeta_t \|_{M^1_2} \le  e^{C \, (1+t)} \, \sqrt{ \max \left\{
        M_{k}(f_0),M_{k}(g_0) \right\} } \, \left\| f_0 - g_0 \right\|_{M^1}^{2-5/k},
$$
from which estimate \eqref{estim:dt-ht} follows.
\qed 

\subsection{\bf  Proof of (A4) uniformly in time.} 
\label{sec:HSA3bis}

Let us start from an auxiliary result. It was proved in
\cite{MouhotCMP} that the nonlinear and linearized Boltzmann
semigroups for hard spheres satisfy
\begin{equation}
  \label{eq:HSdecay}
  \left\| S^{N \! L} _t \right\|_{L^1(m_z ^{-1})} \le C_z \,
  e^{-\lambda \, t}, \quad 
   \left\| e^{\mathcal L \, t} \right\|_{L^1(m_z ^{-1})} \le C_z \,
  e^{-\lambda \, t}
\end{equation}
where $m_z(v) := e^{z \, |v|}$, $z > 0$, $\lambda=\lambda(\EE)$ is the
optimal rate, given by the first non-zero eigenvalue of the linearized
operator $\mathcal L$ in the smaller space $L^2(M^{-1})$ where $M$ is
the maxwellian equilibrium (see \cite[Theorem~1.2]{MouhotCMP}), and
$C_z$ is some constant depending on $z$ and the energy $\EE$.

\begin{lem}
  For any given energy $\EE > 0$, there exists some constants $k_1 \ge
  2$ (depending on $\EE$ and $\delta$) and $C$ (depending on $\EE$)
  and $\eta \in (0,1)$, such that for any $f_0, g_0 \in
  P_{\GG_1}(\R^d)$ satisfying
\[
\forall \, i=1, \dots, d, \quad 
\left\langle f_0, v_i \right\rangle = \left\langle g_0, v_i \right\rangle 
\quad \mbox{ and } \quad \left\langle
f_0, |v|^2 \right\rangle = \left\langle g_0, |v|^2 \right\rangle \le \EE,
\]
we have
\begin{eqnarray}\label{estim:dt-infty}
&&  \left\| S^{N\! L}_t(g_0) - S^{N\! L}_t(f_0) \right\|_{M^1_2} \le 
  e^{C \, - {\lambda \over 2} \, t }  \,  
     \sqrt{ \max \left\{ M_{k_1}(f_0), \, M_{k_1}(g_0) \right\}  }\, 
  \left\| g_0 - f_0 \right\|_{M^1}^{{1+\eta \over 2}},
\\ \label{estim:ht-infty}
&& \left\| \mathcal{L}S^{N \! L}_t[f_0](g_0-f_0) \right\|_{M^1_2} \le 
  e^{C \, - {\lambda \over 2} \, t } \,    \sqrt{  M_{k_1}(f_0) } \, 
  \left\| g_0 - f_0 \right\|_{M^1}^{{1+\eta \over 2}}, 
\\ \label{estim:dt-ht-infty}
&&  \left\| S^{N\! L}_t(g_0) - S^{N\! L}_t(f_0) -
    \mathcal{L}S^{N\! L}_t[f_0](g_0-f_0) \right\|_{M^1} 
    \\ \nonumber
&&\qquad\qquad \le e^{C \, - {\lambda \over 2} \, t }  \,
   \sqrt{ \max \left\{ M_{k_1}(f_0), \, M_{k_1}(g_0) \right\} } \, 
   \| g_0 - f_0 \|_{M^1} ^{1+\eta}.
\end{eqnarray}
\end{lem}


Note that \eqref{estim:dt-ht-infty} implies {\bf (A4)} with
$T=\infty$, $P_{\GG_2} = P_{\GG_1}$.

From \cite{AlCaGaMo}, there exists some constants $z$, $Z$ (only
depending on the collision kernel) such that
$$
\sup_{t \ge 1} \| f_t + g_t + h_t \|_{L^1_{m_{2z}}} \le Z, \qquad
m_z(v) := e^{2 \, z \, |v|}.
$$
We also know from~\eqref{eq:HSdecay} that (possibly increasing $Z$)
$$
\forall \, t \ge 1 \qquad \| f_t - M \|_{L^1_{m_{2z}}} + \| g_t - M
\|_{L^1_{m_{2z}}} \le 2 \, Z \, e^{-\lambda \, t},
$$
where $M := M_{f_0} = M_{g_0}$ stands for the normalized Maxwellian
associated to $f_0$ and $g_0$ and 
\begin{equation}\label{ineq:SpectralGapL1LL}
\left\| e^{\LL \, t} \right\|_{L^1_{m_z}} \le 
C \, e^{- \lambda \, t}, \quad \LL h := 2 \, Q(h,M). 
\end{equation}
We  write 
\[
\partial_t (f_t-g_t) = Q(f_t-g_t,f_t+g_t) = \LL(f_t-g_t) +
Q(f_t-g_t,f_t-M) + Q(f_t-g_t,g_t-M)
\]
and we deduce for
\[
u(t) := \left\| f_t-g_t \right\|_{M^1_{m_z}}
\]
the following differential inequality (starting at some time $T_0$)
\[
u(t) \le e^{-\lambda \, (t-T_0)} \, u(T_0) + C \, \int_{T_0} ^t e^{-\lambda \,
  (t-s)} \, \left\| Q(f_s-g_s,f_s-M) + Q(f_s-g_s,g_s-M)
\right\|_{M^1(m_z)} \, ds
\]
(this formal inequality and next ones can easily be justified
rigorously by a regularizing proceedure and using a uniqueness result
for measure solutions such as \cite{Fo-Mo,EM}). Therefore we obtain
\[ 
u(t) \le e^{-\lambda \, (t-T_0)} \, u(T_0) + C \, \int_{T_0} ^t e^{-\lambda \,
(t-s)} \, \left\| (f_s-M) + (g_s-M) \right\|_{M^1(\langle v \rangle \,
m_z)} \, \left\| f_s - g_s \right\|_{M^1(\langle v \rangle \, m_z)}\,
ds.
\]
By using the control of $M^1(\langle v \rangle \, m_z)$ by
$M^1(m_{2z})$, the decay control~\eqref{eq:HSdecay} and the trivial
estimate 
\[
e^{-\lambda \, (t-T_0)}  \le e^{-\frac{\lambda}2 \, (t-T_0)} \,
e^{-\frac{\lambda}2 \, (s-T_0)}
\]
we get
\[
u(t) \le e^{-\lambda \, (t-T_0)} \, u(T_0) + C \, e^{-\frac{\lambda}2 \, (t-T_0)} \,
\int_{T_0} ^t e^{- \frac{\lambda}2 \, (s-T_0)} \, \left\| f_s - g_s
\right\|_{M^1(\langle v \rangle \, m_z)}\, ds.
\]

We then use the following control for any $a>0$:
\begin{eqnarray*} 
&&\| f-g \|_{M^1_{\langle v \rangle \, m_z}}
= \int |f-g| \, \langle v \rangle \, e^{z \, |v|}  \\
&&\qquad\le a \int_{|v| \le a} |f-g| \, e^{z \, |v|} + e^{-
  z \, a} \int_{|v| \ge a} (f+g) \, e^{2 \, z \, |v|}
\\
&&\qquad\le a \, u + e^{- z \, a} \, Z.
\end{eqnarray*}

Hence we get
\[
\| f-g \|_{M^1_{\langle v \rangle \, m_z}} \le 
\left\{ 
\begin{array}{ll}
 u + e^{- z } \, Z \le (1 + Z) \, u \quad \hbox{when}
\quad u \ge 1, \quad \hbox{(choosing } a := 1 \hbox{)}
\vspace{0.3cm} \\
{2 \over z} \, |\log \, u | \, u + u \, Z
\quad \hbox{when} \quad u \le 1\quad \hbox{(choosing } - {z \over
  2} a := \log \, u \hbox{)}
\end{array}
\right.
\]
and we deduce 
\[
\| f-g \|_{M^1_{\langle v \rangle \, m_z}} \le K \, u \, \left(1 + (\log
u)_-\right), \qquad K := 1 + {2 \over z} + Z.
\]

Then for any $\delta >0$ small, we conclude, by choosing $T_0$ large
enough, with the following integral inequality
\begin{equation}\label{eq:integralu}
u(t) \le e^{-\lambda \, (t-T_0)} \, u(T_0) + \delta \, e^{-\frac{\lambda}4 \, (t-T_0)} \,
\int_{T_0} ^t e^{- \frac{\lambda}2 \, (s-1)} \,  u_s \, 
\left(1 + \left(\log u_s\right)_-\right) \, ds.
\end{equation}

Let us prove that this integral inequality implies 
\begin{equation}\label{eq:stab-tps-long-HS}
\forall \, t \ge 1, \quad 
u(t) \le C \, e^{-\frac{\lambda}4 \, t} \, u(T_0)^{1-\delta}.
\end{equation}

Consider the case of equality in~\eqref{eq:integralu}. Then we have
$u(t) \ge e^{-\lambda \, (t-1)} \, u(1)$ and therefore 
\[
\left(1 + \left(\log u_s\right)_-\right) \le \left(1 + 
\left(\log u(T_0)\right)_- + \lambda \, (s-T_0) \right).
\]
We then have
\begin{multline*}
U(t) := \int_{T_0} ^t e^{- \frac{\lambda}2 \, (s-T_0)} \, u_s \,
\left(1 + \left(\log u_s\right)_-\right) \, ds \le \int_{T_0} ^t e^{-
  \frac{\lambda}2 \, (s-T_0)} \, u_s \, \left(1 + \left(\log
    u(T_0)\right)_- + \lambda \, (s-T_0) \right) \, ds \\ \le \left(1 +
  \left(\log u(T_0)\right)_- \right) \,  \int_{T_0} ^t e^{-
  \frac{\lambda}4 \, (s-T_0)} \, u_s \, ds
\end{multline*}
and we conlude the proof of the claimed inequality
\eqref{eq:stab-tps-long-HS} by a Gronwall-like argument.

Then estimate \eqref{estim:dt-infty} follows by choosing $\delta$
small enough (in relation to $\eta$) and then connecting the last
estimate \eqref{eq:stab-tps-long-HS} from time $T_0$ on together with
the previous finite time estimate \eqref{estim:dt} from time $0$ until
time $T_0$. 

Then estimates \eqref{estim:ht-infty} and \eqref{estim:dt-ht-infty} are
proved in the same way by using the equations
\[
\partial_t h_t = \LL(f_t-g_t) + Q(h_t,f_t-M)
\]
(which is even simpler than the equation for $f_t-g_t$) and
\[
\partial_t \zeta_t = 2 \, \LL(h_t) +
Q(\zeta_t,f_t-M)+Q(\zeta_t,g_t-M) + Q(h_t,d_t).
\]

\subsection{Proof of (A5) uniformly in time}  
Let us prove that for any $\bar z, \MM_{\bar z} \in (0,\infty)$ there
exists some continuous function $\Theta : \R_+ \to \R_+$, $\Theta(0) =
0$, such that for any $f_0$, $g_0 \in P_{m_{\bar z}}(\R^d)$, with
$m_{\bar z}$ defined in the previous section, such that 
\[
\left\| f_0 \right\|_{M^1_{m_{\bar z}}} \le \MM_{\bar z}, \quad 
\left\| g_0 \right\|_{M^1_{m_{\bar z}}} \le \MM_{\bar z},
\]
 there holds
\begin{equation} \label{estim:W1dt} 
\sup_{t \ge 0} W_t \le \Theta (W_0), 
\end{equation}
where $W_t$ stands for the Kantorovich-Rubinstein distance
\[
W_t = W_1\left( S^{N \! L}_t(f_0), \, S^{N \! L}_t(g_0) \right).
\]
As we shall see, we may choose
\begin{equation}\label{estim:W1dt}
 \Theta (w) := \bar\Theta  \min \left\{ 1, \, 
e^{- (\log (w/\theta_0))_-^{1/2}} \right\},
\end{equation}
for some constants $\bar\Theta, \, \theta_0 > 0$ (only depending on
$\bar z$ and $\MM_{\bar z}$).  

We start with 
\begin{equation}\label{ineq:wt2eme}
\forall \, t \ge 0 \quad W_t \le \left\| (f_t  - g_t) |v|
\right\|_{M^1} 
\le  {1\over2} \left\| (f_t  + g_t) \langle v \rangle^2 \right\|_{M^1} 
= 1 + \EE =: \bar \Theta.
\end{equation}

Let us improve this inequality for small value of $W_0$. On the one
hand, it has been proved in \cite[Theorem~2.2 and
Corollary~2.3]{Fo-Mo} that
\begin{equation}\label{ineq:dwlogw-} 
W_t \le W_0 + K \,
\int_0 ^t W_s \, \left(1+(\log W_s)_-\right) \, ds,
\end{equation} 
for some constant $K$. 
Whenever $W_t \le 1/2$, the integral inequality \eqref{ineq:dwlogw-}
implies (possibly increasing the constant $K$)
$$
W_t \le W_0 + K \, \int_0 ^t W_s \, \left(\log W_s\right)_- \, ds.
$$
From the Gronwall lemma we deduce
\begin{equation}\label{ineq:wtFirst}
W_t \le \left(W_0\right)^{\exp({-Kt})} 
\quad\hbox{whenever}\quad W_t \le 1/2.
\end{equation}
On the other hand, from~\cite{Bobylev97} and~\cite{MouhotCMP}, there
exists $\lambda_2, Z > 0$, $z \in (0,\bar z)$ such that
\[
\forall \, t \ge 0 \qquad \| f_t  - M_{f_0} \|_{L^1_{m_{z}}} +  \| g_t  - M_{g_0}
\|_{L^1_{m_{z}}}   \le  Z \, e^{-\lambda \, t},
\]
where $M_{f_0}$ and $M_{g_0}$ stand again for the normalized
Maxwellian associated to $f_0$ and $g_0$. Denoting by $u_{f_0}$ and $u_{g_0}$ the mean velocity of $f_0$ and $g_0$, by 
$\theta_{f_0}$ and $\theta_{g_0}$ the temperature associated to $f_0$ and $g_0$, and by $\EE_{f_0}$ and $\EE_{g_0}$ the energy associated to $f_0$ and $g_0$,  there also holds
\begin{eqnarray*} 
W_1( M_{f_0},M_{g_0} ) 
&\le& C_s \, \left( \|ÊM_{f_0} - M_{g_0} \|_{H^{-s}}^2 \right)^\eta
\\
&\le& C_s \, \left( \int_{\R^d}  
{ \left|e^{- \theta_{f_0} \, {|\xi|^2 \over 2} - i \, u_{f_0} \, \sqrt{\theta_{f_0}} \, \xi } 
- e^{- \theta_{g_0} \, {|\xi|^2 \over 2} - i \, u_{g_0} \, \sqrt{\theta_{g_0}} \, \xi } 
\right|^2 Ê\over \langle \xiÊ\rangle^{2s}}  \, d\xi\right)^\eta
\\
&\le& C_s \, \left( \int_{\R^d}  
{ |\theta_{f_0} - \theta_{g_0}|^2  \, |\xi|^4+ | u_{f_0} \, \sqrt{\theta_{f_0}}  - u_{g_0} \, \sqrt{\theta_{g_0}} |^2 \, |\xi|^2 Ê\over \langle \xiÊ\rangle^{2s}}  \, d\xi\right)^\eta
\\
&\le& C_s \, \left(  |\theta_{f_0} - \theta_{g_0}|^\eta  +  | u_{f_0} - u_{g_0} |^{2\eta} \right) 
\\
&\le& C_s \, \left(  |\EE_{f_0} - \EE_{g_0}|^\eta  +  | u_{f_0} - u_{g_0} |^{2\eta} \right) 
\\
&\le& C_s \, \left( W_2(f_0,g_0)^{2\eta } + W_1(f_0,g_0)^{2\eta} \right) 
\\
&\le& C_s \, \left( W_1(f_0,g_0)^{\eta/2 } + W_1(f_0,g_0)^{2\eta} \right)
\end{eqnarray*}
(see also \cite{Chafai-Malrieu} for more general estimates of the
Wasserstein distance between two gaussians). Gathering these two
estimates, we deduce
\begin{equation}\label{ineq:wt2eme}
\forall \, t \ge 0 \qquad W_t \le Z_1 \, e^{-\lambda \, t} + Z_2 \, W_0^{\eta/2}.
\end{equation}

Let us (implicitly) define $T, \bar W_0 \in (0,\frac{1}{4^{2/\eta}})$ by
\[
Z_1 \, e^{-\lambda \, T} = {1 \over 4} \quad\hbox{and next}\quad 
\left(\bar W_0\right)^{\exp({-K\, T})} = {1\over2}.
\]
Then, for any $W_0 \in (0,\bar W_0)$ we have 
$$
\forall \, t \ge 0 \qquad W_t \le \min \left\{
  \left(W_0\right)^{\exp({-Kt})}; \, Z_1 \, e^{-\lambda \, t} \right\} + Z_2 \, W_0^{\eta/2}.
$$

Let us search for $t^*$ such that the two functions involved in
the minus function are equal:
$$
\phi(t^*) := \left|\log W_0 \right| \, \exp({-K \, t^*}) = \lambda \, t^*
- \log Z =: \psi(t^*).
$$
The time $t^*$ is unique because $\phi$ is decreasing while $\psi$ is
increasing (and it exists for $\log \bar W_0 \le \log Z$). Choosing
\[
t^\sharp := {1 \over K} \, \log \left( \left|\log
    W_0\right|^{1/2}\right),
\]
we find
$$
\phi (t^\sharp) =  \left|\log W_0\right|^{1/2} \ge {\lambda \over K}
\, \log \left( \left|\log W_0\right|^{1/2}\right) - \log Z = 
\psi(t^\sharp),
$$
at least whenever $W_0 \in (0,W_0^\sharp]$, with $W_0^\sharp \in
(0,\bar W_0]$ small enough. As a consequence, for any $W_0 \in
(0,W_0^\sharp]$ we have
$$
\forall \, t \ge 0, \quad W_t \le 2 \,  \phi(t^\sharp) = 2 \, e^{- \left|\log
    W_0\right|^{1/2}}.
$$

This concludes the proof. 


\section{Extensions and complements}
\label{sec:extensions}

In this section first we generalize the chaos propagation
results on the Boltzmann (maxwellian or hard spheres molecules) to the
case where the limit $1$-particle distribution is not compactly
supported. This shall rely on a construction due to Kac of the
$N$-particle initial data together with careful study of the
dependency of the constants in terms of moments of the data. Second,
as a corollary, we use our global in time results to give a new method
for studying chaotic convergence of the $N$-particle equilibria towards the
limit $1$-particle equilibrium. The old question of connecting the
long-time behavior was raised by Kac and it motivated its whole study
of chaos propagation for particle systems. In the case of classical
gas dynamics, we thus recover a well-known computational results
(namely the marginals of the constant probabilities on
$\sqrt{N} \mathbb{S}^{Nd-1}$ converges to products of gaussians)
without any explicit computations, only using the properties of the
Boltzmann semigroups. This new method shall prove highly useful in the
context of granular gases where the steady states or homogeneous
cooling states are not explicitly known.

\subsection{Construction of chaotic initial data $f^N_0$ with
  prescribed energy}
\label{subsec:pN0}

 For the sake of completeness, let us recall, following
 \cite{Kac1956}, how to construct a $f_0$-chaotic sequence of initial
 data $f^N_0$ ({\it i.e.}, which has the {\it ``Boltzmann's
   property''} in the words of Kac).

\begin{lem}\label{lem:Bproperty}
 Consider $f_0 \in P(\R^d)$ with finite energy  $M_{m_e} (f_0) := M_2(f_0) = \EE \in (0,\infty)$ and which fulfills 
 the following moments conditions $M_{m_i}(f_0) = M^{N\!L}_{0,m_i}< \infty$, $i=0, 1, 3$, 
for some radially symmetric and increasing weight functions $m_i$, $i=1,3$, and $m_0(x) := \exp (a \, |x|^2)$, $a > 0$.
Then for any given increasing sequence $(\alpha_N)_{N \ge 1}$ (which
  increases as slow as we wish in general and may be chosen constant
  when $f_0$ has compact support), there exists a sequence $f^N_0 \in
  P(\R^{dN})$, $N \ge 1$, such
  that 
  \begin{itemize}
  \item[(i)] The sequence $(f^N_0)_{N \ge 1}$ is $f_0$-chaotic.
  \item[(ii)] Its support satisfies 
    \[
    \hbox{{\em supp}} \, f^N_0 \subset S^{Nd-1} \left(\sqrt{N \,
      \EE}\right) := \left\{ V \in \R^{dN}; \,\, M^N_2 (V) = \EE \right\} \subset  \mathbb{E}_N.
  \] 
  \item[(iii)] It satisfies the following integral moment bound based on
  $m_2$:
    \[ 
    \forall \, N \in \N^*, \quad \left\langle f^N_0, M_{m_3}^N
    \right\rangle \le \mbox{Cst} \, \left(M^{N\!L}_{0,m_3}\right).
    \]
   \item[(iv)] It satisfies the following moment bound {\em on the
      support} involving $m_3$:
    \[
    \hbox{{\em supp}} \,f^N_0 \subset \KK_N := \left\{ V \in \R^{dN}; \,\,
    M^N_{m_3} (V) \le \alpha_N \right\}.
    \]
\end{itemize}
\end{lem} 

\noindent {\it Sketch of the proof of Lemma~\ref{lem:Bproperty}.}
We essentially recall briefly the key arguments presented in
\cite[Section 5 ``Distributions having Boltzmann's
property'']{Kac1956} and check that the moment conditions required in
the sequel of our paper can be satisfied.  For the sake of simplicity,
we assume with no loss of generality that the energy $\EE = 1$.  We
restrict to the case when $f_0 \in P(\R^d) \cap C(\R^d)$ and we refer
to \cite{CCLLV} for the relaxation of this condition.

Since $f_0 \in C(\R^d)$, we can define
$$
f^N_0(V) := { \prod_{j=1}^N f_0 (v_j) \over F_N(\sqrt{N})}
\Bigg|_{\mathbb{S}^{Nd-1}(\sqrt{N} \, )} \quad \mbox{ with } \quad F_N(r) :=
\int_{ \mathbb{S}^{Nd-1} (r) } \prod_{j=1}^N f_0 (v_j) \, d\omega,
$$
so that (ii) holds. 

From the gaussian moment bound $M_{m_0}(f_0) < \infty$, we obtain from
\cite{Kac1956} that there exists some constants $C> 0$ such that
$$
F_N(\sqrt{N}) \, \sim \, C \, N,
$$
and for $\varphi(v_1,\dots,v_\ell)$, $\ell \le N$: 
\[
\int_{ \mathbb{S}^{Nd-1} (\sqrt{N}\, ) } \varphi(v_1,\dots,v_\ell) \,
\prod_{j=1}^N f_0 (v_j) \, dS(V) \xrightarrow[N \to \infty]{} C \, N
\, \int_{\R^d} \varphi(v_1,\dots,v_\ell) \, df_0 (v_1) \dots
df_0(v_\ell)
\]
which proves the chaos.

We then deduce (i) thanks to \cite[Proposition 2.2]{S6}. Point (iii)
is just a consequence of the above asymptotic with the choice $\varphi
= m_1$, $\psi = 1$. When we assume furthermore that $f_0$ is
compactly supported, say supp$\, f_0 \subset [-A,A]^d$, we deduce
$\hbox{supp} \, f^N_0 \subset \{ V \in \R^{Nd}, \,\, M^N_{m_3} (V) \le
m_3(A) \}$ and (iv) holds.

In the non compactly supported case, for any $k \in \N^*$ and for any
constant $A_k$ we define $f_{0,k} := f_0 \, {\bf 1}_{|v| \le
  A_k}$. It is clear that for any $k \in \N^*$, there exists
$f^N_{0,k}$ such that $f^N_{0,k}$ is $f_{0,k}$-chaotic such that (ii)
and (iii) hold and $ \hbox{supp} \,f^N_{0,k} \subset \{ V \in \R^{dN};
\,\, M^N_{m_3} (V) \le m_3(A_k) \}$. For any given sequence
$(\alpha_N)$ which tends to infinity, we define $k_N$ in such a way
that $m_3(A_{k_N}) = \alpha_N$ so that $k_N \to \infty$ when $N \to
\infty$. We then easily verify that $f^N_0 := f^N_{0,k_N}$ satisfies
the properties (i)--(iv).  \qed

\subsection{Chaos propagation without compact support for the
  Boltzmann equation}
\label{sec:extension-non-cpct}

We may relax the compactly support condition in Theorem~\ref{theo:tMM} 
and Theorem~\ref{theo:HS} thanks to Lemma~\ref{lem:Bproperty}. 
We assume that 
$$
M_{m_0} (f_0) := \int_{\R^d} e^{a \, |v|^2}Ê\, f_0 (dv) < \infty
$$
for some $a \in (0,\infty)$ and we define $f^N_0$ as in Lemma~\ref{lem:Bproperty}. Instead of \eqref{eq:cvgabstract1} in 
 Theorem~\ref{theo:abstract} we have the following bound. 
 For any increasing sequence $\alpha_N \to \infty$, for any $ \varphi  \in  (\FF_1 \cap \FF_2 \cap \FF_3)^{\otimes \ell}$,
 there exists a constant $C_{\ell,\varphi}$ (independent of $N$) such
  that for any $N \in \N^*$, with $N \ge 2 \ell$,
    \begin{eqnarray}
  \label{eq:cvgabstractWithoutCompactS}
  &&\quad \sup_{[0,T)}\left| \left \langle \left( S^N_t(f_0 ^N) - \left(
          S^{N \! L}_t(f_0) \right)^{\otimes N} \right), \varphi 
    \right\rangle \right|  \
  \\ \nonumber 
  &&\quad 
  \le C_{\ell,\varphi} \, \Bigg[ \frac{1}{N} 
  + C^N_{T,m_1} \, C_T^\infty \, \varepsilon_2(N) 
   +  \Theta_{\alpha_N,T}  \left(  \WW_{\!\! \hbox{ \small {\em dist}}_{\GG_3}} \!\! 
    \left(\pi^N_P f^N_0,\delta_{f_0}\right) \right) \Bigg].
  \end{eqnarray}
 By choosing 
  $(\alpha^N)$ appropriately, we will deduce from \eqref{eq:cvgabstractWithoutCompactS} that 
 $S^N_t(f_0 ^N)$ is  again $S^{NL}_t(f_0)$-chaotic uniformly in time. 
  
  \smallskip
In the case of the Boltzmann model for Hard Spheres or true Maxwell Molecules, we can take
$\hbox{dist}_{\GG_3} = W_1$ the usual Monge-Kantorovich-Wasserstein
distance in $P(\R^d)$. We claim that 
\begin{eqnarray}\label{eq:cvgceWW1}
 \WW_{W_1} \!\! 
    \left(\pi^N_P f^N_0,\delta_{f_0}\right) 
    \mathop{\longrightarrow}_{N\to\infty} 0.
\end{eqnarray}
First, thanks to \cite[Proposition 2.2]{S6} and the fact that $f^N_0$ is $f_0$-chaotic from Lemma~\ref{lem:Bproperty},
we deduce that $ \pi^N_P f^N_0$ converges to $\delta_{f_0}$ in the weak sense in $P(P(\R^d))$ (that means taking
duality product with functions of $C(P(\R^d))$). Next, thanks to \cite[Theorem 7.12]{VillaniTOT}, \eqref{eq:cvgceWW1}
boils down to prove  that 
\begin{eqnarray}\label{eq:tightWW1}
\lim_{R \to \infty} \sup_{N \in \N^*} \int_{W_1(\rho,f_0) \ge R} W_1(\rho,f_0) \, \pi^N_P f^N_0(\rho) = 0,
\end{eqnarray}
which will be a straightforward consequence of the following bound 
\begin{eqnarray}\label{eq:bddWW1}
\sup_{N \in \N^*} \int_{E^N} [W_1(\mu^N_V,f_0)]^2 \, f^N_0(dV) < \infty.
\end{eqnarray}
Finally, in order to get \eqref{eq:bddWW1}, we infer that from \cite[Theorem 7.10]{VillaniTOT}
\begin{eqnarray*}
 [W_1(\mu^N_V,f_0)]^2 
 &\le& \|Ê\mu^N_V - f_0 \|_{M^1_1}^2 
 \\
&\le& 2 \|Ê\mu^N_V  \|_{M^1_1}^2 + 2 \, \| f_0 \|_{M^1_1}^2  
\\
&\le& 2 \, [ÊM^N_1(V)]^2 + 2 \, \| f_0 \|_{M^1_1}^2 
\le  2 \, M^N_2(V)+ 2 \, \| f_0 \|_{M^1_1}^2.
\end{eqnarray*}
That implies 
$$
\int_{E^N} [W_1(\mu^N_V,f_0)]^2 \, f^N_0(dV) \le 2 \, \| f_0 \|_{M^1_1}^2 
+ 2 \, \langle f^N_0,M_2 \rangle,
$$
which, together with (ii) in Lemma~\ref{lem:Bproperty}, ends the proof of \eqref{eq:bddWW1} and then 
of \eqref{eq:cvgceWW1}.

\smallskip
From the fact that the $\Theta_{T,A}$ functions exhibited in Theorem~\ref{theo:tMM} 
and Theorem~\ref{theo:HS} are independent of $T$ and satisfy $\Theta_{T,A}(x) = \Theta_{A}(x) \to 0$
when $x \to 0$ for any fixed $A \in (0,\infty)$ we may build (thanks to a diagonal process) a 
sequence $(\alpha^N)$ such that 
$$
\Theta_{\alpha^N} \left (  \WW_{W_1} \!\! 
    \left(\pi^N_P f^N_0,\delta_{f_0}\right) \right)
    \mathop{\longrightarrow}_{N\to\infty} 0.
$$
Coming back to \eqref{eq:cvgabstractWithoutCompactS} we obtain that
$S^N_t(f_0 ^N)$ is  $S^{NL}_t(f_0)$-chaotic uniformly in time.

\subsection{Chaoticity of the sequence of steady states}
\label{sec:entropies}







\begin{theo}[Abstract fluctuation estimate in the infinite time]\label{theo:abstractInfty}
  Consider a sequence of initial datum $f^N_0$ which satisfies {\bf
  (A1)} (with $C^N_{0,m_3} = \alpha^N$ may depend on $N$) and is $f_0$-chaotic with $f_0 \in P_{\GG_1} \cap P_{\GG_3}$.    
   Assume moreover that the assumptions of theorem~\ref{theo:abstract} hold
  with $T = \infty$.  Assume finally that $f^N_t \wto \gamma^N$ when $t\to\infty$ in the
  weak sense of measures in $P(E^N)$ as well as $f_t \wto \gamma$ when $t\to\infty$  in
  the weak sense of measures in $P(E)$. 
  For any $\ell \in \N^*$ and $ \varphi  \in  (\FF_1 \cap \FF_2 \cap \FF_3)^{\otimes \ell}$,
 there exists a constant $C_{\ell,\varphi}$ such that for any $N \in \N^*$, with $N \ge 2 \ell$,
  we have
  \begin{eqnarray}
  \label{eq:cvgabstractinfty1}
  &&\quad \left| \left \langle \Pi_\ell \gamma^N 
      - \gamma^{\otimes\ell} , \varphi 
    \right\rangle \right|  \le \\
   &&\quad  \nonumber
   \le C_{\ell,\varphi} \, \Bigg[ \frac{1}{N} 
  + C^N_{\infty,m_1} \, C_\infty^\infty \, \varepsilon_2(N) 
   +  \Theta_{C^N_{0,m_3},\infty}  \left(  \WW_{\!\! \hbox{ \small {\em dist}}_{\GG_3}} \!\! 
    \left(\pi^N_P f^N_0,\delta_{f_0}\right) \right) \Bigg].
  \end{eqnarray}
  As a consequence, $\gamma^N$ is $\gamma$-chaotic. 
  \end{theo}

  The proof of that result is trivial: we just have to apply
  Theorem~\ref{theo:abstract} and to pass to the limit in the left
  hand side of the inequality \eqref{eq:cvgabstract1} in order to
  gat \eqref{eq:cvgabstractinfty1}. Arguing as in section~\ref{sec:extension-non-cpct}
  we deduce $\gamma^N$ is $\gamma$-chaotic (whenever $C^N_{0,m_3} = \alpha^N$
  grows slowly enough). 
  
  \medskip
  The application
  to the Boltzmann equation is the following. Consider a sequence of
  initial data $f^N_0$ such that supp$\, f^N_0 \subset
  \mathbb{S}^{Nd-1}(\sqrt{N})$, and such that $f^N_0$ is $f_0$-chaotic with
  $\int |v|^2 f_0 = 1$, $\int v_i \, f_0 = 0$ for any $i=1, \dots, d$. On
  the one hand, we know (see \cite{Kac1956,CCLLV}) that $f^N_t$
  converges in the large time asymptotic to $\gamma^N$, the uniform
  distribution on the sphere $\mathbb{S}^{Nd-1}(\sqrt{N})$ (that holds in
  $L^2(\gamma^N)$ with rate $\exp(-\lambda_N \, t)$ whenever $f^N_0
  \in L^2(\gamma^N)$). We also know (see \cite{MouhotCMP} and the
  references therein for the hard sphere case and
  \cite{T1,CGT} or section~\ref{sec:modelEBbounded} for the (true) Maxwell Molecules
  case) that $f_t$ converges in the large time
  asymptotic to $\gamma$, the normalized Gaussian. As a consequence,
  we get \eqref{eq:cvgabstractinfty1}. That result may seem to be
  trivial, and it is in some sense, because an explicit computation
  (which go back at least to Mehler in 1866) yields the same
  result. However, our proof it is not based on an explicit
  computation nor a variationnal/entropy optimization principle. The
  consequence is that it applies to many more situations, in
  particular in the case of some dissipative Boltzmann equation
  (linked to the Granular media), see \cite{MMW}.

  \subsection{On statistical solutions and the non uniqueness of its
    steady states}
\label{sec:entropies}

Let us consider the $N$-particles system associated to the Boltzmann collision process that we do not write in dual for as we have done before, it writes
\begin{equation}\label{eq:MasterN}
\partial_t f^N = {1 \over N} \sum_{i < j }  \int_{S^{d-1}} B(|v_i-v_j|,\cos \theta) \, 
\Big[ f^N(\dots,v'_i, \dots, v'_j,\dots) - f \Big] \, d\sigma.
\end{equation}
We want to describe how the BBGKY (Bogoliubov, Born, Green, Kirkwood and Yvon) method introduced 
to derive Boltzmann's equation from Liouville's equation applies in our simpler space homogeneous context. 
Let us thus also introduce the $k$-th marginal:
$$
f^N_\ell (v_1, ..., v_\ell) := \int_{\R^{d (N-k)}} f^N(v_1, ..., v_\ell, w_{\ell+1}, ..., w_N) \, dw_{\ell+1} \, ...dw_{N}
$$
Integrating the master equation \eqref{eq:MasterN} leads to 
\begin{eqnarray*}
\partial_t f^N_\ell
&=& {1 \over N} \sum_{i,j \le \ell} \ZZ^N_{ij}  \qquad\qquad = \OO (\ell^2/N) \\
&&+ {1 \over N} \sum_{i \le \ell < j}  \ZZ^N_{ij} \\
&&+ {1 \over N} \sum_{i, j > \ell}  \ZZ^N_{ij}  \qquad\qquad = 0,
\end{eqnarray*}
with 
$$
\ZZ^N_{ij} := \int_{\R^{d \, (N-\ell-1)}} \int_{S^{N-1}} B\Big[ f^N(\dots,v'_i, \dots, v'_j,\dots) - f^N \Big] \, d\sigma \, dv_{\ell+1} \, ... \, dv_N .
$$
Only the second term does not vanish in the limit $N\to \infty$, so that assuming that $f^N_\ell \to \pi_\ell$, we find that $(\pi_\ell)$ is a solution to the infinite dimensional system of linear equation (the Boltzmann equation for a system of an infinite number of particles or the statistical Boltzmann equation)
\begin{equation}\label{eq:BBGKY}
\partial f_\ell = A_{\ell+1} (\pi_{\ell+1})
\end{equation}
with  $\pi_\ell = \pi_\ell(t,v_1, ..., v_\ell) \ge 0$ and
$$
V \in \R^{d\ell} \mapsto A_{\ell+1} (\pi_{\ell+1}) (V) = \sum_{j=1}^\ell \int_{S^{d-1} \times \R^d} \Bigl\{ \pi_{\ell+1}(V'_j) - \pi_{\ell+1}(V) \Bigr\} \, b( \sigma \cdot (v_j - v_{\ell+1})) \, dv_{\ell+1} d\sigma,
 $$
with $V'_j = (v_1, ..., v'_j, ..., v_\ell, v'_{\ell+1})$, $v'_j = v'(v_j,v_{\ell+1},\sigma)$, $v'_{\ell+1}= v_*'(v_j,v_{\ell+1},\sigma)$ where
vectors  $v'$ and $v'_*$ are defined by \eqref{eq:rel:vit}. 

\begin{lem}\label{sec6:SolStat} There exists a non chaotic stationary solution to the statistical Boltzmann equation. In other words, 
there exists $\pi \in P(P(\R^d))$ such that $\pi \not= \delta_p$ for some $p \in P(\R^d)$ and $A_{\ell+1} (\pi_{\ell+1}) = 0$ for any $\ell \in \N$. 
\end{lem}

\noindent
{\bf Proof of Lemma~\ref{sec6:SolStat}. } It is clear that any function on the form
$$
V \in \R^{d(\ell+1)} \, \mapsto \, \pi_{\ell+1} (V) = \phi (|V|^2)
$$
is a stationary solution got equation \eqref{eq:BBGKY}, that is $\CC_{\ell+1} (\pi_{\ell+1}) = 0$. Now we define, with $d=1$ for the sake of simplicity, the sequence 
$$
V \in \R^\ell \mapsto \pi_\ell (V) = {c_\ell \over (1+ |V|^2)^{m+\ell/2}}
$$
with  $c_1$ such that  $\pi_1$ is a probability measure and $c_2 = c_1 \, \alpha_2$ with $\alpha_2$ chosen in the following way:
\begin{eqnarray*}
\int_{\R} {\alpha_2 \over (1+v^2 + v_*^2)^{m+1}} \, dv_*
&=&{\alpha_2 \over (1+v^2)^{m+1}} \int_{\R} {1 \over (1 + { v_*^2 \over 1+v^2} )^{m+1}} \, dv_* \\
&=& {\alpha_2 \over (1+v^2)^{m+1/2}} \int_{\R} {1 \over (1 + w_*^2)^2} \, dw_* = 
{1 \over (1+v^2)^{m+1/2}}.
\end{eqnarray*}
By an iterative process we may chose the constants $c_\ell$ in such a way that  $\pi_\ell$ is a solution to $A_\ell (\pi_\ell) = 0$ (because it is a function of the energy) and satisfies the compatibility condition:
$$
\pi_\ell(V) = \int_\R \pi_{\ell+1}(V,v_*) \, dv_*.
$$ 
We have exhibit a solution which is not chaotic. 
\qed

\smallskip
We come back to the abstract setting. We start with the $N$-particles system equation, 
$$
\partial_t f^N = A^N \, f^N,
$$
that we write in dual form
\begin{eqnarray*}
\partial_t \langle  f^N_\ell, \varphi \rangle &= &
\partial_t \langle  f^N, \varphi \otimes {\bf 1}^{N-\ell} \rangle \\
&=&  \langle  f^N, G^N ( \varphi \otimes {\bf 1}^{N-\ell}) \rangle 
=  \langle  f^N_{\ell+1}, G^N_{\ell+1} ( \varphi ) \rangle.
\end{eqnarray*}

\begin{lem}\label{lem:BBGKY} Assume that 

{\bf (A6)} $(f^N_\ell)$ is tight in $P(E^\ell)$ for any $\ell \ge \N$ (or equivalently, $f^N$ is  tight in $P(P(E))$);

{\bf (A7)} $G^N_{\ell+1} \varphi \to G^\infty_{\ell+1} \varphi$ when $N \to \infty$, for any fixed $\varphi \in C_b(E^\ell)$. 

Then, up to extraction a subsequence, $(f^N)$ converges (in the sense of any $\ell$-th marginals) to a solution $\pi = (\pi_{\ell})  \in P(P(E))$
to the infinite Hierachy
$$
\partial_t \pi = A^\infty \pi \quad\hbox{in}\quad P(P(E)),
$$
which simply means 
$$
\partial_t \langle \pi_\ell,\varphi \rangle  = \langle\pi_{\ell+1}, G^\infty_{\ell+1} \varphi \rangle  \quad\hbox{for any }\quad \ell \in \N^*.
$$

\end{lem}

\subsection{Uniqueness of statistical solutions and chaos}
\label{sec:tensor}

Assuming {\bf (A2i) } and 
\begin{itemize}
\item[{\bf (A2iii')}] $[0,\infty) \to P_{\GG_1}$, $t \mapsto S^{N\!L}_t f$ uniform continuously for any $f \in P_{\GG_1}$,
\end{itemize}
 which is a weak version of {\bf (A2iii)}, we have that $(S^{N\!L}_t)$ is a $c_0$-semigroup. As it has proved in step 1 of Lemma~\ref{lem:H0}, for any $\Phi \in C_b(P_{\GG_1},\R)$ we may define $T^\infty \Phi$ by 
$$
(T^\infty \Phi) (f) = \Phi(S^{N\!L}_t f),
$$
and we build in that way a $c_0$-semigroup $(T^\infty_t)$ on $C_b(P_{\GG_1},\R)$. The Hille-Yosida theory imply that there exists an closed operator $G^\infty$ with dense domain $\hbox{dom}(G^\infty)$ in $C_b(P_{\GG_1},\R)$ so that $(T^\infty_t)$ is the semigroup associated to the generator $G^\infty$. 

\smallskip
Now, on the one hand, for any $\pi_0 \in P(P_{\GG_1})$ we may define the semigroup $(S_t)$ on $P(P_{\GG_1})$ and the flow $(\bar \pi_t)$ by setting $\bar \pi_t = S^\infty_t \pi_0$ and (duality formula)
$$
\forall \, \Phi \in C_b (P_{\GG_1};\R) \qquad \langle S^\infty_t \pi_0,\Phi \rangle  = \langle \pi_0 ,  T^\infty _t \Phi \rangle.
$$

Remark (see \cite{ArkerydCI99}) that $S^\infty_t
\pi_0 \in (C_b(P(V))'_+ \not= P(P(E))$.  Under the additional
assumption that $P_{\GG_1,a}$ is compact for any $a$ (that is true in
our application cases when $\| \cdot \|_{\GG_1}$ metrizes
the weak measures topology) that relation defined a unique probability
$S^\infty_t \pi_0 \in P(P_{\GG_1})$, and again $(S^\infty_t)$ is a
$c_0$-semigroup on $P(P_{\GG_1})$. 
 
\smallskip
On the other hand, we say that $\pi_t \in C(\R_+; P(P_{\GG_1}))$ is a solution to the equation 
\begin{equation}\label{sec6:dtpi=A}
\partial_t \pi_t = A^\infty \, \pi_t,
\end{equation}
if for any $\Phi \in  \hbox{dom}(G^\infty)$ there holds
$$
{d \over dt} \langle \pi_t  , \Phi \rangle = \langle \pi_t, G^\infty \Phi \rangle \quad \hbox{in} \,\, \DD'([0,\infty)).
$$

\begin{theo}\label{theo:BBGKYuniq} Assume that {\bf (A2)} and {\bf (A4)} hold, as well as that $C^1(P_{\GG_1};\R)$ is dense in $C_b(P_{\GG_1};\R)$. For any initial datum $\pi_0 \in P(P_{\GG_1})$ the flow $\bar\pi_t$ is the unique solution in $C([0,\infty); P(P_{\GG_1}))$ to \eqref{sec6:dtpi=A} starting from $\pi_0$. Moreover, if $\pi_0$ is $f_0$-chaotic, then $\pi_t$ is $S^{N \! L}_t f_0$-chaotic for any $t \ge 0$. 
\end{theo}

\noindent{\sl Proof of Theorem~\ref{theo:BBGKYuniq}. } {\sl Step 1: Chaos propagation. } From Hewitt-Savage's theorem~\cite{Hewitt-Savage}, for any $\pi \in P(P(E))$ there exists a unique sequence $(\pi^\ell) \in P(E^\ell)$ such that the identities
\begin{eqnarray*}
\langle \pi^\ell , \varphi \rangle 
&=& \int_{P(E)} \langle f^{\otimes\ell},\varphi\rangle \, \pi(df) \\
&=& \int_{P(E)} R^\ell_\varphi (f) \, \pi(df) = \langle \pi, R^\ell_\varphi \rangle,
\end{eqnarray*}
hold for any $\varphi \in C_b(E)^{\otimes\ell}$. As a consequence,  if $\pi_0$ is $f_0$-chaotic,  
\begin{eqnarray*}
\langle \bar\pi^\ell_t, \varphi \rangle
&=&\langle \bar \pi_t, R^\ell_\varphi \rangle = \langle \pi_0, T^\infty_t R^\ell_\varphi \rangle = (T^\infty_t R^\ell_\varphi ) (f_0) \\
&=& R^\ell_\varphi (S^{N \! L}_t f_0)  =  \langle S^{N \! L}_t \, f_0,\varphi_1 \rangle \, ... \, \langle S^{N \! L}_t \, f_0,\varphi_\ell \rangle,
\end{eqnarray*}
which means that $\bar\pi_t^\ell = f^{\otimes\ell}_t$, or equivalently $\bar\pi_t = \delta_{f_t}$, and the statistical solution $\bar\pi_t$ is $f_t$-chaotic.

\smallskip\noindent
 {\sl Step 2: Uniqueness. } For any $t > 0$ and $n \in \N^*$ owe define $\eps := t/n$ and $t_k = \eps \, k$, $s_k = t - t_k$.  
Then for any $\Phi \in C^1_b(P_{\GG_1};\R)$ we define $\Phi_t := T^\infty_t \Phi$. The very fundamental point is
  that thanks to Lemma~\ref{lem:H0} we have $\Phi_t \in C^1_b(P_{\GG_1};\R) \subset \hbox{dom}(G^\infty)$ for any $t \ge 0$.
We write 
 \begin{eqnarray*}
 &&\langle \pi_t,\Phi \rangle -  \langle \bar\pi_t,\Phi \rangle = \langle \pi_t,\Phi \rangle -  \langle \pi_0,\Phi_t \rangle \\
 &&\qquad= \sum_{k=0}^{n-1} \left\{ \left[ \langle \pi_{t_{k+1}}, \Phi_{s_{k+1}} \rangle -  \langle \pi_{t_{k+1}}, \Phi_{s_{k}} \rangle\right]
 +\left[ \langle \pi_{t_{k+1}}, \Phi_{s_{k}} \rangle -  \langle \pi_{t_k}, \Phi_{s_{k}} \rangle\right] \right\} 
\\
&&\qquad=  \TT_1 + \TT_2 = \sum_{k=0}^{n-1} \left\{ \TT_{1,k} + \TT_{2,k} \right\}. 
 \end{eqnarray*}
On the one hand, we have
\begin{eqnarray*}
\TT_{1,k} &=& \langle \pi_{t_{k+1}}, \Phi_{s_{k+1}} - T^\infty_\eps \Phi_{s_{k+1}}\rangle 
= -  \langle \pi_{t_{k+1}}, \int_0^\eps {d \over ds} [T^\infty_s \Phi_{s_{k+1}} ] \, ds \rangle
\\
&=& -  \langle \pi_{t_{k+1}},    \int_0^\eps [G^\infty \Phi_{s_{k+1}+s} ] \, ds \rangle
= -  \int_{s_{k+1}}^{s_k}  \langle \pi_{t- [s+1,\eps]} ,  G^\infty \Phi_s \rangle \, ds,
\end{eqnarray*}
where  $[s,\eps] = [s/\eps]\, \eps$. Passing to the limit  $n\to\infty$, we get
$$
\TT_1 = - \int_0^{t} \langle \pi_{t- [s+1,\eps]} ,  G^\infty \Phi_s \rangle \, ds
\, \mathop{\longrightarrow}_{n \to\infty}\, - \int_0^{t} \langle \pi_{t-s} ,  G^\infty \Phi_{s} \rangle \, ds.
$$
On the other hand, we have 
\begin{eqnarray*}
\TT_{2,k} &=& \int_0^\eps {d \over d\tau} \langle \pi_{t_{k}+\tau}, \Phi_{s_{k}} \rangle \, d\tau
\\
&=& \int_0^\eps \langle \pi_{t_{k}+\tau}, G^\infty \Phi_{s_{k}} \rangle \, d\tau
\\
&=& \int_{t_k}^{t_{k+1}} \langle \pi_{\tau}, G^\infty \Phi_{t - [\tau,\eps]}  \rangle \, d\tau.
\end{eqnarray*}
Passing to the limit  $n\to\infty$, we get
$$
\TT_2 =  \int_0^{t} \langle \pi_{\tau}, G^\infty \Phi_{t - [\tau,\eps]}  \rangle \, d\tau
\, \mathop{\longrightarrow}_{n \to\infty}\,  \int_0^{t} \langle \pi_\tau ,  G^\infty \Phi_{t - \tau} \rangle \, d\tau.
$$
As a conclusion, for any $\Phi \in C^1(P_{\GG_1};\R)$, we have proved
$$
\langle \pi_t,\Phi \rangle = \langle \bar\pi_t,\Phi \rangle.
$$
From a density argument we conclude that $ \pi_t = \bar\pi_t$. 
\qed

\smallskip Gathering Lemma~\ref{lem:BBGKY} and
Theorem~\ref{theo:BBGKYuniq} we obtain a propagation to the chaos
result.

\begin{cor}[Abstract chaos propagation]
  Assume {\bf (A2)}, {\bf (A4)}, {\bf (A6)}, {\bf (A7)} and the
  following compatibility  between $G^\infty$ and the sequence
  $G^\infty_\ell$:
  $$
  \forall \, \ell \in \N^*, \,\, \forall \, \varphi \in \FF_1^{\otimes \ell} 
  \qquad
  \langle G^N , R_\varphi \rangle \, \mathop{\longrightarrow}_{N\to\infty} \, 
  \langle G^\infty , R_\varphi \rangle. 
  $$ 
  Assume furthermore that $f^N_0$ is
  $f_0$-chaotic. Then $f^N_t$ is $S^{N\!L}_t f_0$-chaotic.
\end{cor}


\bibliographystyle{acm}
\bibliography{./meanfield}


\signsm \signcm 

\end{document}